\documentclass[aop,MSNbibl,secthm,nochecklpage,dvips]{arximspdf}

\doi{10.1214/08-AOP399}
\volume{37}
\issue{1}
\pubyear{2009}
\firstpage{143}
\lastpage{188}

\makeatletter
\newproclaim{asm}{Assumption}
\newproclaim{note}{Note}
\newtheorem{lem}[thm]{Lemma}
\newtheorem{cor}[thm]{Corollary}

\def\eqref#1{(\ref{#1})}
\makeatother

\begin{document}
\begin{frontmatter}

\title{Quenched limits for transient, zero speed
one-dimensional random walk in random~environment}
\runtitle{Nonexistence of quenched stable limits for RWRE}

\begin{aug}
\author[A]{\fnms{Jonathon} \snm{Peterson}\thanksref{t1,t2}\corref{}\ead[label=e1]{peterson@math.wisc.edu}} \and
\author[B]{\fnms{Ofer} \snm{Zeitouni}\thanksref{t2}\ead[label=e2]{zeitouni@math.umn.edu}}
\runauthor{J. Peterson and O. Zeitouni}
\thankstext{t1}{Supported in part by a Doctoral Dissertation
Fellowship from the University of Minnesota.}
\thankstext{t2}{Supported in part by NSF Grant DMS-05-03775.}
\affiliation{University of Wisconsin and University of Minnesota}
\address[A]{Department of Mathematics\\
University of Wisconsin\\
480 Lincoln Drive\\
Madison, Wisconsin 53705\\
USA\\
\printead{e1}} 
\address[B]{School of Mathematics\\
University of Minnesota\\
206 Church St. SE\\
Minneapolis, Minnesota 55455\\
and\\
Faculty of Mathematics\\
Weizmann Institute of Science\\
Rehovot 76100\\
Israel\\
\printead{e2}}
\end{aug}

\received{\smonth{6} \syear{2006}}
\revised{\smonth{2} \syear{2008}}

%
\begin{abstract}
We consider a nearest-neighbor,
one dimensional random walk $\{X_n\}_{n\geq0}$ in a random
i.i.d. environment,
in the regime where the walk is transient but with zero speed, so that
$X_n$ is of order $n^{s}$ for some $s<1$.
Under the quenched law (i.e., conditioned on the environment), we show that
no limit laws are possible: There exist sequences $\{n_k\}$ and
$\{x_k\}$ depending on the environment only, such that
$X_{n_k}-x_k=o(\log n_k)^2$
(a {\it localized regime}).
On the other hand, there exist sequences $\{t_m\}$ and $\{s_m\}$
depending on
the environment only, such that $\log s_m/\log t_m\to s<1$ and
$P_\omega(X_{t_m}/s_m\leq x)\to1/2$ for all $x>0$
and~\mbox{$\to0$} for $x\leq0$ (a \textit{spread out regime}).
\end{abstract}

%

\end{frontmatter}

\section{Introduction and statement of main results}\label{s1}

Let $\Omega= [0,1]^\mathbb{Z}$, and let $\mathcal{F}$ be the Borel
$\sigma$-algebra on $\Omega$. A random environment is an $\Omega$-valued
random variable $\omega= \{\omega_i\}_{i\in\mathbb{Z}}$ with
distribution $P$. In
this paper we will assume that the $\omega_i$ are i.i.d.

The \textit{quenched} law $P_\omega^x$ for a random walk $X_n$ in the
environment $\omega$ is defined by
\[
P_\omega^x( X_0 = x ) = 1 \quad\mbox{and} \quad
P_\omega^x ( X_{n+1} = j | X_n = i  ) = \cases{%
\omega_i, &\quad  if $j=i+1$, \cr
1-\omega_i, &\quad  if $j=i-1$.}
\]
$\mathbb{Z}^\mathbb{N}$ is the space for the paths of the random walk
$\{X_n\}_{n\in \mathbb{N}}$,
and $\mathcal{G}$ denotes the $\sigma$-algebra generated by the
cylinder sets.
Note that for each $\omega\in\Omega$, $P_\omega$ is a probability measure
on $\mathcal{G}$, and for each $G\in\mathcal{G}$,
$P_\omega^x(G)\dvtx (\Omega, \mathcal{F}) \rightarrow[0,1]$ is a measurable
function of $\omega$. Expectations under the law $P_\omega^x$ are denoted
$E_\omega^x$.

The \emph{annealed} law for the random walk in random
environment $X_n$ is defined by
\[
\mathbb{P}^x(F\times G) = \int_F P_\omega^x(G)P(d\omega),
\qquad F\in\mathcal{F}, G\in\mathcal{G} .
\]
For ease of notation, we will use $P_\omega$ and $\mathbb{P}$ in place
of $P_\omega^0$ and $\mathbb{P}^0$, respectively. We will also use
$\mathbb{P}^x$ to
refer to the marginal on the space of paths, that is, $\mathbb{P}^x(G)=
\mathbb{P}^x(\Omega\times G) = E_P [ P^x_\omega(G)  ]$ for
$G\in\mathcal{G}$. Expectations under the law $\mathbb{P}$ will be
written~$\mathbb{E}$.

A simple criterion for recurrence and a formula for the speed of
transience was given by Solomon in \cite{sRWRE}. For any integers
$i\leq j$, define
%
\begin{equation}\label{rhodef}
\rho_i := \frac{1-\omega_i}{\omega_i}
\quad\mbox{and}\quad
\Pi_{i,j} := \prod_{k=i}^j \rho_k
\end{equation}
and for $x\in\mathbb{Z}$, define the hitting times
\[
T_x:= \min\{n \geq0\dvtx X_n=x\}.
\]
Then $X_n$ is transient to the right (resp. to the left)
if $E_P(\log\rho_0) < 0$ (resp. $E_P \log\rho_0 > 0$) and recurrent
if $E_P (\log\rho_0) = 0$ (henceforth, we will write $\rho$ instead of
$\rho_0$ in expectations involving only $\rho_0$). In the case where
$E_P \log\rho< 0$ (transience to the right),
Solomon established the following law of large numbers:
\[
v_P:= \lim_{n\rightarrow\infty} \frac{X_n}{n} =
\lim_{n\rightarrow\infty} \frac{n}{T_n}
= \frac{1}{\mathbb{E}T_1}, \qquad\mathbb{P}\mbox{-a.s.}
\]
For any integers $i<j$, define
%
\begin{equation}\label{Wdef}
W_{i,j} := \sum_{k=i}^j \Pi_{k,j}
\quad\mbox{and}
\quad W_j := \sum_{k\leq j} \Pi_{k,j}.
\end{equation}
When $E_P \log\rho< 0$, it was shown in \cite{sRWRE} and
\cite{zRWRE} (remark following Lem\-ma~2.1.12), that
%
\begin{equation}\label{QET}
E_\omega^j T_{j+1} = 1+2W_j < \infty, \qquad P\mbox{-a.s.},
\end{equation}
and thus $v_P = 1/(1+2E_P W_0)$. Since $P$ is a product measure,
$E_P W_0 =\sum_{k=1}^\infty (E_P \rho )^k$.
In particular, $v_P =0$ if $E_P \rho\geq1$.

Kesten, Kozlov and Spitzer \cite{kksStable} determined the
annealed limiting distribution of a RWRE with $E_P \log\rho< 0$, that is,
transient to the right. They derived the limiting distributions for the
walk by first establishing a stable
limit law of index $s$ for $T_n$, where $s$ is defined by the equation
\[
E_P\rho^s = 1.
\]
In particular, they showed that when $s<1$, there exists a $b>0$ such that
\[
\lim_{n\rightarrow\infty} \mathbb{P}
\biggl( \frac{T_n}{n^{1/s}} \leq x  \biggr) = L_{s,b}(x),
\]
and
%
\begin{equation}\label{annealedstable}
\lim_{n\rightarrow\infty} \mathbb{P} \biggl( \frac{X_n}{n^s} \leq x
 \biggr) = 1-L_{s,b}(x^{-1/s}),
\end{equation}
where $L_{s,b}$ is the distribution function for a stable random
variable with characteristic function
%
\begin{equation}\label{char}
\hat{L}_{s,b}(t)= \exp \biggl\{ -b|t|^s  \biggl(
1-i\frac{t}{|t|}\tan(\pi s/2)  \biggr) \biggr\}.
\end{equation}
The value of $b$ was recently identified \cite{ESZ}.
While the annealed limiting distributions for transient
one-dimensional RWRE have been known for quite a while, the corresponding
quenched limiting distributions have remained largely unstudied
until recently. Goldsheid \cite{gQCLT} and Peterson \cite{pThesis}
independently proved that when $s>2$, a quenched CLT holds with a
random (depending on the environment)
centering. A~similar result was given by Rassoul-Agha and Sepp\"al\"ainen
in \cite{rsBFD} under different assumptions on the environment.
Previously, in \cite{kmCLT} and \cite{zRWRE},
it was shown that the limiting statement for
the quenched CLT with random centering holds in probability rather
than almost surely. No other results of quenched limiting
distributions are known
when $s\leq2$.

In this paper, we analyze the quenched limiting distributions of a
one-dimensional transient RWRE in the case $s<1$. One could expect
that the quenched limiting distributions are of the same type as
the annealed limiting distributions since annealed probabilities
are averages of quenched probabilities. However, this turns out not to
be the
case. In fact, a consequence of our main results,
Theorems \ref{refstable}, \ref{local} and \ref{nonlocal} below is that
the annealed stable behavior of $T_n$ comes from
fluctuations in the environment.

Throughout the paper, we will
make the following assumptions.
\begin{asm} \label{essentialasm}
$P$ is an i.i.d. product measure on $\Omega$ such that
%
\begin{equation}\label{zerospeedregime}
E_P \log\rho< 0 \quad\mbox{and}\quad
E_P \rho^s = 1 \qquad \mbox{for some } s>0 .
\end{equation}
\end{asm}
\begin{asm}
The distribution of $\log\rho$ is nonlattice under
$P$ and $E_P \rho^s 
\log\rho<\infty$. \label{techasm}
\end{asm}
\begin{note*}
Since $E_P \rho^\gamma$ is a convex function of
$\gamma$, the two statements in \eqref{zerospeedregime} imply that
$E_P \rho^\gamma< 1$ for all
$\gamma<s$ and $E_P \rho^\gamma> 1$ for all $\gamma> s$.
Assumption \ref{essentialasm}
contains the essential assumption necessary for the walk to be
transient. The main results of this paper are for $s<1$ (the zero-speed
regime), but many statements hold for $s\in(0,2)$ or even $s\in
(0,\infty)$. If no mention is made of bounds on $s$, then it is
assumed that the statement holds for all $s>0$. We recall that the
technical conditions contained in
Assumption \ref{techasm} were also invoked in \cite{kksStable}.
\end{note*}

Define the ``ladder locations'' $\nu_i$ of the environment by
%
\begin{equation}\label{nudef}
\quad
\nu_0 = 0 \quad\mbox{and}\quad\nu_i = \cases{%
\inf\{n > \nu_{i-1}\dvtx \Pi_{\nu_{i-1},n-1} < 1\}, &\quad $i \geq1$,\vspace*{2pt}\cr
\sup\{j < \nu_{i+1}\dvtx \Pi_{k,j-1}<1,\ \forall k<j \}, & \quad $i\leq-1$.}
\end{equation}
%
Throughout the remainder of the paper, we will let $\nu=\nu_1$. We
will sometimes refer to sections of the environment between
$\nu_{i-1}$ and $\nu_i -1$ as ``blocks'' of the environment. Note
that the block between $\nu_{-1}$ and $\nu_0 -1$ is different from
all the other blocks between consecutive ladder locations. Define
the measure $Q$ on environments by $Q(\cdot):=P(\cdot|\mathcal{R})$,
where the event
\[
\mathcal{R}:=\{ \omega\in\Omega\dvtx
\Pi_{-k,-1} < 1,\ \forall k \geq 1\}.
\]
Note that $P(\mathcal{R}) > 0$ since $E_P \log\rho< 0$.
$Q$ is defined so that the blocks of the environment between ladder
locations are i.i.d.
under $Q$, all with distribution the same as that of the block
from $0$ to $\nu-1$ under $P$. In Section \ref{stablecrossing}, we
prove the following annealed theorem.
\begin{thm}\label{refstable}
Let Assumptions \ref{essentialasm} and \ref{techasm} hold, and let
$s<1$. Then there exists a $b'>0$ such that
\[
\lim_{n\rightarrow\infty} Q \biggl( \frac{ E_\omega T_{\nu_n}}
{n^{1/s}} \leq x \biggr) = L_{s,b'}(x).
\]
\end{thm}

We then use Theorem \ref{refstable} to prove the following two theorems
which show that $P$-a.s. there exist two different random sequences of
times (depending on the environment) where the random walk has
different limiting behavior. These are the main results of the paper.
\begin{thm}\label{local}
Let Assumptions \ref{essentialasm} and \ref{techasm} hold, and let
$s<1$. Then
$P$-a.s. there exist random subsequences $t_m=t_m(\omega)$ and
$u_m=u_m(\omega)$, such that for any $\delta > 0$,
\[
\lim_{m\rightarrow\infty} P_\omega
\biggl( \frac{X_{t_m} - u_m}{(\log t_m)^2} \in[-\delta , \delta ]  \biggr) = 1.
\]
\end{thm}
\begin{thm} \label{nonlocal}
Let Assumptions \ref{essentialasm} and \ref{techasm} hold, and let
$s<1$. Then
$P$-a.s. there exists a random subsequence $n_{k_m}=n_{k_m}(\omega)$
of $n_k=2^{2^k}$ and a random sequence $t_m=t_m(\omega)$, such that
\[
\lim_{m\rightarrow\infty} \frac{\log t_m}{\log n_{k_m} }= \frac{1}{s},
\]
and
\[
\lim_{m\rightarrow\infty}
P_\omega \biggl(\frac{X_{t_m}}{n_{k_m}} \leq x  \biggr)
= \cases{%
0, &\quad if $x \leq0$,\cr
\dfrac{1}{2}, &\quad if $0<x<\infty$.}
\]
\end{thm}

Note that Theorems \ref{local} and \ref{nonlocal}
preclude the possibility of a quenched analogue of the annealed
statement \eqref{annealedstable}.
It should be noted that in \cite{gsMVSS}, Gantert and Shi prove
that when $s\leq1$, there exists a random sequence of times $t_m$ at which
the local time of the random walk at a single site
is a positive fraction of $t_m$. This is related to the
statement of Theorem \ref{local}, but we do not see a simple
argument which directly implies Theorem \ref{local} from the
results of \cite{gsMVSS}.

As in \cite{kksStable}, limiting distributions for $X_n$ arise from first
studying limiting distributions for $T_n$. Thus,
to prove Theorem \ref{nonlocal}, we first prove that there
exists random subsequences $x_m=x_m(\omega)$ and $v_{m,\omega}$ in which
\[
\lim_{m\rightarrow\infty} P_\omega
\biggl( \frac{T_{x_m} - E_\omega T_{x_m}}{\sqrt{v_{m,\omega}}}
\leq y  \biggr) = \int_{-\infty}^y
\frac{1}{\sqrt{2\pi}} e^{-t^2/2} \,dt =: \Phi(y).
\]
We actually prove a stronger statement than this in Theorem \ref{gaussianT}
below,
where we prove that all $x_m$ ``near'' a subsequence $n_{k_m}$
of $n_k=2^{2^k}$ have the same Gaussian behavior (what we mean by ``near''
the subsequence $n_{k_m}$ is made precise in the statement of the theorem).

The structure of the paper is as follows.
In Section \ref{introlemmas}, we prove some introductory lemmas
which will be used throughout the paper. Section \ref{stablecrossing} is
devoted to proving Theorem \ref{refstable}. In Section \ref{localization},
we use the latter to prove Theorem \ref{local}. In Section \ref{gaussian},
we prove the existence of random subsequences $\{n_k\}$
where $T_{n_k}$ is approximately
Gaussian, and use this fact to prove Theorem \ref{nonlocal}.
Section \ref{tailofTnu} contains the proof of the following technical
theorem which is used throughout the paper.
\begin{thm}\label{Tnutail}
Let Assumptions \ref{essentialasm} and \ref{techasm} hold. Then there
exists a constant $K_\infty\in(0,\infty)$ such that
\[
Q(E_\omega T_\nu> x ) \sim K_\infty x^{-s}.
\]
\end{thm}

The proof of Theorem \ref{Tnutail} is based on
results from \cite{kRDE} and mimics the proof of tail asymptotics
in \cite{kksStable}.

\section{Introductory lemmas}\label{introlemmas}

Before proceeding with the proofs of the main theorems, we mention a
few easy lemmas which will be used throughout the rest of the paper.
Recall the definitions of $\Pi_{1,k}$ and $W_i$ in \eqref{rhodef} and
\eqref{Wdef}.
\begin{lem}\label{nutail}
For any $c< -E_P \log\rho$, there exist $\delta _c, A_c>0$ such that
%
\begin{equation}\label{LDPrho}
P(\Pi_{1,k} > e^{-c k }) = P  \Biggl( \frac{1}{k}\sum_{i=1}^k \log
\rho_i > -c  \Biggr) \leq A_ce^{-\delta _c k}.
\end{equation}
Also, there exist constant $C_1,C_2>0$ such that $P(\nu> x) \leq C_1
e^{-C_2 x}$ for all $x\geq0$.
\end{lem}
\begin{pf}
First, note that due to Assumption \ref{essentialasm}, $\log\rho$
has negative mean and finite exponential moments in a
neighborhood of zero. If $c < - E_P \log\rho$,
Cram\'{e}r's theorem (\cite{dzLDTA}, Theorem 2.2.3) then yields \eqref{LDPrho}.
By the definition of $\nu$, we have $P(\nu>x) \leq P(\Pi_{0,\lfloor
x \rfloor-1} \geq1)$, which together with
\eqref{LDPrho}, completes the proof of the lemma.
\end{pf}

From \cite{kRDE}, Theorem 5, there exist constants
$K,K_1>0$ such that for all $i$
%
\begin{equation}\label{PWtail}
P(W_i > x) \sim K x^{-s}
\quad\mbox{and}\quad
P(W_i > x) \leq K_1 x^{-s}.
\end{equation}
The tails of $W_{-1}$, however, are different (under the measure $Q$),
as the following lemma shows.
\begin{lem}\label{Wtail}
There exist constants $C_3,C_4>0$ such that
$Q(W_{-1} > x) \leq C_3 e^{-C_4 x}$ for all $x\geq0$.
\end{lem}
\begin{pf}
Since $\Pi_{i,-1} < 1$, $Q$-a.s. we have $W_{-1} < k + \sum_{i
< -k} \Pi_{i,-1}$ for any $k>0$. Also, note that from
\eqref{LDPrho}, we have $Q(\Pi_{-k,-1} > e^{-c k} ) \leq
{A_c} e^{-\delta _c k}/P(\mathcal{R})$. Thus,
\begin{eqnarray*}
Q(W_{-1} > x) &\leq & Q \Biggl( \frac{x}{2} +
\sum_{k={x}/{2}}^\infty e^{-ck} > x  \Biggr)
+ Q \biggl( \Pi_{-k,-1}> e^{-ck} \mbox{, for some } k\geq\frac{x}{2} \biggr)
\\
& \leq& \mathbf{1}_{{x}/{2} + {1}/{(1-e^{-c})} > x}
+ \sum_{k={x}/{2}}^\infty Q(\Pi_{-k,-1}>e^{-ck})
\\
&\leq& \mathbf{1}_{{1}/{(1-e^{-c})} > {x}/{2}}
+ \mathcal{O} (e^{-\delta _c x/2} ).
\end{eqnarray*}\upqed
\end{pf}

We also need a few more definitions
that will be used throughout the paper. For any $i\leq k$,
%
\begin{equation}\label{Rdef}
R_{i,k} := \sum_{j=i}^k \Pi_{i,j} \quad\mbox{and}\quad
R_i:= \sum_{j=i}^\infty\Pi_{i,j}.
\end{equation}
Note that since $P$ is a product measure,
$R_{i,k}$ and $R_i$ have the same distributions as $W_{i,k}$ and
$W_i$ respectively. In particular with $K,K_1$, the same as
in \eqref{PWtail},
%
\begin{equation}\label{Rtail}
P(R_i > x) \sim K x^{-s} \quad\mbox{and}\quad
P(R_i > x) \leq K_1 x^{-s}.
\end{equation}


\section{Stable behavior of expected crossing time}\label{stablecrossing}

Recall from Theorem \ref{Tnutail} that there exists
$K_\infty>0$ such that $Q(E_\omega T_\nu> x) \sim K_\infty x^{-s}$.
Thus, $E_\omega T_\nu$ is in the domain of attraction of a
stable distribution. Also, from the comments after the definition of
$Q$ in the Introduction, it is evident that under $Q$, the environment
$\omega$ is stationary under shifts of the ladder times
$\nu_i$. Thus, under $Q$,
$\{ E_\omega^{\nu_{i-1}} T_{\nu_i} \}_{i\in\mathbb{Z}}$
is a stationary sequence of random variables. Therefore, it is reasonable
to expect that $n^{-1/s}E_\omega T_{\nu_n} = n^{-1/s}
\sum_{i=1}^n E_\omega^{\nu_{i-1}} T_{\nu_i}$ converge in distribution
to a stable distribution of index $s$. The main obstacle to proving
this is that the random variables $E_\omega^{\nu_{i-1}} T_{\nu_i}$ are
not independent. This dependence, however, is rather weak. The strategy
of the proof of Theorem \ref{refstable} is to first show that we need
only consider the blocks where the expected crossing time
$E_\omega^{\nu_{i-1}} T_{\nu_i}$ is relatively large. These blocks will
then be separated enough to make the expected crossing times
essentially independent.

For every $k\in\mathbb{Z}$, define
%
\begin{equation}\label{Mdef}
M_k:=\max \{\Pi_{\nu_{k-1}, j} \dvtx \nu_{k-1}\leq j < \nu_k \}.
\end{equation}
Theorem 1 in
\cite{iEV} gives that there exists a constant $C_5>0$ such that
%
\begin{equation}\label{Mtail}
Q(M_1 > x)\sim C_5 x^{-s}.
\end{equation}
Thus, $M_1$ and $E_\omega T_\nu$ have similar tails under $Q$. We will
now show
that $E_\omega T_\nu$ cannot be too much larger than $M_1$.
From \eqref{QET}, we have that
%
\begin{equation}\label{ETnuexpand}
E_\omega T_\nu= \nu+ 2 \sum_{j=0}^{\nu-1} W_j = \nu+2
W_{-1}R_{0,\nu-1} + 2 \sum_{i=0}^{\nu-1} R_{i,\nu-1}.
\end{equation}
From the definitions of $\nu$ and $M_1$, we have that $R_{i,\nu-1}
\leq(\nu- i) M_1 \leq\nu M_1$ for any $0\leq i < \nu$.
Therefore, $E_\omega T_\nu\leq\nu+ 2 W_{-1}\nu M_1 + 2 \nu^2 M_1$.
Thus, given any $0<\alpha <\beta $ and $\delta >0$, we have
%
\begin{eqnarray}\label{TbigMsmall}
\qquad
Q(E_\omega T_\nu> \delta n^{\beta }, M_1 \leq n^{\alpha })
&\leq & Q(\nu+2 W_{-1}\nu n^{\alpha } + 2\nu^2 n^{\alpha } > \delta
n^{\beta })
\nonumber\\
&\leq & Q\bigl(W_{-1} > n^{(\beta -\alpha )/2}\bigr)
+ Q \bigl(\nu^2 > n^{(\beta -\alpha)/2} \bigr)
\\
&=& o \bigl(e^{-n^{(\beta -\alpha )/5}} \bigr),
\nonumber
\end{eqnarray}
where the second inequality holds for all $n$ large enough and the
last equality is a result of Lemmas \ref{nutail} and \ref{Wtail}.
We now show that only the ladder times with $M_k>n^{(1-\varepsilon)/s}$
contribute to the limiting distribution of $n^{-1/s} E_\omega
T_{\nu_n}$.
\begin{lem}\label{smallblocks}
Assume $s<1$. Then for any $\varepsilon>0$ and any $\delta >0$, there
exists an $\eta> 0$ such that
\[
\lim_{n\rightarrow\infty} Q \Biggl( \sum_{i=1}^n (E_\omega^{\nu
_{i-1}} T_{\nu_i}) \mathbf{1}_{M_i \leq n^{(1-\varepsilon)/s}}
> \delta n^{1/s}  \Biggr) = o(n^{-\eta}).
\]
\end{lem}
\begin{pf}
First note that
\begin{eqnarray*}
&& Q  \Biggl( \sum_{i=1}^n (E_\omega^{\nu_{i-1}} T_{\nu_i} ) \mathbf{1}_{
M_i \leq n^{(1-\varepsilon)/s}} > \delta n^{1/s}  \Biggr)
\\
&&\qquad \leq Q \Biggl(
\sum_{i=1}^n (E_\omega^{\nu_{i-1}} T_{\nu_i} )
\mathbf{1}_{E_\omega^{\nu_{i-1}} T_{\nu_i}
\leq n^{(1-{\varepsilon}/{2} )/s}} > \delta n^{1/s}  \Biggr)
\\
&&\qquad \quad {} + n Q \bigl( E_\omega T_{\nu}
> n^{(1-{\varepsilon}/{2})/s}, M_1 \leq n^{(1-\varepsilon)/s} \bigr).
\end{eqnarray*}
By \eqref{TbigMsmall}, the last term above decreases
faster than any power of $n$. Thus, it is enough to prove that for
any $\delta ,\varepsilon>0$, there exists an $\eta>0$ such that
\[
Q \Biggl( \sum_{i=1}^n (E_\omega^{\nu_{i-1}} T_{\nu_i} )
\mathbf{1}_{E_\omega^{\nu_{i-1}} T_{\nu_i}
\leq n^{(1-\varepsilon)/s}} > \delta n^{1/s}  \Biggr) = o( n^{-\eta}
).
\]
Next, pick $C\in (1, \frac{1}{s}  )$ and let
$J_{C,\varepsilon,k,n}:=  \{ i\leq n \dvtx n^{(1-C^k \varepsilon)/s}
< E_\omega^{\nu_{i-1}} T_{\nu_i} \leq\break n^{(1-C^{k-1}\varepsilon)/s}\}$. Let
$k_0=k_0(C,\varepsilon)$ be the smallest integer such that
$(1-C^k\varepsilon) \leq 0$. Then for any $k < k_0$, we have
\begin{eqnarray*}
Q  \Biggl( \sum_{i\in J_{C,\varepsilon,k,n} } E_\omega^{\nu_{i-1}}
T_{\nu_i} > \delta n^{1/s}  \Biggr)
&\leq& Q  \bigl( \# J_{C,\varepsilon,k,n} > \delta
n^{1/s-(1-C^{k-1}\varepsilon)/s}  \bigr)
\\
&\leq& \frac{n Q( E_\omega T_{\nu} > n^{(1-C^{k}\varepsilon
)/s})}{\delta n^{C^{k-1}\varepsilon/s}} \sim\frac{K_\infty}{\delta }
n^{-C^{k-1}\varepsilon({1}/{s}-C)},
\end{eqnarray*}
where the asymptotics in the last line above is from Theorem
\ref{Tnutail}. Letting $\eta=
\frac{\varepsilon}{2} (\frac{1}{s}-C )$, we have for any
$k < k_0$
that
%
\begin{equation}\label{klessk0}
Q  \Biggl( \sum_{i\in J_{C,\varepsilon,k,n} } E_\omega^{\nu_{i-1}}
T_{\nu_i} > \delta n^{1/s}  \Biggr) = o(n^{-\eta}).
\end{equation}
Finally, note that
%
\begin{equation}\label{k0}
\qquad \quad
Q \Biggl( \sum_{i=1}^n (E_\omega^{\nu_{i-1}} T_{\nu_i})
\mathbf{1}_{E_\omega^{\nu_{i-1}} T_{\nu_i} \leq
n^{(1-C^{k_0-1}\varepsilon)/s}}
\geq\delta n^{1/s}  \Biggr)
\leq \mathbf{1}_{n^{1+(1-C^{k_0-1}\varepsilon)/s} \geq\delta n^{1/s}}.
\end{equation}
However, since $C^{k_0} \varepsilon\geq1 > Cs$, we have
$C^{k_0-1}\varepsilon > s $, which implies that
the right side of \eqref{k0} vanishes for all $n$ large enough.
Therefore, combining \eqref{klessk0} and~\eqref{k0}, we have
\begin{eqnarray*}
&& Q \Biggl( \sum_{i=1}^n (E_\omega^{\nu_{i-1}} T_{\nu_i})
\mathbf{1}_{E_\omega^{\nu_{i-1}} T_{\nu_i}
\leq n^{(1-\varepsilon)/s}} > \delta n^{1/s}  \Biggr)
\\
&&\qquad \leq \sum_{k=1}^{k_0-1} Q
\Biggl( \sum_{i\in J_{C,\varepsilon,k,n} }
E_\omega^{\nu_{i-1}} T_{\nu_i} > \frac{\delta }{k_0}
n^{1/s}  \Biggr)
\\
&&\quad \qquad {} + Q \Biggl( \sum_{i=1}^n (E_\omega^{\nu_{i-1}} T_{\nu_i})
\mathbf{1}_{E_\omega^{\nu_{i-1}} T_{\nu_i} \leq
n^{(1-C^{k_0-1}\varepsilon)/s}}
\geq\frac{\delta }{k_0}
n^{1/s}  \Biggr) = o(n^{-\eta}).
\end{eqnarray*}\upqed
\end{pf}

In order to make the crossing times of the significant blocks essentially
independent, we introduce some reflections to the RWRE. For
$n=1,2,\ldots,$ define
%
\begin{equation}\label{bdef}
b_n:= \lfloor\log^2(n) \rfloor.
\end{equation}
Let $\bar{X}_t^{(n)}$ be the random walk that is the same as $X_t$
with the added condition that after reaching $\nu_k$ the
environment is modified by setting $\omega_{\nu_{k-b_n}} = 1 $, that is,
never allow the walk to backtrack more than $\log^2(n)$ ladder
times.
We couple $\bar{X}_t^{(n)}$ with the random walk $X_t$ in such a way
that $\bar{X}_t^{(n)} \geq X_t$ with equality holding until
the first time $t$ when the walk $\bar{X}_t^{(n)}$ reaches a modified
environment location.
%
Denote by $\bar{T}_{x}^{(n)}$ the corresponding hitting
times for the walk $\bar{X}_t^{(n)}$. The following lemmas show that
we can add reflections to
the random walk without changing the expected crossing time by
very much.
\begin{lem} \label{ETdiff}
There exist $B,\delta ' > 0$ such that for any $x>0$
\[
Q \bigl( E_\omega T_\nu- E_\omega\bar{T}_\nu^{(n)} > x  \bigr)
\leq B (x^{-s}\vee1) e^{-\delta ' b_n}.
\]
\end{lem}
\begin{pf}
First, note that for any $n$ the formula for $E_\omega
\bar{T}_\nu^{(n)}$ is the same as for $E_\omega T_\nu$ in
\eqref{ETnuexpand} except with $\rho_{\nu_{-b_n}} = 0$. Thus,
$E_\omega T_\nu$ can be written as
%
\begin{equation}\label{reflectexpand}
E_\omega T_\nu= E_\omega\bar{T}_\nu^{(n)}
+ 2 (1+W_{\nu_{-b_n}-1}) \Pi_{\nu_{-b_n } , -1}R_{0,\nu-1}.
\end{equation}
Now, since $\nu_{-b_n} \leq-b_n$, we have
\begin{eqnarray*}
Q ( \Pi_{\nu_{-b_n},-1}
> e^{-c b_n}  ) &\leq& \sum_{k=b_n}^\infty Q ( \Pi_{-k,-1}
> e^{-c k}  )
\\
& \leq& \sum_{k=b_n}^\infty
\frac{1}{P(\mathcal{R})} P ( \Pi_{-k,-1} > e^{-c k}  ).
\end{eqnarray*}
Applying \eqref{LDPrho},
we have that for any $0<c<-E_P \log\rho$,
there exist $A',\delta _c > 0$ such that
\[
Q ( \Pi_{\nu_{-b_n},-1}
> e^{-c b_n}  ) \leq A' e^{-\delta _c b_n} .
\]
Therefore, for any
$x>0$,
%
\begin{eqnarray}\label{ETd1}
\qquad
Q \bigl( E_\omega T_{\nu} - E_\omega\bar{T}_{\nu}^{(n)} > x
\bigr) &\leq & Q \bigl( 2(1+W_{\nu_{-b_n}-1})
\Pi_{\nu_{-b_n},-1} R_{0,\nu-1} > x  \bigr)
\nonumber\\
& \leq & Q \bigl( 2(1+W_{\nu_{-b_n}-1}) R_{0,\nu-1} >
x e^{c b_n}  \bigr) + A' e^{-\delta _c b_n}
\\
&=& Q \bigl( 2(1+W_{-1}) R_{0,\nu-1} > x e^{c b_n}  \bigr)
+ A' e^{-\delta _c b_n},
\nonumber
\end{eqnarray}
where the equality in the second line is due to the fact that the
blocks of the environment are i.i.d. under $Q$. Also, from
\eqref{ETnuexpand} and Theorem \ref{Tnutail}, we have
%
\begin{equation}\label{ETd2}
\qquad
Q \bigl( 2(1+W_{-1}) R_{0,\nu-1} > x e^{c b_n}  \bigr)
\leq Q ( E_\omega T_\nu> x e^{c b_n}  ) \sim K_\infty x^{-s} e^{-c s b_n}.
\end{equation}
Combining \eqref{ETd1} and \eqref{ETd2} completes the proof.
\end{pf}
\begin{lem}\label{reftail}
For any $x>0$ and $\varepsilon>0$, we have that
%
\begin{equation}\label{refTnutail}
\lim_{n\rightarrow\infty} n Q
\bigl(E_\omega\bar{T}_\nu^{(n)} > x n^{1/s},
M_1 > n^{(1-\varepsilon)/s}  \bigr) = K_\infty x^{-s}.
\end{equation}
\end{lem}
\begin{pf}
Since adding reflections only decreases the crossing times, we can
get an upper bound using Theorem \ref{Tnutail},
that is,
%
\begin{eqnarray}\label{refTnub}
&& \limsup_{n\rightarrow\infty} n Q \bigl(E_\omega\bar{T}_\nu^{(n)}
> x n^{1/s}, M_1 > n^{(1-\varepsilon)/s}  \bigr)
\nonumber\\[-8pt]
\\[-8pt]
&&\qquad \leq\limsup_{n\rightarrow\infty} n Q(E_\omega T_\nu> x n^{1/s}) =
K_\infty x^{-s}.
\nonumber
\end{eqnarray}
To get a lower bound, we first note that for any $\delta >0$,
%
\begin{eqnarray}\label{r1}
\qquad
Q \bigl( E_\omega T_\nu> (1+\delta )x n^{1/s}  \bigr)
&\leq& Q \bigl(E_\omega \bar{T}_\nu^{(n)} > x n^{1/s},
M_1 > n^{(1-\varepsilon)/s}  \bigr)
\nonumber\\
&&{} + Q \bigl( E_\omega T_\nu- E_\omega\bar{T}_\nu^{(n)} > \delta x
n^{1/s}  \bigr)
\nonumber\\[-8pt]
\\[-8pt]
&& {} + Q \bigl( E_\omega T_\nu> (1+\delta )x n^{1/s}, M_1 \leq
n^{(1-\varepsilon)/s}  \bigr)
\nonumber\\
&\leq & Q \bigl( E_\omega\bar{T}_\nu^{(n)} > x n^{1/s}, M_1 >
n^{(1-\varepsilon)/s}  \bigr) + o(1/n),\hspace*{-6pt}
\nonumber
\end{eqnarray}
where the second inequality is from \eqref{TbigMsmall} and
Lemma~\ref{ETdiff}. Again, using Theorem~\ref{Tnutail}, we have
%
\begin{eqnarray}\label{refTnulb}
&& \liminf_{n\rightarrow\infty} n Q
\bigl(E_\omega\bar{T}_\nu^{(n)} > x n^{1/s},
M_1 > n^{(1-\varepsilon)/s}  \bigr)
\nonumber\\
&&\qquad \geq \liminf_{n\rightarrow\infty} n
Q \bigl( E_\omega T_\nu> (1+\delta )x n^{1/s}  \bigr) - o(1)
\\
&&\qquad = K_\infty(1+\delta )^{-s}x^{-s}.
\nonumber
\end{eqnarray}
Thus, by applying \eqref{refTnub}
and \eqref{refTnulb} and then letting $\delta \rightarrow0$, we get
\eqref{refTnutail}.
\end{pf}

Our general strategy is to show that the partial sums
\[
\frac{1}{n^{1/s}} \sum_{k=1}^n E_\omega^{\nu_{k-1}}
\bar{T}_{\nu_k}^{(n)} \mathbf{1}_{M_k>n^{(1-\varepsilon)/s} }
\]
converge in distribution to a stable law of parameter $s$.
To establish this, we will need bounds on the mixing properties of the sequence
$E_\omega^{\nu_{k-1}}
\bar{T}_{\nu_k}^{(n)}\mathbf{1}_{M_k>n^{(1-\varepsilon)/s} }$. As in
\cite{kGPD}, we say that an array $\{ \xi_{n,k}\dvtx  k\in\mathbb{Z},
n\in \mathbb{N}\}$ which is stationary in rows is $\alpha $-mixing if
$\lim_{k\rightarrow\infty}\limsup_{n\rightarrow\infty} \alpha_n(k) = 0$, where
\begin{eqnarray*}
\alpha _n(k) &:=& \sup \{ |P(A\cap B)-P(A)P(B)|\dvtx
A\in \sigma (\ldots, \xi_{n,-1}, \xi_{n,0}  ),
\\
&&\hspace*{140pt}
B\in\sigma (\xi_{n,k},\xi_{n,k+1}, \ldots )  \}.
\end{eqnarray*}

\begin{lem} \label{alphamixing}
For any $0<\varepsilon< \frac{1}{2}$, under the measure $Q$, the
array of random variables
$\{ E_\omega^{\nu_{k-1}} \bar{T}_{\nu_k}^{(n)}
\mathbf{1}_{M_k>n^{(1-\varepsilon)/s} } \}_{k\in\mathbb{Z}, n\in
\mathbb{N}}$ is $\alpha $-mixing with
\[
\sup_{k\in[1,\log^2 n]}\alpha_n(k)=o(n^{-1+2\epsilon}),
\qquad\alpha_n(k)=0\qquad\forall k>\log^2 n.
\]
\end{lem}
\begin{pf}
Fix $\varepsilon\in(0,\frac{1}{2})$. For ease of notation, define
$\xi_{n,k}:= E_\omega^{\nu_{k-1}} \bar{T}_{\nu_k}^{(n)}\times\break
\mathbf{1}_{M_k>n^{(1-\varepsilon)/s} }$. As we mentioned before,
under $Q$ the environment is stationary
under shifts of the sequence of ladder locations and
thus $\xi_{n,k}$ is stationary in rows under $Q$.

If $k>\log^2(n)$, then because of the reflections,
$\sigma (\ldots, \xi_{n,-1}, \xi_{n,0}  )$ and
$\sigma (\xi_{n,k},\break \xi_{n,k+1}, \ldots )$ are
independent and so $\alpha _n(k)= 0$.
To handle the case when $k\leq\log^2(n)$,
fix $A\in \sigma (\ldots, \xi_{n,-1}, \xi_{n,0}  )$ and
$B\in \sigma (\xi_{n,k},\xi_{n,k+1}, \ldots ) $, and
define the event
\[
C_{n,\varepsilon}:=\bigl\{M_j \leq n^{(1-\varepsilon)/s}, \mbox{ for }
1\leq j \leq b_n \bigr\} = \{\xi_{n,j} = 0, \mbox{ for } 1\leq j \leq b_n \}.
\]
For any $j> b_n$, we have that $\xi_{n,j}$ only depends on the
environment to the right of zero. Thus,
%
\[
Q(A\cap B \cap C_{n,\varepsilon})= Q(A)Q(B\cap C_{n,\varepsilon})
\]
since $B \cap C_{n,\varepsilon} \in\sigma(\omega_0,\omega_1,\ldots)$. Also,
note that by \eqref{Mtail} we have $Q(C_{n,\varepsilon}^c ) \leq b_n Q(M_1 >
n^{(1-\varepsilon)/s}) = o(n^{-1+2\varepsilon})$. Therefore,
\begin{eqnarray*}
|Q(A\cap B)-Q(A)Q(B)| & \leq& |Q(A\cap B) - Q(A\cap B\cap C_{n,\varepsilon})|
\\
&&{} + |Q(A\cap B\cap C_{n,\varepsilon}) - Q(A)Q(B\cap C_{n,\varepsilon})|
\\
&&{} + Q(A)|Q(B\cap C_{n,\varepsilon})-Q(B)|
\\
&\leq& 2 Q(C_{n,\varepsilon}^c) = o(n^{-1+2\varepsilon}).
\end{eqnarray*}\upqed
\end{pf}
\begin{pf*}{Proof of Theorem \ref{refstable}}
First, we show that the partial sums
\[
\frac{1}{n^{1/s}} \sum_{k=1}^n E_\omega^{\nu_{k-1}}
\bar{T}_{\nu_k}^{(n)} \mathbf{1}_{M_k>n^{(1-\varepsilon)/s} }
\]
converge in distribution to a stable random variable of parameter $s$.
To this end, we will
apply \cite{kGPD}, Theorem 5.1(III). We now verify
the conditions of that theorem.
The first condition that needs to be satisfied is
\[
\lim_{n\rightarrow\infty} n Q \bigl( n^{-1/s} E_\omega\bar{T}_{\nu}^{(n)}
\mathbf{1}_{M_1>n^{(1-\varepsilon)/s} } > x  \bigr)
= K_\infty x^{-s}.
\]
However, this is exactly the content of Lemma \ref{reftail}.

Secondly, we need a sequence $m_n$ such that $m_n\rightarrow\infty$,
$m_n=o(n)$ and $n\alpha _n(m_n) \rightarrow0$, and such that for any
$\delta >0$,
%
\begin{eqnarray} \label{mncond}
&& \lim_{n\rightarrow\infty} \sum_{k=1}^{m_n} n Q
\bigl( E_\omega\bar{T}_{\nu}^{(n)}
\mathbf{1}_{M_1>n^{(1-\varepsilon)/s} } > \delta n^{1/s},
\nonumber\\[-8pt]
\\[-8pt]
&&\hspace*{61pt}
E_\omega^{\nu_{k}} \bar{T}_{\nu_{k+1}}^{(n)}
\mathbf{1}_{M_{k+1}>n^{(1-\varepsilon)/s} } > \delta n^{1/s}
\bigr) = 0.
\nonumber
\end{eqnarray}
However, by the independence of $M_1$ and $M_{k+1}$ for any
$k\geq1$, the probability inside the sum is less
than $Q(M_1> n^{(1-\varepsilon)/s})^2$. By \eqref{Mtail}, this last
expression is~$ \sim C_5 n^{-2+2\varepsilon} $.
Thus, letting $m_n = n^{1/2- \varepsilon}$ yields \eqref{mncond}.
[Note that by
Lemma \ref{alphamixing}, $n\alpha _n(m_n) = 0$ for all $n$ large
enough.]

Finally, we need to show that
%
\begin{equation}\label{truncexp}
\lim_{\delta \rightarrow0}\limsup_{n\rightarrow\infty} n E_Q
\bigl [ n^{-1/s} E_\omega\bar{T}_{\nu}^{(n)} \mathbf
{1}_{M_1>n^{(1-\varepsilon
)/s} }
\mathbf{1}_{E_\omega\bar{T}_{\nu}^{(n)} \leq\delta }  \bigr] = 0.
\end{equation}
Now, by \eqref{refTnub}, there exists a constant $C_6>0$ such that
for any $x > 0$,
\[
Q \bigl( E_\omega\bar{T}_{\nu}^{(n)} > x n^{1/s} , M_1
> n^{(1-\varepsilon)/s}  \bigr) \leq C_6 x^{-s}\frac{1}{n}.
\]
Then using this, we have
\begin{eqnarray*}
&& n E_Q \bigl [ n^{-1/s} E_\omega\bar{T}_{\nu}^{(n)} \mathbf
{1}_{M_1>n^{(1-\varepsilon)/s} }
\mathbf{1}_{E_\omega\bar{T}_{\nu}^{(n)} \leq\delta }  \bigr]
\\
&&\qquad = n \int_0^\delta Q \bigl( E_\omega\bar{T}_{\nu}^{(n)} > x n^{1/s} ,
M_1 > n^{(1-\varepsilon)/s}  \bigr) \,dx
\\
&&\qquad \leq  C_6 \int_0^\delta x^{-s} \,dx
= \frac{C_6 \delta ^{1-s}}{1-s},
\end{eqnarray*}
where the last
integral is finite since $s<1$. Equation \eqref{truncexp} follows.

Having checked all its hypotheses, Kobus
(\cite{kGPD}, Theorem 5.1(III)) applies and yields
that there exists a $b'>0$ such that
%
\begin{equation}\label{refdist}
Q \Biggl( \frac{1}{n^{1/s}} \sum_{k=1}^n E_\omega^{\nu_{k-1}}
\bar{T}_{\nu_k}^{(n)} \mathbf{1}_{M_k>n^{(1-\varepsilon)/s} } \leq
x  \Biggr) = L_{s,b'}(x),
\end{equation}
where the characteristic function for the distribution $L_{s,b'}$
is given in \eqref{char}. To get the limiting distribution of
$\frac{1}{n^{1/s}} E_\omega T_{\nu_n}$, we use \eqref{reflectexpand}
and rewrite this as
%
\begin{eqnarray}
\frac{1}{n^{1/s}} E_\omega T_{\nu_n} &= & \frac{1}{n^{1/s}}
\sum_{k=1}^n E_\omega^{\nu_{k-1}} \bar{T}_{\nu_k}^{(n)}
\mathbf{1}_{M_k > n^{(1-\varepsilon)/s}}
\label{bb}\\
&&{} + \frac{1}{n^{1/s}} \sum_{k=1}^n
E_\omega^{\nu_{k-1}} \bar{T}_{\nu_k}^{(n)}
\mathbf{1}_{M_k \leq n^{(1-\varepsilon)/s}}
\label{sb}\\
&&{} + \frac{1}{n^{1/s}}  \bigl( E_\omega T_{\nu_n}
-E_\omega \bar{T}_{\nu_n}^{(n)}  \bigr).
\label{rb}
\end{eqnarray}
Lemma \ref{smallblocks} gives that \eqref{sb} converges in
distribution (under $Q$) to 0. Also, we can use Lemma \ref{ETdiff}
to show that \eqref{rb} converges in distribution to 0 as well.
Indeed, for any $\delta >0$,
\begin{eqnarray*}
Q \bigl( E_\omega T_{\nu_n} - E_\omega\bar{T}_{\nu_n}^{(n)} >
\delta n^{1/s}  \bigr)
&\leq & n Q \bigl( E_\omega T_{\nu} - E_\omega\bar{T}_{\nu}^{(n)}
> \delta n^{1/s-1}  \bigr)
= \mathcal{O} ( n^s e^{-\delta 'b_n}  ).
\end{eqnarray*}
Therefore, $n^{-1/s}E_\omega T_{\nu_n}$ has the same limiting
distribution (under $Q$) as the right side of \eqref{bb}, which by
\eqref{refdist} is an $s$-stable distribution with distribution
function $L_{s,b'}$.
\end{pf*}

%
\section{Localization along a subsequence}\label{localization}

The goal of this section is to show
when $s<1$ that $P$-a.s. there exists a subsequence $t_m= t_m(\omega)$
of times such that the RWRE is essentially located in a section of
the environment of length $\log^2(t_m)$. This will essentially be
done by finding a ladder time whose crossing time is \emph{much}
larger than all the other ladder times before it. As a first step
in this direction, we prove that with strictly positive probability
this happens in the first $n$ ladder locations. Recall the definition
of $M_k$; cf.~\eqref{Mdef}.\looseness=-1
\begin{lem}\label{QBB}
Assume $s<1$. Then for any $C>1$, we have
\[
\liminf_{n\rightarrow\infty} Q \Biggl( \exists k \in[1, n/2]\dvtx
M_k \geq C \sum_{j\in[1,n]\backslash\{k\}}
E_\omega^{\nu _{j-1}} \bar{T}_{\nu_j}^{(n)} \Biggr) > 0.
\]
\end{lem}

\begin{pf}
Recall that $\bar{T}^{(n)}_x$ is the hitting time of $x$ by the RWRE
modified so that it never backtracks $b_n=\lfloor\log^2(n) \rfloor$
ladder locations.

To prove the lemma, first note that since $C>1$ and $E_\omega^{\nu_{k-1}}
\bar{T}^{(n)}_{\nu_k} \geq M_k$ there can only be at most one $k\leq
n$ with $M_k \geq C \sum_{k\neq j \leq n} E_\omega^{\nu_{j-1}} \bar
{T}_{\nu_j}^{(n)}$. Therefore,
%
\begin{eqnarray}\label{onebigblock}
&& Q \Biggl( \exists k\in[1,n/2]\dvtx  M_k \geq C        \sum_{j\in
[1,n]\backslash\{k\}}        E_\omega^{\nu_{j-1}} \bar
{T}^{(n)}_{\nu_j}  \Biggr)
\nonumber\\[-8pt]
\\[-8pt]
&&\qquad = \sum_{k=1}^{n/2} Q \Biggl( M_k \geq C
\sum_{j\in[1,n]\backslash\{k\}}
 E_\omega^{\nu_{j-1}} \bar{T}_{\nu_j}^{(n)} \Biggr)
 \nonumber
\end{eqnarray}
Now, define the events
%
\begin{eqnarray}\label{FGevents}
F_{n} &:=& \{\nu_j-\nu_{j-1} \leq b_n , \ \forall j\in(-b_n,n] \},
\nonumber\\[-8pt]
\\[-8pt]
G_{k,n,\varepsilon}&:=& \bigl\{ M_j \leq n^{(1-\varepsilon)/s},
\ \forall j\in(k,k+b_n] \bigr\}.
\nonumber
\end{eqnarray}
$F_n$ and $G_{k,n,\varepsilon}$ are both \emph{typical} events. Indeed,
from Lemma \ref{nutail}, $Q(F_{n}^c) \leq(b_n+n) Q(\nu> b_n) =
\mathcal{O}(n e^{-C_2 b_n})$, and from \eqref{Mtail}, we have
$Q(G_{k,n,\varepsilon}^c) \leq b_n Q(M_1 > n^{(1-\varepsilon)/s}) =
o(n^{-1+2\varepsilon})$.
Now, from \eqref{QET}, adjusted for reflections, we have for
any $j\in[1,n]$ that
\begin{eqnarray*}
&& E_\omega^{\nu_{j-1}} \bar{T}_{\nu_j}^{(n)}
\\
&&\qquad = (\nu_j - \nu_{j-1}) + 2\sum_{l=\nu_{j-1}}^{\nu_j-1} W_{\nu_{j-1-b_n},l}
\\
&&\qquad = (\nu_j - \nu_{j-1}) + 2\sum_{ \nu_{j-1} \leq i\leq l < \nu_j}\Pi_{i,l}
+  2\sum_{\nu_{j-1-b_n}< i < \nu_{j-1} \leq l < \nu_j}
\Pi_{i,\nu_{j-1}-1}\Pi_{\nu_{j-1}, l}
\\
&&\qquad \leq (\nu_j - \nu_{j-1}) + 2 (\nu_j-\nu_{j-1} )^2 M_j
+ 2(\nu_j-\nu_{j-1})(\nu_{j-1}-\nu_{j-1-b_n})M_j,
\end{eqnarray*}
where in the last inequality we used the facts that $\Pi_{\nu_{j-1},
i-1} \geq1$ for $\nu_{j-1} < i < \nu_j$ and $\Pi_{i, \nu_{j-1}-1}
< 1$ for all $i<\nu_{j-1}$. Then on the event $F_n\cap
G_{k,n,\varepsilon}$, we have for $k+1\leq j\leq k+b_n$ that
\begin{eqnarray*}
E_\omega^{\nu_{j-1}} \bar{T}_{\nu_j}^{(n)} \leq b_n + 2 b_n^2
n^{(1-\varepsilon)/s} + 2 b_n^3 n^{(1-\varepsilon)/s} \leq5 b_n^3
n^{(1-\varepsilon)/s},
\end{eqnarray*}
where for the first inequality we used that on the event
$F_n\cap G_{k,n,\varepsilon}$ we have $\nu_j-\nu_{j-1} \leq b_n$ and
$M_1\leq n^{(1-\varepsilon)/s}$. Then using this, we get
\begin{eqnarray*}
&& Q \Biggl( M_k \geq C        \sum_{j\in[1,n]\backslash\{k\}}
     E_\omega^{\nu_{j-1}} \bar{T}_{\nu_j}^{(n)} \Biggr)
\\
&&\qquad  \geq  Q \bigl( M_k \geq C  \bigl( E_\omega
\bar{T}_{\nu_{k-1}}^{(n)} + 5b_n^4 n^{(1-\varepsilon)/s}
+ E_\omega^{\nu_{k+b_n}}\bar{T}_{\nu_n}^{(n)}  \bigr), F_n,
G_{k,n,\varepsilon}  \bigr)
\\
&&\qquad \geq Q ( M_k \geq C n^{1/s}, \nu_{k}-\nu_{k-1} \leq b_n )
\\
&&\qquad \quad {} \times Q \bigl( E_\omega\bar{T}_{\nu_{k-1}}^{(n)} + 5b_n^4
n^{(1-\varepsilon)/s} + E_\omega^{\nu_{k+b_n}}\bar{T}_{\nu_n}^{(n)}
\leq n^{1/s}, \tilde{F}_n, G_{k,n,\varepsilon}  \bigr),
\end{eqnarray*}
where $\tilde{F}_n := \{\nu_j-\nu_{j-1} \leq b_n, \quad\forall j\in
(-b_n, n]\backslash\{k\} \} \supset F_n $.
In the last inequality, we used the fact that
$E_\omega^{\nu_{j-1}} \bar{T}_{\nu_j}^{(n)}$ is independent of
$M_k$ for $j<k$ or $j>k+b_n$. Note that we can replace
$\tilde{F}_n$ by $F_n$ in the last line above because it will only
make the probability smaller. Then using the above and the fact that
$E_\omega\bar{T}_{\nu_{k-1}}^{(n)} + E_\omega^{\nu_{k+b_n}}\bar{T}_{\nu_n}^{(n)}
\leq E_\omega T_{\nu_n,}$ we have
\begin{eqnarray*}
&& Q \Biggl( M_k \geq C        \sum_{j\in[1,n]\backslash\{k\}}
       E_\omega^{\nu_{j-1}} \bar{T}_{\nu_j}^{(n)} \Biggr)
\\
&&\qquad \geq Q ( M_k \geq C n^{1/s},
\nu_{k}-\nu_{k-1} \leq b_n  )
\\
&&\qquad \quad {}\times
Q \bigl(E_\omega T_{\nu_n} \leq n^{1/s}-5 b_n^4
n^{(1-\varepsilon)/s}, F_n, G_{k,n,\varepsilon} \bigr)
\\
&&\qquad \geq \bigl( Q(M_1\geq Cn^{1/s}) - Q(\nu> b_n)  \bigr)
\\
&&\qquad \quad {}\times
\bigl( Q\bigl( E_\omega T_{\nu_n}\leq n^{1/s}(1-5b_n^4 n^{-\varepsilon/s}) \bigr)
- Q(F_n^c) - Q(G_{k,n,\varepsilon}^c)  \bigr)
\\
&&\qquad \sim C_5 C^{-s} L_{s,b'}(1) \frac{1}{n},
\end{eqnarray*}
where the asymptotics in the last line are from
\eqref{Mtail} and Theorem \ref{refstable}. Combining the last
display and \eqref{onebigblock} proves the lemma.
\end{pf}

In Section \ref{stablecrossing}, we showed that the proper scaling
for $E_\omega T_{\nu_n}$ (or $E_\omega\bar{T}^{(n)}_{\nu_n}$) was
$n^{-1/s}$. The following lemma gives a bound on the moderate
deviations under the measure $P$.
\begin{lem}\label{mdevu}
Assume $s\leq1$. Then for any $ \delta >0$,
\[
P ( E_\omega T_{\nu_n} \geq n^{1/s + \delta }  )
= o(n^{-\delta s/2}).
\]
\end{lem}
\begin{pf}
First, note that
%
\begin{equation}\label{md1}
P(E_\omega T_{\nu_n} \geq n^{1/s+ \delta })
\leq P(E_\omega T_{2\bar \nu n} \geq
n^{1/s+\delta }) + P(\nu_n \geq2 \bar\nu n ),
\end{equation}
where $\bar\nu:= E_P \nu$. To handle the second term on the
right side of \eqref{md1} we note that $\nu_n$ is the sum of
$n$ i.i.d. copies of $\nu$, and that $\nu$ has exponential tails (by
Lemma~\ref{nutail}).
Therefore, Cram\'er's theorem (\cite{dzLDTA}, Theorem 2.2.3)
gives that $P({\nu_n}/{n} \geq 2\bar\nu) = \mathcal{O}( e^{-\delta ' n} )$
for some $\delta ' > 0$.

To handle the first term on the right side of \eqref{md1}, we
note that for any $\gamma< s$ we have
$E_P(E_\omega T_1)^\gamma< \infty$. This follows from the fact that
$P(E_\omega T_1 > x) = P(1 + 2W_0 > x) \sim K 2^s x^{-s}$ by
\eqref{QET} and \eqref{PWtail}. Then by Chebyshev's inequality and
the fact that $\gamma< s \leq1$, we have
%
\begin{equation}\label{Tnbig}
P ( E_\omega T_{2\bar\nu n} \geq n^{1/s+\delta }  )
\leq \frac{E_P ( \sum_{k=1}^{2 \bar\nu n} E_\omega^{k-1} T_k
)^\gamma}{n^{\gamma(1/s+\delta )}}
\leq\frac{2 \bar\nu n E_P(E_\omega T_1)^\gamma}{n^{\gamma(1/s+\delta )}}.
\end{equation}
Then choosing $\gamma$ arbitrarily close to $s$, we can have that this
last term is $o(n^{-\delta s/2})$.
\end{pf}

Throughout the remainder of the paper, we will use the following
subsequences of integers:
%
\begin{equation}\label{nkdkdef}
n_k := 2^{2^k},\qquad d_k := n_k-n_{k-1}
\end{equation}
Note that $n_{k-1} = \sqrt{n_k}$ and so $d_k \sim n_k$ as
$k\rightarrow\infty$.
\begin{cor}\label{subseq1}
For any $k$, define
\[
\mu_k := \max \bigl\{ E_\omega^{\nu_{j-1}} \bar{T}_{\nu
_j}^{(d_k)} \dvtx n_{k-1} < j \leq n_k  \bigr\}.
\]
If $s<1$, then
\[
\lim_{k\rightarrow\infty} \frac{ E_\omega^{\nu_{n_{k-1}}}
\bar{T}_{\nu _{n_{k}}}^{(d_k)} - \mu_k }
{E_\omega\bar{T}_{\nu_{n_k}}^{(d_k)} -\mu_k} = 1, \qquad P\mbox{-a.s.}
\]
\end{cor}
\begin{pf}
Let $\varepsilon>0$. Then
%
\begin{eqnarray}\label{shortlong}
&& P \biggl( \frac{ E_\omega^{\nu_{n_{k-1}}} \bar{T}_{\nu_{n_{k}}}^{(d_k)}
- \mu_k }{E_\omega\bar{T}_{\nu_{n_k}}^{(d_k)} - \mu_k} \leq
1-\varepsilon \biggr)
\nonumber\\
&&\qquad = P \biggl( \frac{E_\omega\bar{T}_{\nu_{n_{k-1}}}^{(d_k)}}{E_\omega
\bar{T}_{\nu_{n_k}}^{(d_k)} - \mu_k} \geq\varepsilon \biggr)
\\
&&\qquad \leq P \bigl( E_\omega\bar{T}_{\nu_{n_{k-1}}}^{(d_k)}
\geq n_{k-1}^{1/s+\delta }  \bigr)
+ P \bigl( E_\omega\bar{T}_{\nu_{n_k}}^{(d_k)} -
\mu_k \leq\varepsilon^{-1} n_{k-1}^{1/s+\delta }  \bigr).
\nonumber
\end{eqnarray}
Lemma \ref{mdevu} gives that $P ( E_\omega\bar{T}_{\nu
_{n_{k-1}}}^{(d_k)} \geq n_{k-1}^{1/s+\delta }  ) \leq P (
E_\omega T_{\nu_{n_{k-1}}} \geq n_{k-1}^{1/s+\delta }  ) =\break
o(n_{k-1}^{-\delta s/2})$.
To handle the second term in the right side of \eqref{shortlong},
note that if $\delta <\frac{1}{3s}$, then the subsequence
$n_k$ grows fast enough such that for all $k$ large
enough $n_k^{1/s-\delta } \geq\varepsilon^{-1} n_{k-1}^{1/s+\delta }$.
Therefore, for $k$ sufficiently large and $\delta < \frac{1}{3s}$, we have
\[
P \bigl( E_\omega\bar{T}_{\nu_{n_k}}^{(d_k)} - \mu_k \leq
\varepsilon^{-1} n_{k-1}^{1/s+\delta }  \bigr)
\leq P \bigl(E_\omega\bar{T}_{\nu_{n_k}}^{(d_k)}
- \mu_k \leq n_{k}^{1/s-\delta } \bigr).
\]
However, $ E_\omega\bar{T}^{(d_k)}_{\nu_{n_k}} - \mu_k \leq n_{k}^{1/s-\delta } $
implies that $M_j < E_\omega^{\nu_{j-1}} \bar{T}_{\nu_j}^{(d_k)}
\leq n_k^{1/s-\delta }$ for at least $n_k-1$ of the $j\leq n_k$. Thus,
since $P(M_1> n_k^{1/s-\delta })\sim C_5 n_k^{-1+\delta s}$, we have that
%
\begin{eqnarray}\label{longsmall}
P \bigl( E_\omega\bar{T}_{\nu_{n_k}}^{(d_k)} - \mu_k \leq
\varepsilon^{-1} n_{k-1}^{1/s+\delta }  \bigr)
&\leq & n_k  \bigl( 1 - P ( M_1 > n_k^{1/s-\delta }  )  \bigr)^{n_k-1}
\nonumber\\[-8pt]
\\[-8pt]
&=&  o(e^{-n_k^{\delta s/2}}).
\nonumber
\end{eqnarray}
Therefore, for any $\varepsilon>0$ and $\delta <\frac{1}{3s}$, we
have that
\[
P \biggl( \frac{ E_\omega^{\nu_{n_{k-1}}} \bar{T}_{\nu_{n_{k}}}^{(d_k)}
- \mu_k }{E_\omega\bar{T}_{\nu_{n_k}}^{(d_k)} - \mu_k}
\leq1-\varepsilon \biggr) = o (n_{k-1}^{-\delta s/2} ).
\]
By our choice of $n_k$, the sequence $n_{k-1}^{-\delta s/2}$ is
summable in $k$. Applying the Borel--Cantelli lemma completes the proof.
\end{pf}
\begin{cor}\label{ssdom}
Assume $s<1$. Then $P$-a.s. there exists a random subsequence
$j_m= j_m(\omega)$ such that
\[
M_{j_m} \geq m^2 E_\omega\bar{T}_{\nu_{j_m-1}}^{(j_m)}.
\]
\end{cor}
\begin{pf}
Recall the definitions of $n_k$ and $d_k$ in \eqref{nkdkdef}. Then for
any $C>1$, define the event
\[
D_{k,C}:= \bigl\{\exists j\in(n_{k-1}, n_{k-1}+ d_k/2]\dvtx
M_j \geq C \bigl(E_\omega^{\nu_{n_{k-1}}}\bar{T}_{\nu_{j-1}}^{(d_k)}
+ E_\omega^{\nu_j}\bar{T}_{\nu_{n_k}}^{(d_k)}  \bigr)  \bigr\}.
\]
Note that due to the reflections, the event $D_{k,C}$ depends only on
the environment from $\nu_{n_{k-1}-b_{d_k}}$ to $\nu_{n_{k}}-1$. Then
since $n_{k-1} - b_{d_k} > n_{k-2}$ for all $k \geq4$,
we have that the events $\{ D_{2k,C} \}_{k=2}^{\infty}$
are all independent. Also, since the events do not involve the
environment to the left of $0$, they have the same probability under
$Q$ as under $P$. Then since $Q$ is stationary under shifts of
$\nu_i$, we have that for $k\geq4$,
\begin{eqnarray*}
P(D_{k,C}) = Q(D_{k,C})
= Q \bigl(\exists j\in[1,d_k/2]\dvtx
M_j \geq C \bigl(E_\omega\bar{T}_{\nu_{j-1}}^{(d_k)}
+ E_\omega^{\nu_j}\bar{T}_{\nu_{d_k}}^{(d_k)}  \bigr)  \bigr).
\end{eqnarray*}
Thus, for any $C>1$, we have by Lemma \ref{QBB} that
$\liminf_{k\rightarrow\infty} P(D_{k,C}) > 0$.
This combined with the fact that the events $\{D_{2k,C}\}_{k=2}^\infty$
are independent gives that for any $C>1$ infinitely many of
the events $D_{2k,C}$ occur $P$-a.s. Therefore, there exists a
subsequence $k_m$ of integers such that for each $m$,
there exists $j_m \in(n_{k_m-1} , n_{k_m-1}+d_{k_m}/2 ]$ such that
\[
M_{j_m} \geq2m^2  \bigl( E_{\omega}^{\nu_{n_{k_m-1}}}
\bar{T}_{\nu_{j_m-1}}^{(d_{k_m})} + E_\omega^{\nu_{j_m}}
\bar{T}_{\nu_{n_{k_m}}}^{(d_{k_m})}  \bigr)
= 2m^2 \bigl(E_\omega^{\nu_{n_{k_m-1}}}
\bar{T}_{\nu_{n_{k_m}}}^{(d_{k_m})} - \mu_{k_m} \bigr) ,
\]
where the second equality holds due to our choice of $j_m$, which
implies that $\mu_{k_m} = E_\omega^{\nu_{j_m -1}}
\bar{T}_{\nu_{j_m}}^{(d_{k_m})}$. Then
by Corollary \ref{subseq1}, we have that for all $m$ large enough
\[
M_{j_m}\geq2m^2 \bigl ( E_\omega^{\nu_{k_m -1}}
\bar{T}_{\nu_{n_{k_m}}}^{(d_{k_m})} - \mu_{k_m}  \bigr)
\geq m^2  \bigl( E_\omega\bar{T}_{\nu_{n_{k_m}}}^{(d_{k_m})}
- \mu_{k_m}  \bigr)
\geq m^2 E_\omega \bar{T}_{\nu_{j_m-1}}^{(d_{k_m})},
\]
where the last inequality is because
$\mu_{k_m} = E_\omega^{\nu_{j_m -1}} \bar{T}_{\nu_{j_m}}^{(d_{k_m})}$.
Now, for all $k$ large enough, we have $n_{k-1} + d_k/2 < d_k$.
Thus, we may assume (by possibly choosing a further subsequence)
that $j_m < d_{k_m}$ as well, and since allowing less
backtracking only decreases the crossing time we have
\[
M_{j_m} \geq m^2 E_\omega\bar{T}_{\nu_{j_m-1}}^{(d_{k_m})}
\geq m^2 E_\omega\bar{T}_{\nu_{j_m-1}}^{(j_m)}.
\]\upqed
\end{pf}

The following lemma shows that the reflections that we have
been using this whole time really do not affect the random walk.
Recall the coupling of $X_t$ and $\bar{X}_t^{(n)}$
introduced after \eqref{bdef}.
%
\begin{lem} \label{bt}
\[
\lim_{n\rightarrow\infty} P_\omega
\bigl( T_{\nu_{n-1}} \neq\bar{T}_{\nu_{n-1}}^{(n)}  \bigr) = 0,
\qquad P\mbox{-a.s.}
\]
\end{lem}
\begin{pf}
Let $\varepsilon>0$. By Chebyshev's inequality,
\[
P \bigl( P_\omega \bigl( T_{\nu_{n-1}} \neq\bar{T}_{\nu_{n-1}}^{(n)} \bigr)
> \varepsilon \bigr) \leq\varepsilon^{-1} \mathbb{P}
\bigl( T_{\nu_{n-1}} \neq\bar{T}_{\nu_{n-1}}^{(n)} \bigr) .
\]
Thus, by the Borel--Cantelli lemma, it is
enough to prove that $\mathbb{P} ( T_{\nu_{n-1}} \neq\bar
{T}_{\nu_{n-1}}^{(n)}
 )$ is summable. Now, the event
$T_{\nu_{n-1}} \neq\bar{T}_{\nu_{n-1}}^{(n)} $ implies that there is
an $i < \nu_{n-1}$ such that after reaching $i$ for the first time,
the random walk then backtracks a distance of $b_{n}$.
Thus, again letting $\bar\nu= E_P \nu$, we have
\begin{eqnarray*}
\mathbb{P} \bigl( T_{\nu_{n-1}} \neq\bar{T}_{\nu_{n-1}}^{(n)}
 \bigr) &\leq & P\bigl(\nu_{n-1} \geq2\bar\nu(n-1)\bigr)
 + \sum_{i=0}^{2\bar\nu (n-1)} \mathbb{P}^i(T_{i-b_{n}}< \infty)
 \\
&=& P\bigl(\nu_{n-1} \geq2\bar\nu(n-1)\bigr)
+ 2\bar\nu(n-1) \mathbb{P}(T_{-b_{n}}< \infty).
\end{eqnarray*}
As noted in Lemma \ref{mdevu}, $P(\nu_{n-1} \geq2\bar\nu(n-1)) =
\mathcal{O}( e^{-\delta ' n})$, so we need only to show that
$n\mathbb{P}(T_{-b_{n}} < \infty)$
is summable. However, \cite{gsMVSS}, Lemma 3.3, gives that
there exists a constant $C_7$ such that for any $k\geq1$,
%
\begin{equation}\label{backtracktail}
\mathbb{P}( T_{-k} < \infty) \leq e^{-C_7 k}.
\end{equation}
Thus, $n\mathbb{P}(T_{-b_{n}} < \infty) \leq n
e^{-C_7 b_{n}}$ which is summable by the definition of $b_n$.
\end{pf}

We define the random variable $N_t:= \max\{k: \exists n\leq t,
X_n = \nu_k\}$ to be the maximum number of ladder
locations crossed by the random walk by time~$t$.

\begin{lem}\label{seperation}
\[
\lim_{t\rightarrow\infty} \frac{\nu_{N_t}-X_t}{\log^2(t)} = 0,
\qquad\mathbb{P}\mbox{-a.s.}
\]
\end{lem}
\begin{pf}
Let $\delta > 0$. If we can show that
$\sum_{t=1}^\infty\mathbb{P}(|N_t-X_t|\geq \delta \log^2 t) < \infty$,
then by the Borel--Cantelli lemma, we will be done.
Now, the only way that $N_t$ and $X_t$ can differ by more than $\delta
\log
^2 t$ is if either one of the gaps between the first~$t$ ladder times
is larger than $\delta \log^2 t$ or if for some $i<t$ the random walk
backtracks $\delta \log^2 t$ steps after first reaching $i$. Thus,
%
\begin{eqnarray}\label{backtrack}
\quad
&& \mathbb{P}(|N_t-X_t|\geq\delta \log^2 t)
\nonumber\\[-8pt]
\\[-8pt]
&&\qquad \leq P (\exists j\in [1,t+1]\dvtx \nu
_j-\nu_{j-1} > \delta \log^2 t  )
+ t \mathbb{P}\bigl( T_{-\lceil\delta \log^2 t \rceil} < T_1\bigr)
\nonumber
\end{eqnarray}
So, we need only to show that the two terms on the right side are
summable. For the first term, we use Lemma \ref{nutail} and note that
\begin{eqnarray*}
P (\exists j\in[1,t+1]\dvtx \nu_j-\nu_{j-1} > \delta \log^2 t )
&\leq& (t+1)P(\nu> \delta \log^2 t)
\\
&\leq& (t+1)C_1e^{-C_2\delta \log^2 t},
\end{eqnarray*}
which is summable in $t$. By \eqref{backtracktail}, the second
term on the right side of \eqref{backtrack} is also summable.
\end{pf}
\begin{pf*}{Proof of Theorem \ref{local}}
By Corollary \ref{ssdom}, $P$-a.s. there exists a subsequence
$j_m(\omega )$ such that
$M_{j_m} \geq m^2 E_\omega\bar{T}_{\nu_{j_m-1}}^{(j_m)}$. Define
$ t_m = t_m(\omega) = \frac{1}{m} M_{j_m}$
and $u_m=u_m(\omega)= \nu_{j_m-1}$. Then
\[
P_\omega \biggl( \frac{X_{t_m} - u_m}{\log^2 t_m} \notin[-\delta ,
\delta ]  \biggr) \leq P_\omega(N_{t_m}\neq j_m-1 )
+ P_\omega(|\nu_{N_{t_m}}-X_{t_m}| > \delta \log^2 t_m ).
\]
From Lemma \ref{seperation}, the second term goes to zero as
$m\rightarrow\infty$. Thus, we only need to show that
%
\begin{equation}\label{loc}
\lim_{m\rightarrow\infty} P_\omega( N_{t_m} = j_m -1 ) = 1.
\end{equation}
To see this, first note that
\begin{eqnarray*}
P_\omega ( N_{t_m} < j_m -1  )
&=& P_\omega ( T_{\nu_{j_m-1}}> t_m )
\\
&\leq & P_\omega \bigl( T_{\nu_{j_m-1}} \neq \bar{T}_{\nu_{j_m-1}}^{(j_m)}  \bigr)
+ P_\omega \bigl( \bar{T}_{\nu_{j_m-1}}^{(j_m)} > t_m  \bigr).
\end{eqnarray*}
By Lemma \ref{bt},
$P_\omega ( T_{\nu_{j_m-1}} \neq\bar{T}_{\nu_{j_m-1}}^{(j_m)} )
\rightarrow0$ as $m\rightarrow\infty$, $P$-a.s. Also, by our
definition of $t_m$ and
our choice of the subsequence $j_m$, we have
\[
P_\omega \bigl( \bar{T}_{\nu_{j_m-1}}^{(j_m)} > t_m  \bigr)
\leq\frac{E_\omega\bar{T}_{\nu_{j_m-1}}^{(j_m)}}{t_m} =
\frac{m E_\omega\bar{T}_{\nu_{j_m-1}}^{(j_m)}}{M_{j_m}}\leq
\frac{1}{m} \mathop{\longrightarrow}_{m\rightarrow\infty} 0.
\]
It still remains to show $\lim_{m\rightarrow\infty} P_\omega (
N_{t_m} < j_m  ) = 1$. To prove this, first define the
stopping times $T_x^+:= \min\{n > 0\dvtx X_n =x \} $. Then
\begin{eqnarray*}
P_\omega ( N_{t_m} < j_m  )
&=& P_\omega( T_{\nu_{j_m}} > t_m )
\geq P_\omega^{\nu_{j_m -1}}
\biggl( T_{\nu_{j_m}} > \frac{1}{m} M_{j_m} \biggr)
\\
&\geq& P_\omega^{\nu_{j_m -1}}  ( T_{\nu_{j_m-1}}^+ <
T_{\nu_{j_m}}  ) ^{({1}/{m}) M_{j_m} } .
\end{eqnarray*}
Then using the hitting time calculations given in
\cite{zRWRE}, (2.1.4), we have that
\begin{eqnarray*}
P_\omega^{\nu_{j_m -1}}  ( T_{\nu_{j_m-1}}^+ < T_{\nu_{j_m}} )
= 1-\frac{1-\omega_{\nu_{j_m-1}}}{ R_{\nu_{j_m-1},\nu_{j_m}-1} }.
\end{eqnarray*}
Therefore, since $M_{j_m} \leq R_{\nu_{j_m-1},\nu_{j_m}-1}$, we have
\[
P_\omega ( N_{t_m} < j_m  ) \geq
\biggl ( 1-\frac{1-\omega_{\nu_{j_m-1}}}{ R_{\nu_{j_m-1},\nu
_{j_m}-1} } \biggr)^{({1}/{m}) M_{j_m}}
\geq \biggl( 1-\frac{1}{M_{j_m}} \biggr)^{({1}/{m}) M_{j_m}}
\mathop{\longrightarrow}_{m\rightarrow\infty} 1,
\]
thus proving \eqref{loc} and, therefore, the theorem.
\end{pf*}


\section{Nonlocal behavior on a random subsequence} \label{gaussian}

There are two main goals of this section. The first is to prove the
existence of random subsequences $x_m$ where the hitting times $T_{x_m}$
are approximately Gaussian random variables. This result is then
used to prove the existence of random times $t_m(\omega)$ in which the
scaling for the random walk is of the order $t_m^s$ instead
of $\log^2 t_m$ as in Theorem \ref{local}. However, before we can begin
proving a quenched CLT for the hitting times $T_n$ (at least
along a random subsequence), we first need to understand the
tail asymptotics of
$\operatorname{Var}_\omega T_\nu:= E_\omega((T_\nu-E_w T_\nu)^2)$,
the quenched variance of $T_\nu$.


\subsection{Tail asymptotics of $Q( \operatorname{Var}_\omega T_\nu> x)$}

The goal of this subsection is to prove the following theorem.
\begin{thm}\label{qVartail}
Let Assumptions \ref{essentialasm} and \ref{techasm} hold. Then with
$K_\infty>0$ the same as in Theorem \ref{Tnutail}, we have
%
\begin{equation}\label{vartail}
Q ( \operatorname{Var}_\omega T_\nu> x  )
\sim Q \bigl( (E_\omega T_\nu)^2 > x \bigr)
\sim K_\infty x^{-s/2}\qquad\mbox{as } x\rightarrow\infty,
\end{equation}
and for any $\varepsilon> 0$ and $x>0$,
%
\begin{equation}\label{rvartail}
Q \bigl( \operatorname{Var}_\omega\bar{T}_\nu^{(n)} > x n^{2/s} ,
M_1 > n^{(1-\varepsilon)/s}  \bigr)
\sim K_\infty x^{-s/2} \frac{1}{n}\qquad \mbox{as } n\rightarrow\infty.
\end{equation}
Consequently,
%
\begin{equation}\label{VarbigMsmall}
Q \bigl( \operatorname{Var}_\omega T_\nu> \delta n^{1/s}, M
_1\leq n^{(1-\varepsilon)/s} \bigr) = o(n^{-1}).
\end{equation}
\end{thm}

A formula for the quenched variance of crossing times is given in
\cite{gQCLT}, (2.2). Translating to our notation and simplifying, we have
the formula
%
\begin{equation}\label{qvar}
\qquad
\operatorname{Var}_\omega T_1
:= E_\omega(T_1 - E_\omega T_1)^2 = 4(W_{0}+ W_{0}^2)+ 8 \sum_{i<0}
\Pi_{i+1,0}(W_{i}+W_{i}^2).
\end{equation}
Now, given the environment the crossing times $T_j-T_{j-1}$ are independent.
Thus, we get the formula
%
\begin{eqnarray}\label{VarTnuexpand}
\operatorname{Var}_\omega T_\nu&=& 4\sum_{j=0}^{\nu-1}(W_{j}+W_{j}^2)
+ 8\sum_{j=0}^{\nu-1}\sum_{i< j} \Pi_{i+1,j}(W_{i}+W_{i}^2)
\nonumber\\
&=& 4\sum_{j=0}^{\nu-1}(W_{j}+W_{j}^2)
\nonumber\\[-8pt]
\\[-8pt]
&&{} + 8 R_{0,\nu-1}
\Biggl(W_{-1}+W_{-1}^2+ \sum_{i < -1} \Pi_{i+1,-1}(W_{i}+W_{i}^2) \Biggr)
\nonumber\\
&&{} + 8\sum_{0\leq i < j < \nu}\Pi_{i+1,j}(W_{i}+W_{i}^2).
\nonumber
\end{eqnarray}
%
We want to analyze the
tails of $\operatorname{Var}_\omega T_\nu$ by comparison
with $(E_\omega T_\nu)^2$. Using \eqref{ETnuexpand}, we have
\[
(E_\omega T_\nu)^2 =  \Biggl( \nu+ 2 \sum_{j=0}^{\nu-1} W_{j}  \Biggr)^2
= \nu^2 + 4 \nu\sum_{j=0}^{\nu-1} W_j + 4\sum_{j=0}^{\nu-1}W_j^2
+ 8\sum_{0\leq i<j<\nu} W_i W_j.
\]
Thus, we have
%
\begin{eqnarray}
\qquad
&& (E_\omega T_\nu)^2 - \operatorname{Var}_\omega T_\nu
\nonumber\\
&& \qquad = \nu^2 +4(\nu-1)\sum_{j=0}^{\nu-1} W_{j}
+ 8\sum_{0\leq i<j<\nu} W_{i} ( W_j - \Pi_{i+1,j}-\Pi_{i+1,j}W_i  )
\label{above}\\
&&\qquad \quad {} - 8 R_{0,\nu-1} \Biggl(W_{-1}+W_{-1}^2+
\sum_{i < -1} \Pi_{i+1,-1}(W_{i}+W_{i}^2) \Biggr)
\label{below}\\
&&\qquad =: D^+(\omega) - 8 R_{0,\nu-1} D^-(\omega).
\label{Ddef}
\end{eqnarray}
Note that $D^-(\omega)$ and $D^+(\omega)$ are nonnegative random variables.
The next few lemmas show that the tails of $D^+(\omega)$ and $R_{0,\nu-1}
D^-(\omega)$ are much smaller than the tails of $(E_\omega T_\nu)^2$.
\begin{lem}\label{abound}
For any $\varepsilon>0$, we have $Q ( D^+(\omega) > x  )
= o(x^{-s+\varepsilon})$.
\end{lem}
\begin{pf}
Notice first that from \eqref{ETnuexpand} we have
$\nu^2+4(\nu-1)\sum_{j=0}^{\nu-1} W_j \leq2\nu E_\omega T_\nu$. Also
we can rewrite $ W_j - \Pi_{i+1,j}-\Pi_{i+1,j}W_i = W_{i+2,j}$
when $i < j-1$ (this term is zero when $i=j-1$). Therefore,
\[
Q \bigl( D^+(\omega) > x  \bigr) \leq Q( 2\nu E_\omega T_\nu> x/2 ) +
Q \Biggl( 8 \sum_{i=0}^{\nu-3} \sum_{j=i+2}^{\nu-1} W_i W_{i+2,j}
> x/2  \Biggr).
\]
Lemma \ref{nutail} and Theorem \ref{Tnutail} give that
$Q (2\nu E_\omega T_\nu> x  ) \leq Q(2\nu> \log^2(x))
+ Q (E_\omega T_\nu> \frac{ x}{\log^2(x)}  ) = o(x^{-s+\varepsilon})$
for any $\varepsilon>0$. Thus, we need only prove that
$Q ( \sum_{i=0}^{\nu-3}\sum_{j=i+2}^{\nu-1} W_i W_{i+2,j} > x  )
= o(x^{-s+\varepsilon })$ for
any $\varepsilon>0$. Note that for $i<\nu$, we have $W_i= W_{0,i} +
\Pi_{0,i}W_{-1} \leq\Pi_{0,i}(i+1+W_{-1})$, thus,
%
\begin{eqnarray}
&& Q \Biggl( \sum_{i=0}^{\nu-3} \sum_{j=i+2}^{\nu-1} W_i W_{i+2,j} > x \Biggr)
 \nonumber\\
&&\qquad \leq Q \Biggl((\nu+W_{-1}) \sum_{i=0}^{\nu-3} \sum_{j=i+2}^{\nu
-1}\Pi_{0,i} W_{i+2,j} > x  \Biggr)
\nonumber\\
&&\qquad \leq Q\bigl(\nu> \log^2(x)/2\bigr)
+ Q\bigl(W_{-1}> \log^2(x)/2\bigr)
\label{nuplusW}\\
&&\qquad \quad {} + \sum_{i=0}^{\log^2(x)-3}   \sum_{j=i+2}^{\log
^2(x)-1} P \biggl( \Pi_{0,i} W_{i+2,j} > \frac{x}{\log^6 (x)}
\biggr) \label{pitimesW},
\end{eqnarray}
where we were able to switch to $P$ instead of $Q$ in the last
line because the event inside the probability only concerns the
environment to the right of $0$. Now, Lemmas \ref{nutail} and
\ref{Wtail} give that \eqref{nuplusW} is $o(x^{-s+\varepsilon})$ for any
$\varepsilon>0$, so we need only to consider~\eqref{pitimesW}. Under
the measure $P$, we have that $\Pi_{0,i}$ and $W_{i+2,j}$ are
independent, and by~\eqref{PWtail} we have $P(W_{i+2,j}> x)\leq
P(W_{j}>x) \leq K_1 x^{-s}$. Thus,
\begin{eqnarray*}
P \biggl( \Pi_{0,i} W_{i+2,j} > \frac{x}{\log^6 (x)}  \biggr)
&=& E_P \biggl[ P  \biggl( W_{i+2,j} > \frac{x}{\log^6(x) \Pi_{0,i}}
\Big| \Pi_{0,i}  \biggr)  \biggr]
\\
&\leq & K_1 \log^{6s}(x) x^{-s} E_P[\Pi_{0,i}^s].
\end{eqnarray*}
Then because $E_P \Pi_{0,i}^s = ( E_P \rho^s )^{i+1} = 1$ by
Assumption \ref{essentialasm}, we have
\begin{eqnarray*}
\sum_{i=0}^{\log^2(x)-3}   \sum_{j=i+2}^{\log^2(x)-1}
P \biggl( \Pi_{0,i} W_{i+2,j} > \frac{x}{\log^6 (x)}  \biggr)
&\leq & K_1 \log^{4+6s}(x) x^{-s}
\\
&=& o(x^{-s+\varepsilon}).
\end{eqnarray*}\upqed
\end{pf}
\begin{lem}\label{piW2}
For any $\varepsilon>0$,
%
\begin{equation}\label{piW2tail}
Q \bigl( D^-(\omega) > x  \bigr) = o(x^{-s+\varepsilon}),
\end{equation}
and thus for any $\gamma< s$,
%
\begin{equation}\label{EpiW2}
E_Q D^-(\omega)^\gamma< \infty.
\end{equation}
\end{lem}
\begin{pf}
It is obvious that \eqref{piW2tail} implies \eqref{EpiW2} and so
we will only prove the former. For any $i$, we may expand $W_i+W_i^2$ as
\begin{eqnarray*}
W_i+W_i^2 &=& \sum_{k\leq i} \Pi_{k,i}
+  \Biggl( \sum_{k\leq i} \Pi_{k,i}  \Biggr)^2
= \sum_{k\leq i} \Pi_{k,i} + \sum_{k\leq i} \Pi_{k,i}^2
+ 2\sum_{k\leq i}\sum_{l < k} \Pi_{k,i}\Pi_{l,i}
\\
&=& \sum_{k\leq i} \Pi_{k,i}  \Biggl( 1 + \Pi_{k,i} + 2 \sum_{l<k}
\Pi_{l,i}  \Biggr).
\end{eqnarray*}
Therefore, we may rewrite
%
\begin{eqnarray}\label{piW2expand}
D^-(\omega) &=&  W_{-1}+W_{-1}^2+ \sum_{i<-1}\Pi_{i+1,-1}(W_i+W_i^2)
\nonumber\\[-8pt]
\\[-8pt]
&=& \sum_{i\leq-1}\sum_{k\leq i} \Pi_{k,-1}
\Biggl( 1 + \Pi_{k,i} + 2 \sum_{l < k} \Pi_{l,i}  \Biggr).
\nonumber
\end{eqnarray}
Next, for any $c>0$ and $n\in\mathbb{N}$ define the event
\begin{eqnarray*}
E_{c,n} :\!\!&=&  \bigl\{ \Pi_{j,i} \leq e^{-c(i-j+1)},
\forall{-}n\leq i \leq-1 ,\forall j \leq i-n \bigr\}
\\
&=& \bigcap_{-n\leq i\leq-1}
\bigcap_{j \leq i-n} \bigl\{\Pi_{j,i} \leq e^{-c(i-j+1)} \bigr\}.
\end{eqnarray*}
Now, under the measure $Q$, we have that $\Pi_{k,-1} < 1$ for all
$k\leq-1$, and thus on the event $E_{c,n}$ we have using the
representation in \eqref{piW2expand} that
%
\begin{eqnarray}\label{piW2ub}
\qquad
D^-(\omega) &=& \sum_{ i\leq-1}\sum_{k\leq i} \Pi_{k,-1}
\Biggl( 1 + \Pi_{k,i}+ 2 \sum_{l < k} \Pi_{l,i}  \Biggr)
\nonumber\\
& \leq& \sum_{-n\leq i\leq-1}  \Biggl(
\sum_{k\leq i} \Pi_{k,i}(\Pi_{i+1,-1}+\Pi_{k,-1})
\nonumber\\
&&\hspace*{44pt}{}
+ 2 \sum_{i-n<k\leq i}  \sum_{l<k} \Pi_{l,i}
+ 2 \sum_{l<k \leq i-n} e^{ck} \Pi_{l,i}  \Biggr)
\nonumber\\
&&{} + \sum_{i<-n}  \Biggl(\sum_{k\leq i } e^{ck}
+ \sum_{k\leq i} e^{ck}\Pi_{k,i}
+ 2 \sum_{l<k\leq i } e^{ck} \Pi_{l,i}  \Biggr)
\nonumber\\
&\leq& \sum_{-n\leq i\leq-1}
\Biggl( (2+n)W_i + 2 \sum_{l<k \leq i-n}  e^{ck}e^{-c(i-l+1)}  \Biggr)
\\
&&{} + \sum_{i<-n} \Biggl ( \frac{e^{c(i+1)}}{e^c-1} + e^{ci} W_i +
\frac{2 e^{c(i+1)}}{e^c-1} \sum_{l<i} \Pi_{l,i}  \Biggr)
\nonumber\\
&\leq& (2+n)\sum_{-n\leq i\leq-1} W_i
+ \frac{2e^{-c(2n-1)}}{(e^c-1)^3(e^c+1)}
+ \frac{e^{-c(n-1)}}{(e^c-1)^2}
\nonumber\\
&&{} + \sum_{i<-n} e^{ci} W_i  \biggl(1+\frac{2 e^c}{e^c-1} \biggr)
\nonumber\\
& \leq& (2+n)\sum_{-n\leq i\leq-1} W_i + \frac{ e^c(1+e^{2c})}
{(e^c-1)^3(e^c+1)}+ \frac{3 e^c - 1}{e^c-1} \sum_{i<-n} e^{ci} W_i.
\nonumber
\end{eqnarray}
Then using \eqref{piW2ub} with $n$ replaced by $\lfloor\log^2 x
\rfloor=b_x$ we have
%
\begin{eqnarray}\label{onE}
Q \bigl( D^-(\omega) > x  \bigr)
 &\leq& Q ( E_{c,b_x}^c  ) + \mathbf{1}_{\{
{ e^c(1+e^{2c}) }/{((e^c-1)^3(e^c+1))} > x/3\} }
\nonumber\\
&&{}
+ Q \Biggl( \sum_{-b_x \leq i\leq-1} W_i > \frac{x}{3(2+ b_x)}  \Biggr)
\\
&&{} + Q \Biggl( \sum_{i<-1} e^{ci} W_i > \frac{(e^c-1)x}{3(3e^c-1)}  \Biggr).
\nonumber
\end{eqnarray}
Now, for any $0<c<-E_P \log\rho$, Lemma \ref{nutail} gives that
$Q(\Pi_{i,j}
> e^{-c (j-i+1)}) \leq\frac{A_c}{P(\mathcal{R})} e^{-\delta _c (j-i+1)}$
for some $\delta _c,A_c>0$. Therefore,
%
\begin{eqnarray}\label{QEc}
Q(E_{c,n}^c) &\leq& \sum_{-n\leq i\leq-1} \sum_{j \leq i-n}
Q\bigl(\Pi_{j,i} > e^{-c(i-j+1)}\bigr)
\nonumber\\[-8pt]
\\[-8pt]
&\leq&
\frac{n A_c e^{-\delta _c n}}{P(\mathcal{R})(e^{\delta _c}-1)} =
o(e^{-\delta _c n/2}).
\nonumber
\end{eqnarray}
Thus, for any $0<c<- E_P \log\rho$, we have that
the first two terms on the right side of \eqref{onE} are
decreasing in $x$ of order
$o(e^{-\delta _c b_x /2}) = o(x^{-s+\varepsilon})$. To handle last two
terms in the right side of \eqref{onE},
note first that from \eqref{PWtail},
$Q ( W_i > x  ) \leq\frac{1}{P(\mathcal{R})} P( W_i > x)
\leq
\frac{K_1}{P(\mathcal{R})} x^{-s}$ for any $x>0$ and any $i$. Thus,
\begin{eqnarray*}
Q \Biggl( \sum_{-b_x \leq i\leq-1} W_i > \frac{x}{3(2+ b_x)}  \Biggr)
&\leq& \sum_{-b_x\leq i\leq-1} Q \biggl( W_i > \frac{x}{3(2+ b_x)b_x}
 \biggr)
 \\
& =& o(x^{-s+\varepsilon}),
\end{eqnarray*}
and since $\sum_{i=1}^\infty e^{-c i/2} = (e^{c/2} -1)^{-1}$, we
have
\begin{eqnarray*}
&& Q \Biggl( \sum_{i<-1} e^{ci} W_i > \frac{(e^c-1)x}{9e^c-3}  \Biggr)
\\
&&\qquad \leq  Q \Biggl(\sum_{i=1}^\infty e^{-c i }W_{-i}
> \frac{(e^c-1)x}{9e^c-3}(e^{c/2} - 1)
\sum_{i=1}^{\infty} e^{-c i / 2}  \Biggr)
\\
&&\qquad \leq \sum_{i=1}^\infty Q
\biggl( W_{-i} > \frac{(e^c-1)(e^{c/2} -1)}{9e^c-3} x e^{c i/2}  \biggr)
\\
&&\qquad \leq \frac{K_1 (9e^c-3)^s }{P(\mathcal{R})(e^c-1)^s(e^{c/2}-1)^s} x^{-s}
\sum_{i=1}^\infty e^{-c s i /2} = \mathcal{O}(x^{-s}).
\end{eqnarray*}\upqed
\end{pf}
\begin{cor}\label{bbound}
For any $\varepsilon>0$, $Q (R_{0,\nu-1} D^-(\omega) > x) =
o(x^{-s + \varepsilon})$.
\end{cor}
\begin{pf}
From \eqref{Rtail}, it is easy to see that for any $\gamma< s$ there
exists a $K_\gamma> 0$ such that
$P(R_{0,\nu-1} > x) \leq P(R_{0} > x) \leq K_\gamma x^{-\gamma}$.
Then letting $\mathcal{F}_{-1} = \sigma(\ldots,\omega_{-2},
\omega_{-1})$, we have that
\begin{eqnarray*}
Q \bigl(R_{0,\nu-1}D^-(\omega) > x  \bigr)
&=& E_Q \biggl[ Q\biggl(R_{0,\nu-1} > \frac{x}{D^-(\omega)}
\Big| \mathcal{F}_{-1}  \biggr)  \biggr]
\\
&\leq & K_\gamma x^{-\gamma} E_Q( D^-(\omega) )^{\gamma}.
\end{eqnarray*}
Since $\gamma< s$, the expectation in the last expression
is finite by \eqref{EpiW2}. Choosing $\gamma= s-\frac{\varepsilon}{2}$
completes the proof.
\end{pf}
\begin{pf*}{Proof of Theorem \ref{qVartail}}
Recall from \eqref{Ddef} that
%
\begin{equation}\label{compare}
 ( E_\omega T_\nu )^2 - D^+(\omega)
\leq \operatorname{Var}_\omega T_\nu
\leq  ( E_\omega T_\nu )^2 + 8R_{0,\nu-1}D^-(\omega).
\end{equation}
The lower bound in \eqref{compare} gives that for any $\delta >0$,
\[
Q(\operatorname{Var}_\omega T_\nu> x)
\geq Q\bigl ( (E_\omega T_\nu)^2 > (1+\delta)x  \bigr)
- Q \bigl(D^+(\omega) > \delta x  \bigr).
\]
Thus, from Lemma \ref{abound} and Theorem \ref{Tnutail}, we have
that
%
\begin{eqnarray}\label{comparel}
\liminf_{x\rightarrow\infty} x^{s/2}
Q(\operatorname{Var}_\omega T_\nu> x) \geq K_\infty
(1+\delta )^{-s/2}.
\end{eqnarray}
Similarly, the upper bound in \eqref{compare} and Corollary
\ref{bbound} give that for any $\delta > 0$,
\begin{eqnarray*}
Q(\operatorname{Var}_\omega T_\nu> x)
&\leq& Q \bigl( (E_\omega T_\nu)^2 > (1-\delta)x  \bigr)
+ Q \bigl( 8 R_{0,\nu-1} D^-(\omega) > \delta x  \bigr),
\end{eqnarray*}
and then Corollary \ref{bbound} and Theorem \ref{Tnutail} give
%
\begin{equation}\label{compareu}
\limsup_{x\rightarrow\infty} x^{s/2}
Q(\operatorname{Var}_\omega T_\nu> x) \leq K_\infty(1-\delta )^{-s/2} .
\end{equation}
Letting $\delta \rightarrow0$ in \eqref{comparel} and
\eqref{compareu} completes the proof of \eqref{vartail}.

Essentially the same proof works for \eqref{rvartail}. The
difference is that when evaluating the difference $(E_\omega
\bar{T}_\nu^{(n)})^2 - \operatorname{Var}_\omega\bar{T}_\nu^{(n)}$
the upper and
lower bounds in \eqref{above} and \eqref{below} are smaller in
absolute value. This is because every instance of $W_i$ is
replaced by $W_{\nu_{-b_n}+1,i} \leq W_i$ and the sum in
\eqref{below} is taken only over $\nu_{-b_n} < i <-1$. Therefore,
the following bounds still hold:
%
\begin{equation}\label{comparer}
\bigl ( E_\omega\bar{T}^{(n)}_\nu \bigr)^2 - D^+(\omega)
\leq \operatorname{Var}_\omega
\bar{T}^{(n)}_\nu\leq
\bigl( E_\omega\bar{T}^{(n)}_\nu \bigr)^2 + 8R_{0,\nu-1} D^-(\omega).
\end{equation}
The rest of the proof then follows in the same manner, noting that
from Lemma \ref{reftail}, we have $Q (  (E_\omega
\bar{T}_\nu^{(n)} )^2 > x n^{2/s} , M_1 > n^{(1-\varepsilon)/s}
 )\sim K_\infty x^{-s/2}\frac{1}{n}$, as $n\rightarrow\infty$.
\end{pf*}


\subsection{Existence of random subsequence of nonlocalized behavior}

Introduce the notation:
%
\begin{eqnarray} \label{def}
\mu_{i,n,\omega} &:=& E_\omega^{\nu_{i-1}} \bar{T}_{\nu_i}^{(n)},
\nonumber\\[-8pt]
\\[-8pt]
\sigma_{i,n,\omega}^2 &:=& E_\omega^{\nu_{i-1}}
\bigl( \bar{T}_{\nu_i}^{(n)} - \mu_{i,n,\omega}  \bigr)^2
= \operatorname{Var}_\omega \bigl(\bar{T}_{\nu_i}^{(n)}
- \bar{T}_{\nu_{i-1}}^{(n)} \bigr).
\nonumber
\end{eqnarray}
It is obvious (from the coupling of $\bar{X}_t^{(n)}$ and $X_t$) that
$\mu_{i,n,\omega} \nearrow E_\omega^{\nu_{i-1}} T_{\nu_i}$ as
$n\rightarrow\infty$. It is also true,
although not as obvious, that $\sigma_{i,n,\omega}^2$ is
increasing in $n$ to $\operatorname{Var}_\omega (T_{\nu_i} - T_{\nu_{i-1}})$.
Therefore, we will use the notation $\mu_{i,\infty,\omega}:=
E_\omega ^{\nu_{i-1}} T_{\nu_i}$ and $\sigma_{i,\infty,\omega}^2 :=
\operatorname{Var}_\omega (T_{\nu_i} - T_{\nu_{i-1}}  )$.
To see that $\sigma_{i,n,\omega}^2$ is increasing in $n$, note that the
expansion for $\operatorname{Var}_\omega\bar{T}^{(n)}_\nu$
is the same as the expansion for $\operatorname{Var}_\omega T_\nu$
given in \eqref{VarTnuexpand} but with each
$W_i$ replaced by $W_{\nu_{-b_n}+1, i}$ and with the final sum in
the second line restricted to $\nu_{-b_n} < i < -1$.

The first goal of this subsection is to prove a CLT (along random
subsequences) for the hitting times $T_n$. We begin by showing
that for any $\varepsilon>0$ only\vadjust{\goodbreak}
the crossing times of ladder times with
$M_k > n^{(1-\varepsilon)/s}$ are relevant in the limiting
distribution, at
least along a sparse enough subsequence.

%
\begin{lem}\label{Vsmall}
Assume $s<2$. Then for any $\varepsilon,\delta >0$, there exists an
$\eta>0$ and a sequence $c_n = o(n^{-\eta})$ such that for any $m \leq
\infty$
\begin{eqnarray*}
Q  \Biggl( \sum_{i=1}^n \sigma_{i,m,\omega}^2 \mathbf{1}_{ M_i \leq
n^{(1-\varepsilon)/s}} > \delta n^{2/s}  \Biggr) \leq c_n.
\end{eqnarray*}
\end{lem}
\begin{pf}
Since $\sigma_{i,m,\omega}^2 \leq\sigma_{i,\infty,\omega}^2$, it
is enough to
consider only the case $m=\infty$ (that is, the walk without reflections).
First, we need a bound on the probability of
$\sigma_{i,\infty,\omega} ^2= \operatorname{Var}_\omega(T_{\nu_i}
- T_{\nu_{i-1}})$ being much larger than $M_i^2$. Note that from
\eqref{compare}, we have $\operatorname{Var}_\omega T_\nu\leq(E_\omega
T_\nu)^2 + 8 R_{0,\nu-1}D^-(\omega)$. Then since
$R_{0,\nu-1} \leq\nu M_1$, we have for any $\alpha ,\beta >0$ that
\begin{eqnarray*}
&& Q ( \operatorname{Var}_\omega T_\nu> n^{2\beta }, M_1 \leq n^{\alpha } )
\\
&&\qquad  \leq  Q \biggl( E_\omega T_\nu>
\frac{n^\beta }{\sqrt{2}}, M_1 \leq n^{\alpha }  \biggr)
+ Q \biggl(8 \nu D^-(\omega) > \frac{n^{2\beta -\alpha }}{2}  \biggr).
\end{eqnarray*}
By \eqref{TbigMsmall}, the first term on the right is
$o(e^{-n^{(\beta -\alpha )/5}})$. To bound the second term on the
right, we use Lemma \ref{nutail} and Lemma \ref{piW2} to get that for any
$\alpha < \beta $
\begin{eqnarray*}
Q \biggl( 8 \nu D^-(\omega) > \frac{n^{2\beta -\alpha }}{2}  \biggr)
&\leq & Q(\nu> \log^2 n)
+ Q \biggl( D^-(\omega) > \frac{ n^{2\beta - \alpha}}{16\log^2 n}
 \biggr)
 \\
 & =& o\bigl(n^{-({s}/{2})(3\beta -\alpha )}\bigr).
\end{eqnarray*}
Therefore, similarly to \eqref{TbigMsmall}, we have the bound
%
\begin{eqnarray}\label{VbigMsmall}
Q ( \operatorname{Var}_\omega T_\nu> n^{2\beta }, M_1 \leq n^{\alpha } )
= o\bigl(n^{-({s}/{2})(3\beta -\alpha )}\bigr) .
\end{eqnarray}
The rest of the proof is similar to the proof of Lemma
\ref{smallblocks}. First, from \eqref{VbigMsmall},
\begin{eqnarray*}
&& Q  \Biggl( \sum_{i=1}^n \sigma_{i,\infty,\omega} ^2 \mathbf{1}_{
M_i \leq n^{(1-\varepsilon)/s}} > \delta n^{2/s}  \Biggr)
\\
&&\qquad \leq Q \Biggl( \sum_{i=1}^n \sigma_{i,\infty,\omega} ^2
\mathbf{1}_{\sigma_{i,\infty,\omega} ^2
\leq n^{2(1-{\varepsilon}/{4} )/s}} > \delta n^{2/s}  \Biggr)
\\
&&\qquad \quad {} + n Q \bigl( \operatorname{Var}_\omega T_\nu> n^{2(1-{\varepsilon}/{4})/s},
M_1 \leq n^{(1-\varepsilon)/s}  \bigr)
\\
&&\qquad  = Q \Biggl( \sum_{i=1}^n \sigma_{i,\infty,\omega} ^2
\mathbf{1}_{\sigma_{i,\infty,\omega} ^2 \leq
n^{2(1-{\varepsilon}/{4} )/s}} > \delta n^{2/s}  \Biggr)
+ o(n^{-\varepsilon/8}).
\end{eqnarray*}
%
Therefore, it is enough to prove that for any $\delta ,\varepsilon
>0$, there
exists $\eta>0$ such that
\[
Q \Biggl( \sum_{i=1}^n \sigma_{i,\infty,\omega} ^2
\mathbf{1}_{\sigma_{i,\infty,\omega} ^2 \leq
n^{2(1-{\varepsilon}/{4} )/s}} > \delta n^{2/s}  \Biggr)
= o(n^{-\eta} ).
\]
We prove the above statement by choosing
$C\in(1,\frac{2}{s})$, since $s>2$, and then using Theorem \ref{qVartail}
to get bounds on the size of the set
$ \{i\leq n\dvtx \operatorname{Var}_\omega (T_{\nu_i} - T_{\nu_{i-1}}  )
\in ( n^{2(1-\varepsilon C^k)/s},
n^{2(1-\varepsilon C^{k-1})/s}  ]  \}$ for all $k$ small enough
so that $\varepsilon C^k < 1$. This portion of the proof is similar to
that of Lemma
\ref{smallblocks} and thus will be omitted.
\end{pf}
\begin{cor}\label{Vsdiff}
Assume $s<2$.
Then for any $\delta >0$, there exists an $\eta'>0$ and a sequence
$c_n'=o(n^{-\eta'})$ such that for any $m\leq\infty$
%
\[
Q \Biggl( \Biggl | \sum_{i=1}^n  ( \sigma_{i,m,\omega}^2
-\mu_{i,m,\omega}^2  ) \Biggr |
\geq\delta n^{2/s}  \Biggr) \leq c_n' .
\]
%
\end{cor}
\begin{pf}
For any $\varepsilon> 0$,
%
\begin{eqnarray}
&& Q \Biggl(  \Biggl| \sum_{i=1}^n
( \sigma_{i,m,\omega}^2-\mu_{i,m,\omega}^2 ) \Biggr |
\geq\delta n^{2/s}  \Biggr)
\nonumber\\
&&\qquad \leq  Q \Biggl( \sum_{i=1}^n \sigma_{i,m,\omega}^2
\mathbf{1}_{M_i \leq n^{(1-\varepsilon)/s}}
\geq\frac{\delta}{3}n^{2/s}  \Biggr)
\label{smM1}\\
&&\qquad \quad {} + Q \Biggl( \sum_{i=1}^n \mu_{i,m,\omega}^2\mathbf{1}_{M_i
\leq n^{(1-\varepsilon)/s}} \geq\frac{\delta }{3}n^{2/s}  \Biggr)
\label{smM2}\\
&&\qquad \quad {} + Q \Biggl( \sum_{i=1}^n  | \sigma_{i,m,\omega}^2 -
\mu_{i,m,\omega}^2  | \mathbf{1}_{M_i > n^{(1-\varepsilon)/s}}
\geq \frac{\delta }{3}n^{2/s}  \Biggr).
\label{lgM}
\end{eqnarray}
Lemma \ref{Vsmall} gives that \eqref{smM1} decreases
polynomially in $n$ (with a bound not depending on $m$).
Also, essentially the same proof as in Lemmas \ref{Vsmall} and
\ref{smallblocks} can be used to show that \eqref{smM2} also decreases
polynomially in $n$ (again with a bound not depending on $m$).
Finally, \eqref{lgM} is bounded above by
\begin{eqnarray*}
&& Q \bigl( \#  \bigl\{ i\leq n\dvtx M_i > n^{(1-\varepsilon)/s} \bigr\}
> n^{2\varepsilon} \bigr)
\\
&&\qquad {}
+ n Q \biggl( \bigl |\operatorname{Var}_\omega\bar{T}_\nu^{(m)}
- \bigl(E_\omega\bar{T}_\nu^{(m)}\bigr)^2  \bigr|
\geq\frac{\delta }{3} n^{2/s -2\varepsilon} \biggr) ,
\end{eqnarray*}
and since by \eqref{Mtail}, $Q ( \#  \{ i\leq n\dvtx M_i >
n^{(1-\varepsilon)/s}  \} > n^{2\varepsilon}  ) \leq \frac{ n Q(M_1 >
n^{(1-\varepsilon)/s}) }{n^{2\varepsilon}} \sim C_5 n^{-\varepsilon}$
we need only show that
for some $\varepsilon>0$ the second term above is decreasing faster
than a power of $n$.
However, from \eqref{comparer}, we have $ |\operatorname{Var}_\omega
\bar{T}_\nu^{(m)} - (E_\omega\bar{T}_\nu^{(m)})^2  | \leq
D^+(\omega) + 8 R_{0,\nu-1} D^-(\omega)$. Thus,
\begin{eqnarray*}
&& n Q \biggl( \bigl |\operatorname{Var}_\omega\bar{T}_\nu^{(m)}
- \bigl(E_\omega\bar{T}_\nu^{(m)}\bigr)^2 \bigr |
\geq\frac{\delta }{3} n^{2/s - 2\varepsilon}  \biggr)
\\
&&\qquad
\leq n Q \biggl( D^+(\omega) + 8 R_{0,\nu-1} D^-(\omega)>
\frac{\delta}{3} n^{2/s - 2\varepsilon}  \biggr),
\end{eqnarray*}
and for any $\varepsilon<\frac{1}{2s}$ Lemma \ref{abound} and
Corollary \ref{bbound} give that the last term above decreases faster
than some power of $n$.
%
%
\end{pf}

Since $T_{\nu_n} = \sum_{i=1}^n (T_{\nu_i} - T_{\nu_{i-1}})$
is the sum of independent (quenched) random variables, in order to
prove a CLT we cannot have any of the first $n$ crossing times of
blocks dominating all the others (note this is exactly what
happens in the localization behavior we saw in Section
\ref{localization}). Thus, we look for a random subsequence where
none of the crossing times of blocks are dominant.
Now, for any $\delta \in(0,1]$ and any positive integer $a < n/2$,
define the event
\[
\mathcal{S}_{\delta ,n,a} :=  \bigl\{ \# \{ i\leq\delta n\dvtx
\mu_{i,n,\omega}^2 \in[n^{2/s}, 2n^{2/s})  \} = 2a,
\mu_{j,n,\omega}^2 < 2 n^{2/s} \ \forall j\leq\delta n  \bigr\}.
\]
On the event $\mathcal{S}_{\delta , n, a}$, $2a$ of the first $\delta
n$ crossings times from $\nu_{i-1}$ to $\nu_i$ have roughly the
same size expected crossing times $\mu_{i,n,\omega}$, and the rest are
all smaller (we work with $\mu_{i,n,\omega}^2$ instead of
$\mu_{i,n,\omega}$ so that comparisons with $\sigma_{i,n,\omega}^2$
are slightly easier). We want a lower bound on the probability of
$\mathcal{S}_{\delta ,n,a}$. The difficulty in getting a lower bound is
that the $\mu^2_{i,n,\omega}$ are not independent. However, we can
force all the large crossing times to be independent by forcing
them to be separated by at least $b_n$ ladder locations.

Let $\mathcal{I}_{\delta ,n,a}$ be the collection of all subsets $I$ of
$[1,\delta n]\cap\mathbb{Z}$ of size $2 a$ with the property that
any two
distinct points in $I$ are separated by at least $2b_n$. Also,
define the event
\[
A_{i,n}:=  \{ \mu_{i,n,\omega}^2 \in [ n^{2/s},2 n^{2/s} )  \}.
\]
Then we begin with a simple lower bound:
%
\begin{eqnarray}\label{simlb}
Q( \mathcal{S}_{\delta ,n,a} ) &\geq& Q \Biggl( \bigcup_{ I\in
\mathcal{I}_{\delta ,n,a}}  \Biggl( \bigcap_{i\in I} A_{i,n}
\bigcap_{j\in[1,\delta n]\backslash I}  \{ \mu_{j,n,\omega}^2 <
n^{2/s}  \}  \Biggr)  \Biggr)
\nonumber\\[-8pt]
\\[-8pt]
& =& \sum_{I\in\mathcal{I}_{\delta ,n,a}} Q\Biggl ( \bigcap_{i\in I}
A_{i,n} \bigcap_{j\in[1,\delta n]\backslash I}  \{ \mu
_{j,n,\omega}^2 < n^{2/s}  \}  \Biggr).
\nonumber
\end{eqnarray}
Now, recall the definition of the event $G_{i,n,\varepsilon}$ from
\eqref{FGevents}, and define the event
\[
H_{i,n,\varepsilon} :=  \bigl\{ M_j \leq n^{(1-\varepsilon)/s}
\mbox{ for all } j\in [i-b_n,i) \bigr\}.
\]
Also, for any $I\subset\mathbb{Z}$ let $d(j,I):= \min\{ |j-i| : i\in I
\}$ be the minimum distance from $j$ to the set $I$. Then with
minimal cost, we can assume that for any $I\in
\mathcal{I}_{\delta ,n,a}$ and any $\varepsilon>0$ that all $j\notin
I$ such
that $d(j,I)\leq b_n$ have $M_j \leq n^{(1-\varepsilon)/s}$. Indeed,
%
\begin{eqnarray}\label{Ilb}
\qquad
&& Q \Biggl( \bigcap_{i\in I} A_{i,n}
\bigcap_{j\in [1,\delta n]\backslash I}
\{ \mu_{j,n,\omega}^2 < n^{2/s}  \}  \Biggr)
\nonumber\\
&&\qquad \geq  Q \Biggl( \bigcap_{i\in I}  ( A_{i,n} \cap
G_{i,n,\varepsilon} \cap
H_{i,n,\varepsilon}  ) \bigcap_{j\in[1,\delta n]:d(j,I)>b_n}
\{ \mu_{j,n,\omega}^2 < n^{2/s}  \}  \Biggr)
\nonumber\\
&&\qquad \quad {} - Q \Biggl(\bigcup_{j\notin{I}, d(j,I)\leq b_n}
\bigl\{\mu_{j,n,\omega}^2 \geq n^{2/s} ,
M_j \leq n^{(1-\varepsilon)/s}  \bigr\}  \Biggr)
\\
&&\qquad \geq \prod_{i\in I} Q(A_{i,n} \cap H_{i,n,\varepsilon} )
Q \Biggl(\bigcap_{i\in I} G_{i,n,\varepsilon} \bigcap_{j\in[1,\delta
n]:d(j,I)>b_n}  \{\mu_{j,n,\omega}^2 < n^{2/s}  \}  \Biggr)
\nonumber\\
&&\qquad \quad {}- 4 a b_n Q \bigl( E_\omega T_\nu\geq n^{1/s}, M_1 \leq
n^{(1-\varepsilon)/s}  \bigr).
\nonumber
\end{eqnarray}
From Theorem \ref{Tnutail} and Lemma \ref{reftail}, we have
$Q(A_{i,n}) \sim K_\infty(1-2^{-s/2}) n^{-1}$. We wish to show
the same asymptotics are true for $Q(A_{i,n} \cap H_{i,n,\varepsilon})$
as well. From \eqref{Mtail}, we have $Q(H_{i,n,\varepsilon}^c) \leq b_n
Q(M_1 > n^{(1-\varepsilon)/s}) = o(n^{-1+2\varepsilon})$. Applying this,
along with
\eqref{Mtail} and \eqref{TbigMsmall}, gives that for $\varepsilon>0$,
\begin{eqnarray*}
Q(A_{i,n}) &\leq & Q( A_{i,n} \cap H_{i,n,\varepsilon} )
+ Q \bigl(M_1 >n^{(1-\varepsilon)/s} \bigr)
Q(H_{i,n,\varepsilon}^c)
\\
&&{} + Q \bigl(E_\omega T_\nu> n^{1/s},
M_1 \leq n^{(1-\varepsilon)/s}  \bigr)
\\
& =& Q( A_{i,n} \cap H_{i,n,\varepsilon})
+ o(n^{-2+3\varepsilon}) + o\bigl(e^{-n^{\varepsilon/(5s)}}\bigr).
\end{eqnarray*}
Thus, for any $\varepsilon<\frac{1}{3}$, there exists a
$C_\varepsilon> 0$ such that
%
\begin{equation}\label{AHlb}
Q(A_{i,n}\cap H_{i,n,\varepsilon}) \geq C_\varepsilon n^{-1}.
\end{equation}
To handle the next probability in \eqref{Ilb}, note that
%
\begin{eqnarray}\label{Glb}
&& Q \Biggl( \bigcap_{i\in I} G_{i,n,\varepsilon} \bigcap_{j\in
[1,\delta n]:d(j,I)>b_n}  \{ \mu_{j,n,\omega}^2 < n^{2/s}  \}  \Biggr)
\nonumber\\
&&\qquad \geq  Q \Biggl(\bigcap_{j\in[1,\delta n]}  \{ \mu_{j,n,\omega}^2
< n^{2/s}  \}  \Biggr) - Q \Biggl(
\bigcup_{i\in I} G_{i,n,\varepsilon}^c  \Biggr)
\nonumber\\[-8pt]
\\[-8pt]
&&\qquad \geq Q ( E_\omega T_{\nu_n} < n^{1/s}  )
- 2 a Q(G_{i,n,\varepsilon}^c )
\nonumber\\
&&\qquad = Q ( E_\omega T_{\nu_n} < n^{1/s}  ) - a
o(n^{-1+2\varepsilon}).
\nonumber
\end{eqnarray}
Finally, from \eqref{TbigMsmall}, we have\vadjust{\goodbreak} $ 4 a b_n Q ( E_\omega
T_\nu\geq n^{1/s}, M_1 \leq n^{(1-\varepsilon)/s}  ) = a o (
e^{-n^{\varepsilon/(6 s)}}  ) $. This along with \eqref{AHlb} and
\eqref{Glb} applied to \eqref{simlb} gives
\begin{eqnarray*}
&& Q ( \mathcal{S}_{\delta ,n,a}  )
\\
&&\qquad \geq\#(\mathcal
{I}_{\delta ,n,a}) \bigl[ ( C_\varepsilon n^{-1}  )^{2a}
\bigl( Q (E_\omega T_{\nu_n}
< n^{1/s}  ) - a o(n^{-1+2\varepsilon})  \bigr)
- a o \bigl(e^{-n^{\varepsilon/(6 s)}}  \bigr)  \bigr].
\end{eqnarray*}
An obvious upper bound for $\#(\mathcal{I}_{\delta ,n,a})$ is
${\delta n \choose 2a} \leq\frac{(\delta n)^{2a}}{(2a)!}$. To get a lower
bound on $\#(\mathcal{I}_{\delta ,n,a})$, we note that any set
$I\in\mathcal{I}_{\delta ,n,a}$ can be chosen in the following way:
first choose an integer $i_1 \in[1,\delta n]$ ($\delta n$ ways to do
this). Then choose an integer $i_2 \in[1,\delta n] \backslash\{
j\in\mathbb{Z}\dvtx |j-i_1| \leq2b_n \}$ (at least $\delta n - 1 -4b_n$ ways to
do this). Continue this process until $2a$ integers have been
chosen. When choosing $i_j$, there will be at least $\delta n -
(j-1)(1+4b_n)$ integers available. Then since there are $(2a)!$
orders in which to choose each set if $2a$ integers, we have
\begin{eqnarray*}
\frac{(\delta n)^{2a}}{(2a)!} \geq\#(\mathcal{I}_{\delta ,n,a})
&\geq & \frac{1}{(2a)!} \prod_{j=1}^{2a}
\bigl(\delta n - (j-1)(1+4b_n) \bigr)
\\
&\geq&\frac{(\delta n)^{2a}}{(2a)!}
\biggl( 1 - \frac{(2a-1)(1+4b_n)}{\delta n}  \biggr)^{2a}.
\end{eqnarray*}
Therefore, applying the upper and lower bounds on
$\#(\mathcal{I}_{\delta ,n,a})$, we get
\begin{eqnarray*}
Q ( \mathcal{S}_{\delta ,n,a}  )
& \geq & \frac{(\delta C_\varepsilon)^{2a}}{(2a)!}
\biggl( 1 - \frac{(2a-1)(1+4b_n)}{\delta n} \biggr)^{2a}
\\
&&{}\times
\bigl( Q ( E_\omega T_{\nu_n} < n^{1/s}  )
- a o(n^{-1+2\varepsilon})  \bigr)
\\
&&{} - \frac{(\delta n)^{2a}}{(2a)!} a o
\bigl( e^{-n^{\varepsilon/(6 s)}} \bigr).
\end{eqnarray*}
Recall the definitions of $d_k$ in \eqref{nkdkdef} and define
%
\begin{equation}\label{adeltadef}
a_k:=\lfloor\log\log k\rfloor\vee1
\quad\mbox{and}\quad \delta _k:=a_k^{-1}.
\end{equation}
Now, replacing $\delta ,n$ and $a$ in the above by $\delta _k, d_k$ and
$a_k$, respectively, we have
%
\begin{eqnarray}\label{sublb}
Q ( \mathcal{S}_{\delta _k,d_k,a_k}  )
& \geq& \frac{(\delta _kC_\varepsilon)^{2a_k}}{(2a_k)!}
\biggl( 1 - \frac{(2a_k-1)(1+4b_{d_k})}{\delta _kd_k}  \biggr)^{2a_k}
\nonumber\\
&&{}\times
\bigl ( Q ( E_\omega T_{\nu_{d_k}} < d_k^{1/s} )
- a_ko(d_k^{-1+2\varepsilon})  \bigr)
\nonumber\\[-8pt]
\\[-8pt]
&&{} - \frac{(\delta _k d_k)^{2a_k}}{(2a_k)!} a_k o
\bigl(e^{-d_k^{\varepsilon/(6 s)}}  \bigr)
\nonumber\\
& \geq& \frac{(\delta _k C_\varepsilon)^{2a_k}}{(2a_k)!}
\bigl( 1+o(1) \bigr)\bigl ( L_{s,b'}(1) - o(1)  \bigr) - o(1/k).
\nonumber
\end{eqnarray}
The last inequality is a result of the definitions of $\delta _k, a_k$,
and $d_k$ (it's enough to recall that $d_k \geq2^{2^{k-1}}$, $a_k
\sim\log\log k$, and $\delta _k \sim\frac{1}{\log\log k} $), as well
as Theorem \ref{refstable}. Also, since $\delta _k = a_k^{-1}$, we get
from Stirling's formula that $ \frac{(\delta _k C_\varepsilon)^{2a_k}}{(2a_k)!}
\sim\frac{ ( C_\varepsilon e / 2)^{2a_k}}{\sqrt{2\pi a_k}}$. Thus, since
$a_k\sim\log\log k$, we have that $\frac{1}{k} = o (
\frac{(\delta _k C_\varepsilon)^{2a_k}}{(2a_k)!}  )$. This,
along with \eqref{sublb}, gives that $Q ( \mathcal{S}_{\delta _k,d_k,a_k}
 ) > \frac{1}{k}$ for all $k$ large enough.

We now have a good lower bound on the probability of not having
any of the crossing times of the first $\delta _k d_k$ blocks dominating
all the others. However, for the purpose of proving Theorem
\ref{nonlocal}, we need a little bit more. We also need that none
of the crossing times of succeeding blocks are too large either.
Thus, for any $0<\delta <c$ and $n\in\mathbb{N}$, define the events
\[
U_{\delta ,n,c}:= \Biggl\{ \sum_{i=\delta n + 1}^{cn} \mu
_{i,n,\omega} \leq2 n^{1/s}  \Biggr\},
\qquad
\tilde{U}_{\delta ,n,c}:= \Biggl\{ \sum_{i=\delta n
+ b_n+1}^{cn} \mu_{i,n,\omega} \leq n^{1/s}  \Biggr\}.
\]
\begin{lem} \label{smallblocklemma}
Assume $s<1$. Then there exists a sequence $c_k \rightarrow\infty$,
$c_k = o(\log a_k )$ such that
\[
\sum_{k=1}^{\infty} Q ( \mathcal{S}_{\delta _k,d_k,a_k} \cap
U_{\delta _k,d_k,c_k}  ) = \infty.
\]
\end{lem}
\begin{pf}
For any $\delta <c$ and $a<n/2$, we have
%
\begin{eqnarray}\label{SU}
\qquad
Q ( \mathcal{S}_{\delta ,n,a} \cap U_{\delta ,n,c}  )
&\geq & Q ( \mathcal{S}_{\delta ,n,a} ) Q ( \tilde{U}_{\delta ,n,c} )
- Q \Biggl( \sum_{i=1}^{b_n} \mu_{i,n,\omega}
> n^{1/s}  \Biggr)
\nonumber\\
& \geq & Q ( \mathcal{S}_{\delta ,n,a} )
Q ( E_\omega T_{\nu_{cn}} \leq n^{1/s} )
- b_n Q \biggl( E_\omega T_\nu>\frac{n^{1/s}}{ b_n } \biggr)
\\
&\geq & Q ( \mathcal{S}_{\delta ,n,a} ) Q ( E_\omega
T_{\nu_{cn}} \leq n^{1/s}  ) - o(n^{-1/2}),
\nonumber
\end{eqnarray}
where the last inequality is from Theorem \ref{Tnutail}. Now,
define $c_1=1$ and for $k>1$ let
\[
c_k' := \max \biggl\{ c \in\mathbb{N}\dvtx
Q ( E_\omega T_{\nu_{c d_k}} \leq
d_k^{1/s}  ) \geq\frac{1}{\log k}  \biggr\} \vee1.
\]
Note that by Theorem \ref{refstable} we have that $c_k'\rightarrow
\infty$,
and so we can define $c_k = c_k' \wedge\log\log(a_k)$. Then
applying \eqref{SU} with this choice of $c_k$ we have
\begin{eqnarray*}
&& \sum_{k=1}^\infty Q ( \mathcal{S}_{\delta _k,d_k,a_k} \cap
U_{\delta _k,d_k,c_k}  )
\\
&&\qquad \geq\sum_{k=1}^\infty [Q ( \mathcal{S}_{\delta _k,d_k,a_k} )
Q ( E_\omega T_{\nu_{c_k d_k}}
\leq d_k^{1/s}  ) - o(d_k^{-1/2}) ] = \infty,
\end{eqnarray*}
and the last
sum is infinite because $d_k^{-1/2}$ is summable and for
all $k$ large enough we have
\[
Q ( \mathcal{S}_{\delta _k,d_k,a_k} ) Q ( E_\omega
T_{\nu_{c_k
d_k}} \leq d_k^{1/s}  ) \geq\frac{1}{k\log k}.
\]\upqed
\end{pf}
\begin{cor}\label{smallblocks2}
Assume $s<1$, and let $c_k$ be as in Lemma \ref{smallblocklemma}. Then
$P$-a.s. there exists a random subsequence $n_{k_m}=n_{k_m}(\omega)$
of $n_k=2^{2^k}$ such that for the sequences $\alpha _m, \beta _m$, and
$\gamma_m$ defined by
%
\begin{eqnarray}\label{abgdef}
\alpha _m &:=& n_{k_m-1},
\nonumber\\
\beta _m &:=& n_{k_m-1} +\delta _{k_m} d_{k_m} ,
\\
\gamma_m&:=& n_{k_m-1} + c_{k_m} d_{k_m},
\nonumber
\end{eqnarray}
we have that for all $m$
%
\begin{equation}\label{smbk}
\max_{ i\in(\alpha _m, \beta _m] } \mu_{i,d_{k_m},\omega}^2
\leq 2d_{k_m}^{2/s} \leq\frac{1}{a_{k_m}}
\sum_{i=\alpha _m+1}^{\beta _m} \mu_{i,d_{k_m},\omega}^2
\end{equation}
and
\[
\sum_{\beta _m + 1}^{\gamma_m} \mu_{i,d_{k_m},\omega}
\leq2 d_{k_m}^{1/s}.
\]
\end{cor}
\begin{pf}
Define the events
\begin{eqnarray*}
\mathcal{S}_k' &:=&  \bigl\{ \# \{ i \in(n_{k-1}, n_{k-1}+\delta _k
d_k] \dvtx \mu_{i,d_k,\omega}^2 \in[d_k^{2/s}, 2d_k^{2/s})  \} = 2a_k\bigr\}
\\
&&\hspace*{0pt}{} \cap \{ \mu_{j,d_k,\omega}^2 < 2 d_k^{2/s}
\ \forall j\in(n_{k-1}, n_{k-1}+\delta _k d_k]  \},
\\
U_k' &:= & \Biggl\{ \sum_{n_{k-1}+ \delta _k d_k + 1}^{n_{k-1}+c_k d_k}
\mu_{i,d_{k},\omega} \leq2d_{k}^{1/s}  \Biggr\}.
\end{eqnarray*}
Note that due to the reflections of the random walk, the event
$\mathcal{S}_k' \cap U_k' $ depends on the environment between
ladder locations $n_{k-1}-b_{d_k}$ and $n_{k-1}+c_kd_k$. Thus, for
$k_0$ large enough $\{\mathcal{S}_{2k}' \cap U_{2k}'
\}_{k=k_0}^{\infty}$ is an independent sequence of events.
Similarly, for~$k$ large enough $\mathcal{S}_k'\cap U_k'$ does not
depend on the environment to left of the origin. Thus,
\[
P(\mathcal{S}_k'\cap U_k') =Q(\mathcal{S}_k'\cap U_k')
= Q (\mathcal{S}_{\delta _k,d_k,a_k} \cap U_{\delta _k,d_k,c_k} )
\]
for all $k$ large enough. Lemma \ref{smallblocklemma} then gives
that $\sum_{k=1}^\infty P(\mathcal{S}_{2k}'\cap U_{2k}') =
\infty$, and the Borel--Cantelli lemma then implies that infinitely
many of the events $\mathcal{S}_{2k}' \cap U_{2k}'$ occur $P$-a.s.
Finally, note that $\mathcal{S}_{k_m}'$ implies the event in
\eqref{smbk}.
\end{pf}
%

Before proving a quenched CLT (along a subsequence) for the hitting times
$T_n$, we need one more lemma that gives us some control on the
quenched tails of crossing times of blocks. We can get this from
an application of Kac's moment formula. Let $\bar{T}_y$ be the
hitting time of $y$ when we add a reflection at the starting point
of the random walk. Then Kac's moment formula (\cite{fpKMF}, (6))
and the Markov property give that
$E_\omega^x (\bar{T}_y)^j \leq j!  ( E_\omega^x
\bar{T}_y  )^j$ (note that because
of the reflection at $x$, $E_\omega^x(\bar{T}_y)
\geq E_\omega^{x'}(\bar{T}_y)$ for any $x'\in(x,y)$). Thus,
%
\begin{eqnarray}\label{Kac}
E_\omega^{\nu_{i-1}} \bigl(\bar{T}_{\nu_i}^{(n)}\bigr)^j
&\leq &E_\omega^{\nu_{i-1-b_n}} (\bar{T}_{\nu_i} )^j
\leq j!  (E_\omega^{\nu_{i-1-b_n}} \bar{T}_{\nu_i}  )^j
\nonumber\\[-8pt]
\\[-8pt]
&\leq &j!  (E_\omega^{\nu_{i-1-b_n}} \bar{T}_{\nu_{i-1}} + \mu_{i,n,\omega} )^j.
\nonumber
\end{eqnarray}
\begin{lem} \label{momentbound}
For any $\varepsilon< \frac{1}{3}$, there exists an $\eta>0$ such that
\[
Q \bigl( \exists i\leq n, j\in\mathbb{N}\dvtx M_i > n^{(1-\varepsilon)/s},
E_\omega^{\nu_{i-1}} \bigl( \bar{T}_{\nu_i}^{(n)} \bigr)^j
> j! 2^j \mu_{i,n,\omega}^j  \bigr) = o(n^{-\eta}).
\]
\end{lem}
\begin{pf}
We use \eqref{Kac} to get
\begin{eqnarray*}
&& Q \bigl( \exists i\leq n, j\in\mathbb{N}\dvtx
M_i > n^{(1-\varepsilon)/s},
E_\omega^{\nu_{i-1}} \bigl( \bar{T}_{\nu_i}^{(n)} \bigr)^j
> j! 2^j \mu_{i,n,\omega}^j  \bigr)
\\
&&\qquad \leq Q \bigl( \exists i\leq n \dvtx M_i > n^{(1-\varepsilon)/s},
E_\omega^{\nu_{i-1-b_n}} \bar{T}_{\nu_{i-1}} > \mu_{i,n,\omega}  \bigr)
\\
&&\qquad \leq n Q \bigl( M_1 > n^{(1-\varepsilon)/s} ,
E_\omega^{\nu_{-b_n}} T_0> n^{(1-\varepsilon)/s}  \bigr)
\\
&&\qquad = n Q \bigl( M_1 > n^{(1-\varepsilon)/s}  \bigr)
Q \bigl( E_\omega^{\nu_{-b_n}}T_0 > n^{(1-\varepsilon)/s} \bigr),
\end{eqnarray*}
where the second inequality is due to a union bound and the fact
that $\mu_{i,n,\omega} > M_i$. Now, by \eqref{Mtail}, we have
$nQ ( M_1 > n^{(1-\varepsilon)/s}  ) \sim C_5 n^{\varepsilon}$, and by
Theorem \ref{Tnutail},
\[
Q \bigl( E_\omega^{\nu_{-b_n}} T_0 > n^{(1-\varepsilon)/s} \bigr)
\leq b_n
Q \biggl( E_\omega T_\nu> \frac{n^{(1-\varepsilon)/s}}{b_n}  \biggr)
\sim K_\infty b_n^{1+s} n^{-1+\varepsilon}.
\]
Therefore, $Q ( \exists i\leq n, j\in\mathbb{N}\dvtx M_i >
n^{(1-\varepsilon)/s}, E_\omega^{\nu_{i-1}} ( \bar{T}_{\nu
_i}^{(n)} )^j > j! 2^j \mu_{i,n,\omega}^j  ) =\break o(n^{-1+3\varepsilon})$.
\end{pf}
\begin{thm} \label{gaussianT}
Let Assumptions \ref{essentialasm} and \ref{techasm} hold, and let
$s<1$. Then $P$-a.s. there exists a random subsequence
$n_{k_m}=n_{k_m}(\omega)$
of $n_k=2^{2^k}$ such that for $\alpha _m$, $\beta _m$ and $\gamma
_m$ as in \eqref{abgdef} and any sequence $x_m \in[\nu_{\beta _m},
\nu_{\gamma_m} ]$, we have
%
\begin{equation}
\lim_{m\rightarrow\infty} P_\omega \biggl( \frac{ T_{x_m} -
E_\omega T_{x_m} } {\sqrt{v_{m,\omega}} } \leq y  \biggr) = \Phi(y),
\end{equation}
where
\[
v_{m,\omega} := \sum_{i=\alpha _m+1}^{\beta _m} \mu_{i,d_{k_m},
\omega}^2 .
\]
\end{thm}
\begin{pf}
Let $n_{k_m}(\omega)$ be the random subsequence specified in Corollary
\ref{smallblocks2}. For ease of notation, set
$\tilde{a}_m = a_{k_m}$ and $ \tilde{d}_m = d_{k_m}$.
We have
\[
\max_{i\in(\alpha _m, \beta _m]} \mu_{i,\tilde{d}_m,\omega}^2
\leq 2 \tilde{d}_m^{2/s} \leq
\frac{1}{\tilde{a}_m}\sum_{i=\alpha _m+1}^{\beta _m}
\mu_{i,\tilde{d}_m,\omega}^2 = \frac{v_{m,\omega}}{\tilde{a}_m}
\]
and
\[
\sum_{i=\beta _m+1}^{\gamma_m} \mu_{i,\tilde{d}_m,\omega}
\leq2 \tilde{d}_m^{1/s}.
\]
Now, let $\{ x_m \}_{m=1}^\infty$ be any sequence of integers
(even depending on $\omega$) such that $x_m \in[\nu_{\beta _m},
\nu_{\gamma_m}]$. Then since $(T_{x_m} - E_\omega T_{x_m}) =
(T_{\nu_{\alpha _m}} - E_\omega T_{\nu_{\alpha _m}}) +
(T_{x_m} - T_{\nu_{\alpha _m}} - E_\omega^{\nu_{\alpha _m}}
T_{x_m})$,
it is enough to prove
%
\begin{equation}\label{irrelevantb}
\qquad
\frac{T_{\nu_{\alpha _m}} - E_\omega T_{\nu_{\alpha _m}}}
{\sqrt{v_{m,\omega}}}\stackrel{\mathcal{D}_\omega}{\longrightarrow}0
\quad\mbox{and}\quad
\frac{T_{x_m} - T_{\nu_{\alpha _m}} -
E_\omega^{\nu_{\alpha _m}} T_{x_m}}{\sqrt{v_{m,\omega}}}
\stackrel{\mathcal{D}_\omega}{\longrightarrow}Z\sim N(0,1),
\hspace*{-5pt}
\end{equation}
where we use the notation $Z_n \stackrel{\mathcal{D}_\omega
}{\longrightarrow}Z$ to denote quenched
convergence in distribution, that is $\lim_{n\rightarrow\infty}
P_\omega( Z_n\leq z ) = P_\omega( Z \leq z )$, $P$-a.s. For the first term in
\eqref{irrelevantb}, note that for any $\varepsilon> 0$, we have from
Chebyshev's inequality and $v_{m,\omega} \geq
\tilde{d}_m^{2/s}$ that
\[
P_\omega \biggl(  \biggl| \frac{T_{\nu_{\alpha _m}} - E_\omega
T_{\nu_{\alpha _m}}}{\sqrt{v_{m,\omega}}}  \biggr| \geq\varepsilon
 \biggr) \leq
\frac{\operatorname{Var}_\omega T_{\nu_{\alpha _m}}}
{\varepsilon^2 v_{m,\omega} } \leq\frac
{\operatorname{Var}_\omega T_{\nu_{\alpha _m}}}{\varepsilon^2 \tilde{d}_m^{2/s}}.
\]
Thus, the first claim in \eqref{irrelevantb} will be proved if we can
show that
$\operatorname{Var}_\omega T_{\nu_{\alpha _m}} = o(\tilde{d}_m^{2/s})$.
For this, we need the following lemma.
\begin{lem} \label{Varasymp}
Assume $s\leq2$. Then for any $\delta > 0$,
\[
P \bigl( \operatorname{Var}_\omega T_{\nu_n} \geq n^{2/s + \delta }  \bigr)
= o(n^{-\delta s/4} ).
\]
\end{lem}
\begin{pf}
First, we claim that
%
\begin{equation}\label{eq-100407a}
E_P (\operatorname{Var}_\omega T_1)^\gamma< \infty  \mbox{ for any }
\gamma< \frac{s}{2}.
\end{equation}
Indeed, from
\eqref{qvar}, we have that for any $\gamma< \frac{s}{2} \leq1$,
\begin{eqnarray*}
E_P (\operatorname{Var}_\omega T_1)^\gamma
&\leq& 4^\gamma E_P( W_0 + W_0^2 )^\gamma
+ 8^\gamma\sum_{i<0} E_P
\bigl ( \Pi_{i+1,0}^\gamma(W_i+W_i^{2})^\gamma \bigr)
\\
&=& 4^\gamma E_P( W_0 + W_0^2 )^\gamma+ 8^\gamma
\sum_{i=1}^{\infty} (E_P \rho_0^\gamma)^i E_P ( W_0 + W_0^2)^\gamma,
\end{eqnarray*}
where we used that $P$ is i.i.d. in the last equality.
Since $E_P \rho_0^\gamma< 1$ for
any $\gamma\in(0,s)$, we have that \eqref{eq-100407a} follows
as soon as $E_P(W_0+W_0^2)^\gamma< \infty$. However,
from \eqref{PWtail}, we
get that $E_P(W_0+W_0^2)^\gamma< \infty$
when $\gamma< \frac{s}{2}$.

As in Lemma \ref{mdevu} let $\bar\nu= E_P \nu$. Then,
\[
P ( \operatorname{Var}_\omega T_{\nu_n} \geq n^{2/s + \delta }  )
\leq P(\operatorname{Var}_\omega
T_{2\bar{\nu} n} \geq n^{2/s+\delta } ) + P(\nu_n \geq2\bar{\nu}n).
\]
As in Lemma \ref{mdevu}, the second term is $\mathcal{O}(e^{-\delta ' n} )$
for some $\delta '>0$. To handle the first term on the right side, we
note that for any $\gamma< \frac{s}{2} \leq1$,
%
\begin{eqnarray}\label{eq-100407b}
P( \operatorname{Var}_\omega T_{2\bar{\nu} n} \geq n^{2/s+\delta } )
&\leq& \frac{E_P (\sum_{k=1}^{2\bar{\nu} n}
\operatorname{Var}_\omega(T_k-T_{k-1})  )^\gamma}{n^{\gamma(2/s + \delta )}}
\nonumber\\[-8pt]
\\[-8pt]
&\leq& \frac{2\bar{\nu}n E_P(\operatorname{Var}_\omega T_1)^\gamma}
{n^{\gamma(2/s+\delta )}}.
\nonumber
\end{eqnarray}
Then since $E_P (\operatorname{Var}_\omega T_1)^\gamma< \infty$ for any $\gamma<
\frac{s}{2}$, we can choose $\gamma$ arbitrarily close to
$\frac{s}{2}$ so that the last term on the right of \eqref{eq-100407b}
is $o(n^{-\delta s/4})$.
\end{pf}

As a result of Lemma \ref{Varasymp} and the Borel--Cantelli
lemma, we have that $\operatorname{Var}_\omega T_{\nu_{n_k}} = o(n_k^{2/s+\delta })$ for
any $\delta > 0$. Therefore, for any $\delta \in(0,\frac{2}{s})$,
we have $\operatorname{Var}_\omega T_{\nu_{\alpha _m}} = o( \alpha _m^{2/s
+\delta }) = o(n_{k_m-1}^{2/s+\delta }) = o(\tilde{d}_m^{2/s})$
(in the last equality we use that $d_k\sim n_k$ to grow much faster
than exponentially in $k$).

For the next step in the proof, we show that reflections can be
added without changing the limiting distribution. Specifically, we show
that it is enough to prove the following lemma, whose proof we postpone.
\begin{lem}\label{lem-latenight}
With notation as in Theorem \ref{gaussianT}, we have
%
\begin{equation}\label{refgaussian}
\lim_{m\rightarrow\infty} P_\omega^{\nu_{\alpha _m}}
\biggl(\frac{\bar{T}^{(\tilde{d}_m)}_{x_m} - E_\omega^{\nu_{\alpha_m}}
\bar{T}^{(\tilde{d}_m)}_{x_m}} { \sqrt{v_{m,\omega}} } \leq y  \biggr) = \Phi(y).
\end{equation}
\end{lem}

Assuming Lemma \ref{lem-latenight}, we complete the
proof of Theorem \ref{gaussianT}. It is enough to show that
%
\begin{equation}\label{eq-maybelast}
\lim_{m\rightarrow\infty} P_\omega^{\nu_{\alpha _m}}
\bigl(\bar{T}^{(\tilde{d}_m)}_{x_{k_m}} \neq T_{x_m} \bigr) = 0
\quad\mbox{and}\quad
\lim_{m\rightarrow\infty} E_\omega^{\nu_{\alpha _m}}
\bigl(T_{x_m}- \bar{T}^{(\tilde{d}_m)}_{x_{k_m}} \bigr) = 0.
\end{equation}
Recall that the coupling introduced after \eqref{bdef} gives that $
T_{x_m} - \bar{T}^{(\tilde{d}_m)}_{x_m}
\geq0$. Thus,
\[
P_\omega^{\nu_{\alpha _m}}\bigl( \bar{T}^{(\tilde{d}_m)}_{x_m}
\neq T_{x_m} \bigr) = P_\omega^{\nu_{\alpha _m}} \bigl( T_{x_m} -
\bar{T}^{(\tilde{d}_m)}_{x_m} \geq1  \bigr) \leq E_\omega^{\nu
_{\alpha_m}} \bigl(T_{x_m} - \bar{T}^{(\tilde{d}_m)}_{x_m} \bigr).
\]
Then since
$x_m \leq\nu_{\gamma_m}$ and $
\gamma_m=n_{k_m-1}+c_{k_m}\tilde{d}_m \leq n_{k_m+1}$ for all $m$
large enough, \eqref{eq-maybelast} will follow from
%
\begin{equation}\label{ETETbar}
\lim_{k\rightarrow\infty} E_\omega^{\nu_{n_{k-1}}}
\bigl( T_{\nu_{n_{k+1}}} - \bar{T}^{(d_k)}_{\nu_{n_{k+1}}}  \bigr) = 0,
\qquad P\mbox{-a.s.}
\end{equation}
To prove \eqref{ETETbar}, we argue as follows.
From Lemma \ref{ETdiff} we have that for any $\varepsilon>0$
\begin{eqnarray*}
Q \bigl(  E_\omega^{\nu_{n_{k-1}}}  \bigl(  T_{\nu_{n_{k+1}}} -
\bar{T}^{(d_k)}_{\nu_{n_{k+1}}}   \bigr) > \varepsilon  \bigr)
&\leq & n_{k+1}
Q \biggl(  E_\omega T_{\nu} - E_\omega\bar{T}^{(d_k)}_{\nu} >
\frac{\varepsilon}{n_{k+1}}   \biggr)
\\
&=& n_{k+1} \mathcal{O} (n_{k+1}^s e^{-\delta 'b_{d_k}}   ).
\end{eqnarray*}
Since $n_k \sim d_k$, the last term on the
right is summable. Therefore, by the Borel--Cantelli lemma,
%
\begin{equation}\label{QETETbar}
\lim_{k\rightarrow\infty} E_\omega^{\nu_{n_{k-1}}}  \bigl( T_{\nu
_{n_{k+1}}} -
\bar{T}^{(d_k)}_{\nu_{n_{k+1}}}  \bigr) = 0, \qquad Q\mbox{-a.s.}
\end{equation}
This is almost the same as \eqref{ETETbar}, but with $Q$ instead
of $P$. To use this to prove \eqref{ETETbar}, note that for $i>b_n$
using \eqref{reflectexpand}, we can write
\[
E_\omega^{\nu_{i-1}} T_{\nu_i} - E_\omega^{\nu_{i-1}} \bar
{T}^{(n)}_{\nu_i} = A_{i,n}(\omega) + B_{i,n}(\omega) W_{-1},
\]
where $A_{i,n}(\omega)$ and $B_{i,n}(\omega)$ are nonnegative random
variables depending only
on the environment to the right of 0. Thus, $E_\omega^{\nu_{n_{k-1}}}
 ( T_{\nu_{n_{k+1}}} - \bar{T}^{(d_k)}_{\nu_{n_{k+1}}}  )
= A_{d_k}(\omega) + B_{d_k}(\omega) W_{-1}$ where $A_{d_k}(\omega)$ and
$B_{d_k}(\omega)$ are nonnegative and only depend on
the environment to the right of zero (so $A_{d_k}$ and $B_{d_k}$ have
the same distribution under~$P$ as under $Q$). Therefore, \eqref{ETETbar}
follows from \eqref{QETETbar}, which completes the proof of the theorem.
\end{pf}
%
%
\begin{pf*}{Proof of Lemma \ref{lem-latenight}}
Clearly, it suffices to show
the following claims:
%
\begin{equation}\label{irrelevant}
\frac{\bar{T}^{(\tilde{d}_m)}_{x_m}-\bar{T}^{(\tilde{d}_m)}_{\nu_{\beta _m}}
- E_\omega^{\nu_{\beta _m}}\bar{T}^{(\tilde{d}_m)}_{x_m} }{
\sqrt{v_{m,\omega}} } \stackrel{\mathcal{D}_\omega}{\longrightarrow}0
\end{equation}
and
%
\begin{equation}\label{relevant}
\frac{\bar{T}^{(\tilde{d}_m)}_{\nu_{\beta _m}}
- \bar{T}^{(\tilde{d}_m)}_{\nu_{\alpha _m}}
- E_\omega^{\nu_{\alpha _m}}\bar{T}^{(\tilde{d}_m)}_{\nu_{\beta_m}} }
{\sqrt{v_{m,\omega}} } \stackrel{\mathcal{D}_\omega}{\longrightarrow}
Z \sim N(0,1) .
\end{equation}
To prove \eqref{irrelevant}, we note that
\begin{eqnarray*}
P_\omega \biggl( \biggl | \frac{\bar{T}^{(\tilde{d}_m)}_{x_m} -
\bar{T}^{(\tilde{d}_m)}_{\nu_{\beta _m}} - E_\omega^{\nu_{\beta _m}}
\bar{T}^{(\tilde{d}_m)}_{x_m}}{\sqrt{v_{m,\omega}}}  \biggr|
\geq \varepsilon  \biggr) & \leq\frac{ \operatorname{Var}_\omega
(\bar{T}^{(\tilde{d}_m)}_{x_m}-\bar{T}^{(\tilde{d}_m)}_{\nu_{\beta _m}})
}{\varepsilon^2 v_{m,\omega} }
\leq\frac{\sum_{i=\beta_m+1}^{\gamma_m}
\sigma_{i,\tilde{d}_m,\omega} ^2 }{\varepsilon^2 \tilde{a}_m
\tilde{d}_m^{2/s}},
\end{eqnarray*}
where the last inequality is because $x_m \leq\nu_{\gamma_m}$ and
$v_{m,\omega} \geq\tilde{a}_m \tilde{d}_m^{2/s}$. However, by
Corollary \ref{Vsdiff} and the Borel--Cantelli lemma,
\begin{eqnarray*}
\sum_{i=\beta _m+1}^{\gamma_m} \sigma_{i,\tilde{d}_m,\omega} ^2 =
\sum_{i=\beta _m+1}^{\gamma_m} \mu_{i,\tilde{d}_m,\omega} ^2 +
o ((c_{k_m} \tilde{d}_m)^{2/s} ).
\end{eqnarray*}
The application of Corollary \ref{Vsdiff} uses the fact that for
$k$ large enough the reflections ensure that the events in
question do not involve the environment to the left of zero, and
thus have the same probability under $P$ or $Q$. (This type of
argument will be used a few more times in the remainder of the
proof without mention.) By our choice of the subsequence
$n_{k_m}$, we have
\[
\sum_{i=\beta _m+1}^{\gamma_m} \mu_{i,\tilde{d}_m,\omega} ^2 \leq
\Biggl (\sum_{i=\beta _m+1}^{\gamma_m} \mu_{i,\tilde{d}_m,\omega}
\Biggr)^2 \leq4 \tilde{d}_m^{2/s}.
\]
Therefore,
\begin{eqnarray*}
&& \lim_{m\rightarrow\infty} P_\omega \biggl(  \biggl|
\frac{\bar{T}^{(\tilde{d}_m)}_{x_m} -
\bar{T}^{(\tilde{d}_m)}_{\nu_{\beta _m}} - E_\omega^{\nu_{\beta _m}}
\bar{T}^{(\tilde{d}_m)}_{x_m}}{\sqrt{v_{m,\omega}}}  \biggr|\geq \varepsilon
 \biggr)
\\
&&\qquad \leq\lim_{m\rightarrow\infty} \frac{ 4 \tilde{d}_m^{2/s} +
o ((c_{k_m} \tilde{d}_m)^{2/s} ) }{\varepsilon^2 \tilde{a}_m
\tilde{d}_m^{2/s}} = 0, \qquad P\mbox{-a.s.}
\end{eqnarray*}
where the last limit equals zero because $c_k = o(\log a_k)$.

It only remains to prove
\eqref{relevant}. Rewriting, we express
\[
\bar{T}^{(\tilde{d}_m)}_{\nu_{\beta _m}}-\bar{T}^{(\tilde
{d}_m)}_{\nu_{\alpha _m}}
- E_\omega^{\nu_{\alpha _m}}\bar{T}^{(\tilde{d}_m)}_{\nu_{\beta_m}}
= \sum_{i=\alpha _m+1}^{\beta _m}  \bigl( \bigl(\bar{T}^{(\tilde{d}_m)}_{\nu_i}
- \bar{T}^{(\tilde{d}_m)}_{\nu_{i-1}}\bigr)
- \mu_{i,\tilde{d}_m,\omega} \bigr)
\]
as the sum of independent, zero-mean random variables (quenched),
and thus we need only show the Lindberg--Feller condition. That is, we need
to show
%
\begin{equation}\label{LF1}
\lim_{m\rightarrow\infty} \frac{1}{v_{m,\omega}} \sum_{i=\alpha
_m+1}^{\beta _m} \sigma_{i,\tilde{d}_m,\omega}^2 = 1, \qquad
P\mbox{-a.s.},
\end{equation}
and for all $\varepsilon> 0$
%
\begin{eqnarray}\label{LF2}
&& \lim_{m\rightarrow\infty} \frac{1}{v_{m,\omega}}
\sum_{i=\alpha_m+1}^{\beta _m}E_\omega^{\nu_{i-1}}
\bigl[  \bigl(\bar{T}^{(\tilde{d}_m)}_{\nu_i}-\mu_{i,\tilde{d}_m,\omega}  \bigr)^2
\mathbf{1}_{ |\bar{T}^{(\tilde{d}_m)}_{\nu_i}-\mu_{i,\tilde{d}_m,\omega} | >
\varepsilon\sqrt{v_{m,\omega}})}\bigr  ]
\nonumber\\[-8pt]
\\[-8pt]
&&\qquad  = 0, \qquad P\mbox{-a.s.}
\nonumber
\end{eqnarray}
To prove \eqref{LF1}, note that
\begin{eqnarray*}
\frac{1}{v_{m,\omega}} \sum_{i=\alpha _m+1}^{\beta _m} \sigma
_{i,\tilde{d}_m,\omega}^2
= 1 + \frac{\sum_{i=\alpha _m+1}^{\beta _m}  ( \sigma
_{i,\tilde{d}_m,\omega}^2
- \mu_{i,\tilde{d}_m,\omega}^2  ) }{v_{m,\omega}}.
\end{eqnarray*}
However, again by Corollary \ref{Vsdiff} and the Borel--Cantelli lemma,
we have\break $\sum_{i=\alpha _m+1}^{\beta _m} (\sigma_{i,\tilde
{d}_m,\omega}^2 -
\mu_{i,\tilde{d}_m,\omega}^2) = o ((\delta _{k_m} \tilde
{d}_m)^{2/s} )$. Recalling that
$v_{m,\omega} \geq\tilde{a}_m \tilde{d}_m^{2/s}$ we have that
\eqref{LF1} is proved.

To prove \eqref{LF2}, we break the sum up into two parts depending
on whether $M_i$ is ``small'' or ``large.'' Specifically, for
$\varepsilon'\in(0,\frac{1}{3})$,
we decompose the sum as
%
\begin{eqnarray}
\qquad
&& \frac{1}{v_{m,\omega}} \sum_{i=\alpha _m+1}^{\beta _m} E_\omega^{\nu_{i-1}}
\bigl [  \bigl( \bar{T}^{(\tilde{d}_m)}_{\nu_i}-\mu_{i,\tilde{d}_m,\omega}
\bigr)^2 \mathbf{1}_{ | \bar{T}^{(\tilde{d}_m)}_{\nu_i}
-\mu_{i,\tilde{d}_m,\omega} | > \varepsilon\sqrt{v_{m,\omega}})}\bigr]
\mathbf{1}_{M_i \leq\tilde{d}_m^{(1-\varepsilon')/s}}
\label{Msmall} \\
&&\qquad {} + \frac{1}{v_{m,\omega}} \sum_{i=\alpha _m+1}^{\beta _m}
E_\omega^{\nu_{i-1}} \bigl [  \bigl(
\bar{T}^{(\tilde{d}_m)}_{\nu_i}-\mu_{i,\tilde{d}_m,\omega}  \bigr)^2
\mathbf{1}_{ |
\bar{T}^{(\tilde{d}_m)}_{\nu_i}-\mu_{i,\tilde{d}_m,\omega} |
> \varepsilon \sqrt{v_{m,\omega}}}  \bigr]
\nonumber\\[-8pt]\label{Mlarge}
\\[-8pt]
&&\hspace*{92pt}{}\times \mathbf{1}_{M_i
> \tilde{d}_m^{(1-\varepsilon')/s}} .
\nonumber
\end{eqnarray}
We get an upper bound for \eqref{Msmall} by first omitting the
indicator function inside the expectation,
and then expanding the sum to be up
to $n_{k_m} \geq\beta _m$. Thus, \eqref{Msmall} is bounded above by
\[
\frac{1}{v_{m,\omega}}
\sum_{i=\alpha _m+1}^{\beta _m} \sigma_{i,\tilde{d}_m,\omega}^2
\mathbf{1}_{M_i
\leq\tilde{d}_m^{(1-\varepsilon')/s}} \leq\frac{1}{v_{m,\omega}}
\sum_{i=n_{k_m-1}+1}^{n_{k_m}} \sigma_{i,\tilde{d}_m,\omega}^2
\mathbf{1}_{M_i \leq\tilde{d}_m^{(1-\varepsilon')/s}} .
\]
However, since $d_k$ grows exponentially fast, the Borel--Cantelli
lemma and Lemma~\ref{Vsmall} give that
%
\begin{equation}\label{smV}
\sum_{i=n_{k-1}+1}^{n_k} \sigma_{i,d_k,\omega}^2 \mathbf{1}_{M_i
\leq d_k^{(1-\varepsilon')/s}} = o(d_k^{2/s}).
\end{equation}
Therefore, since our choice of the subsequence $n_{k_m}$ gives
that $v_{m,\omega} \geq\tilde{d}_m^{2/s}$, we have
that \eqref{Msmall} tends to zero as $m\rightarrow\infty$.

To get an upper bound for \eqref{Mlarge}, first note that our
choice of the subsequence $n_{k_m}$ gives that $\varepsilon\sqrt
{v_{m,\omega}}
\geq\varepsilon\sqrt{\tilde{a}_m} \mu_{i,\tilde{d}_m,\omega}$ for
any $i\in
(\alpha _m, \beta _m]$. Thus, for $m$ large enough, we can replace the
indicators inside the expectations in \eqref{Mlarge} by the
indicators of the events $\{ \bar{T}_{\nu_i}^{(\tilde{d}_m)} >
(1+\varepsilon\sqrt{\tilde{a}_m}) \mu_{i,\tilde{d}_m,\omega} \}$.
Thus, for
$m$ large enough and $i\in(\alpha _m, \beta _m]$, we have
%
\begin{eqnarray}\label{truncmoment}
&& E_\omega^{\nu_{i-1}}  \bigl[  \bigl(
\bar{T}^{(\tilde{d}_m)}_{\nu_i}-\mu_{i,\tilde{d}_m,\omega}  \bigr)^2
\mathbf{1}_{ | \bar{T}^{(\tilde{d}_m)}_{\nu_i}-\mu_{i,\tilde{d}_m,\omega} | >
\varepsilon\sqrt{v_{m,\omega}}}  \bigr]
\nonumber\\
&&\qquad\leq E_\omega^{\nu_{i-1}}  \bigl[  \bigl(
\bar{T}^{(\tilde{d}_m)}_{\nu_i}-\mu_{i,\tilde{d}_m,\omega}  \bigr)^2
\mathbf{1}_{ \bar{T}^{(\tilde{d}_m)}_{\nu_i}> (1+\varepsilon
\sqrt{\tilde{a}_m})
\mu_{i,\tilde{d}_m,\omega} }  \bigr]
\nonumber\\[-8pt]
\\[-8pt]
&&\qquad= \varepsilon^2\tilde{a}_m \mu_{i,\tilde{d}_m,\omega}^2
P_\omega^{\nu_{i-1}}
\bigl ( \bar{T}_{\nu_i}^{(\tilde{d}_m)} > \bigl(1+\varepsilon\sqrt
{\tilde{a}_m}\bigr) \mu_{i,\tilde{d}_m,\omega} \bigr)
\nonumber\\
&&\qquad\quad{} + \int_{1+\varepsilon\sqrt{\tilde{a}_m}}^{\infty}
P_\omega^{\nu_{i-1}}
\bigl ( \bar{T}_{\nu_i}^{(\tilde{d}_m)} > x \mu_{i,\tilde
{d}_m,\omega} \bigr) 2(x-1)\mu_{i,\tilde{d}_m,\omega}^2 \,dx .
\nonumber
\end{eqnarray}
We want to use Lemma \ref{momentbound} get an upper bounds on the
probabilities in the last line above. Lemma \ref{momentbound} and the
Borel--Cantelli lemma give
that for $k$ large enough, $E_\omega^{\nu_{i-1}} (
\bar{T}^{(d_k)}_{\nu_i}  )^j \leq2^j j! \mu_{i,d_k,\omega}^j$,
for all $n_{k-1}<i\leq n_k$ such that $M_i > d_k^{(1-\varepsilon')/s}$.
Multiplying by $(4\mu_{i,d_k,\omega})^{-j}$ and summing over $j$ gives
that $E_\omega^{\nu_{i-1}} e^{ \bar{T}^{(d_k)}_{\nu_i} /(4
\mu_{i,d_k,\omega}) } \leq2$. Therefore, Chebyshev's inequality
gives that
\begin{eqnarray*}
P_\omega^{\nu_{i-1}}
\bigl( \bar{T}_{\nu_i}^{(d_k)} > x \mu_{i,d_k,\omega} \bigr)
&\leq&  e^{- x/4} E_\omega^{\nu_{i-1}} e^{
\bar{T}^{(d_k)}_{\nu_i} / (4 \mu_{i,d_k,\omega})}
\\
&\leq& 2 e^{-x/4}.
\end{eqnarray*}
Thus, for all $m$ large enough and for all $i$ with\vadjust{\goodbreak}
$\alpha _m<i\leq \beta _m
\leq n_{k_m}$ and $M_i > \tilde{d}_m^{(1-\varepsilon')/s}$, we have
from \eqref{truncmoment} that
\begin{eqnarray*}
&& E_\omega^{\nu_{i-1}}  \bigl[  \bigl( \bar{T}^{(\tilde{d}_m)}_{\nu
_i}-\mu_{i,\tilde{d}_m,\omega}  \bigr)^2
\mathbf{1}_{ |\bar{T}^{(\tilde{d}_m)}_{\nu_i}-\mu_{i,\tilde
{d}_m,\omega} | > \varepsilon\sqrt{v_{m,\omega}}}  \bigr]
\\
&&\qquad\leq\varepsilon^2\tilde{a}_m \mu_{i,\tilde{d}_m,\omega}^2 2
e^{-(1+\varepsilon\sqrt{\tilde{a}_m})/4}
+ \int_{1+\varepsilon\sqrt{\tilde{a}_m}}^{\infty} 2e^{-x/4}
2(x-1)\mu_{i,\tilde{d}_m,\omega}^2 \,dx
\\
&&\qquad=  \bigl( 2\varepsilon^2\tilde{a}_m + 16(4+\varepsilon\sqrt
{\tilde{a}_m})  \bigr) e^{-(1+\varepsilon\sqrt{\tilde{a}_m})/4}
\mu_{i,\tilde{d}_m,\omega}^2.
\end{eqnarray*}
%
%
Recalling the definition of $v_{m,\omega}=\sum_{i=\alpha_m+1}^{\beta _m}
\mu_{i,\tilde{d}_m,\omega}^2$, we have that as $m\rightarrow\infty$,
\eqref{Mlarge}~is bounded above by
\begin{eqnarray*}
&& \lim_{m\rightarrow\infty} \frac{1}{v_{m,\omega}}
\sum_{i=\alpha_m+1}^{\beta _m}
\bigl ( 2\varepsilon^2\tilde{a}_m + 16\bigl(4+\varepsilon\sqrt{\tilde
{a}_m}\bigr)  \bigr) e^{-(1+\varepsilon\sqrt{\tilde{a}_m})/4}
\mu_{i,\tilde{d}_m,\omega}^2 \mathbf{1}_{M_i > \tilde
{d}_m^{(1-\varepsilon')/s}}
\\
&&\qquad\leq\lim_{m\rightarrow\infty}  \bigl( 2\varepsilon^2\tilde{a}_m
+ 16\bigl(4+\varepsilon\sqrt{\tilde{a}_m}\bigr)  \bigr)
e^{-(1+\varepsilon\sqrt{\tilde{a}_m})/4} = 0.
\end{eqnarray*}
This completes the proof of \eqref{LF2}, and thus of Lemma \ref{lem-latenight}.
\end{pf*}
\begin{pf*}{Proof of Theorem \ref{nonlocal}}
Note first that from Lemma \ref{mdevu} and the
Borel--Cantelli lemma, we have that for any $\varepsilon>0$, $E_\omega
T_{\nu_{n_k}} = o(n_k^{(1+\varepsilon)/s})$, $P$-a.s. This is
equivalent to
%
\begin{equation}\label{eq-100407d}
\limsup_{k\rightarrow\infty} \frac{\log E_\omega T_{\nu
_{n_k}}}{\log n_k} \leq \frac{1}{s}, \qquad P\mbox{-a.s.}
\end{equation}
We can also get bounds on the
probability of $E_\omega T_{\nu_n}$ being small. Since $E_\omega^{\nu_{i-1}}
T_{\nu_i} \geq M_i$, we have
\[
P \bigl(E_\omega T_{\nu_n} \leq n^{(1-\varepsilon)/s} \bigr)
\leq P \bigl( M_i \leq n^{(1-\varepsilon)/s}, \ \forall i\leq
n  \bigr)
\leq \bigl( 1 - P \bigl(M_1 > n^{(1-\varepsilon)/s}  \bigr) \bigr)^n,
\]
and since $P(M_1 > n^{(1-\varepsilon)/s}) \sim C_5 n^{-1+\varepsilon}$; see
\eqref{Mtail}, we have
$P (E_\omega T_{\nu_n} \leq\break n^{(1-\varepsilon)/s} ) \leq
e^{-n^{\varepsilon/2}}$. Thus, by the Borel--Cantelli lemma, for any
$\varepsilon>0$,
we have that $E_\omega T_{\nu_{n_k}} \geq n_k^{(1-\varepsilon)/s}$ for
all $k$ large enough, $P$-a.s., or equivalently
%
\begin{equation}\label{eq-100407c}
\liminf_{k\rightarrow\infty} \frac{\log E_\omega T_{\nu
_{n_k}}}{\log n_k} \geq \frac{1}{s}, \qquad P\mbox{-a.s.}
\end{equation}
Let $n_{k_m}$ be the subsequence specified in Theorem
\ref{gaussianT} and define $t_m:= E_\omega T_{n_{k_m}}$. Then
by \eqref{eq-100407d} and \eqref{eq-100407c},
$\lim_{m\rightarrow\infty} \frac{\log t_m}{\log n_{k_m}} = {1}/{s}$.

For any $t$, define $X_t^*:= \max\{ X_n \dvtx n\leq t \}$. Then for
any $x\in(0,\infty)$, we have
\begin{eqnarray*}
P_\omega \biggl( \frac{X_{t_m}^*}{n_{k_m}} < x  \biggr)
&=& P (X_{t_m}^* < x n_{k_m}  )= P_\omega ( T_{x n_{k_m}} > t_m  )
\\
&=& P_\omega \biggl( \frac{T_{x n_{k_m}} - E_\omega T_{x
n_{k_m}}}{\sqrt{v_{m,\omega}}}
> \frac{ E_\omega T_{n_{k_m}}-E_\omega T_{x n_{k_m}}}{\sqrt
{v_{m,\omega}}} \biggr).
\end{eqnarray*}
Now, with notation as in Theorem \ref{gaussianT}, we have that for
all $m$ large enough\break $\nu_{\beta _m} < x n_{k_m} < \nu_{\gamma_m}$
(note that this also uses the fact that $\nu_n/n \rightarrow E_P
\nu$,\break
\mbox{$P$-a.s.}). Thus, $\frac{T_{x n_{k_m}} - E_\omega T_{x
n_{k_m}}}{\sqrt{v_{m,\omega}}} \stackrel{\mathcal{D}_\omega
}{\longrightarrow}Z\sim N(0,1)$. Then we will have
proved that $ \lim_{m\rightarrow\infty} P_\omega ( \frac
{X_{t_m}^*}{n_{k_m}} <
x  ) = \frac{1}{2}$ for any $x\in(0,\infty)$, if we can show
%
\begin{equation}\label{sd}
\lim_{m\rightarrow\infty} \frac{ E_\omega T_{n_{k_m}}-E_\omega T_{x
n_{k_m}}}{\sqrt{v_{m,\omega}}} = 0,\qquad P\mbox{-a.s.}
\end{equation}
For $m$ large enough, we have $n_{k_m}, x n_{k_m} \in(\nu_{\beta _m} ,
\nu_{\gamma_m})$.
Thus, for $m$ large enough,
\begin{eqnarray*}
\biggl |\frac{ E_\omega T_{x n_{k_m}}-E_\omega T_{n_{k_m}}}{\sqrt{v_{m,\omega
}}} \biggr|
&\leq &\frac{E_\omega^{\nu_{\beta _m}} T_{\nu_{\gamma_m}} }{\sqrt
{v_{m,\omega}}}
\\
&=& \frac{1}{\sqrt{v_{m,\omega}}}
\Biggl( E_\omega^{\nu_{\beta _m}} \bigl(
T_{\nu_{\gamma_m}} - \bar{T}_{\nu_{\gamma_m}}^{(\tilde{d}_m)}
 \bigr)
 + \sum_{i=\beta _m+1}^{\gamma_m} \mu_{i,\tilde{d}_m,\omega} \Biggr).
\end{eqnarray*}
Since $\alpha _m \leq\beta _m \leq\gamma_m \leq n_{k_m+1}$ for all
$m$ large enough, we can apply \eqref{ETETbar} to get
\[
\lim_{m\rightarrow\infty} E_\omega^{\nu_{\beta _m}}
\bigl(T_{\nu_{\gamma_m}} -
\bar{T}_{\nu_{\gamma_m}}^{(\tilde{d}_m)}  \bigr) \leq
\lim_{m\rightarrow\infty} E_\omega^{\nu_{\alpha _m}}
\bigl ( T_{\nu_{n_{k_m+1}}} -
\bar{T}_{\nu_{n_{k_m+1}}}^{(\tilde{d}_m)}  \bigr)=0.
\]
Also, from our choice of $n_{k_m}$ we have that
$\sum_{i=\beta _m+1}^{\gamma_m} \mu_{i,\tilde{d}_m,\omega} \leq2
\tilde{d}_m^{1/s}$ and $v_{m,\omega} \geq\tilde{a}_m
\tilde{d}_m^{2/s}$. Thus \eqref{sd} is proved. Therefore
\[
\lim_{m\rightarrow\infty} P_\omega \biggl(\frac{X_{t_m}^*}{n_{k_m}}
\leq x  \biggr) = \frac{1}{2} \qquad\forall x\in(0,\infty),
\]
and obviously $\lim_{m\rightarrow\infty}
P_\omega (\frac{X_{t_m}^*}{n_{k_m}} < 0  ) = 0$ since
$X_n$ is transient to the right $\mathbb{P}$-a.s. due to Assumption
\ref{essentialasm}. Finally, note that
\[
\frac{X_t^* - X_t}{\log^2 t} = \frac{X_t^* - \nu_{N_t}}{\log^2 t}
+ \frac{\nu_{N_t} - X_t}{\log^2 t} \leq\frac{\max_{i\leq t}
(\nu_i - \nu_{i-1}) }{\log^2 t} + \frac{\nu_{N_t} - X_t}{\log^2 t}.
\]
However, Lemma \ref{seperation} and an easy application of
Lemma \ref{nutail} and the Borel--Cantelli lemma gives that
\[
\lim_{t\rightarrow\infty} \frac{X_t^* - X_t}{\log^2 t} = 0, \qquad
\mathbb{P}\mbox{-a.s.}
\]
This completes the proof of the theorem.
\end{pf*}


\section{Asymptotics of the tail of $E_\omega T_\nu$} \label{tailofTnu}

Recall that
$
E_\omega T_\nu= \nu+ 2\sum_{j=0}^{\nu-1} W_j
= \nu+ 2 \sum_{i\leq j, 0\leq j < \nu} \Pi_{i,j} ,
$
and for any $A>1$ define
\[
\sigma= \sigma_A = \inf\{ n \geq1\dvtx \Pi_{0,n-1} \geq A \}.
\]
Note that $\sigma-1$ is a stopping time for the sequence $\Pi_{0,k}$. For
any $A>1$, $\{\sigma> \nu\} = \{ M_1 < A \}$. Thus, we have by
\eqref{TbigMsmall} that for any $A>1$,
%
\begin{equation}\label{long}
Q(E_\omega T_\nu> x, \sigma> \nu) = Q(E_\omega T_\nu> x, M_1 < A) =
o(x^{-s}).
\end{equation}
Thus, we may focus on the tail estimates $Q(E_\omega T_\nu> x, \sigma
< \nu)$
in which case we can use the following expansion of $E_\omega T_\nu$:
%
\begin{eqnarray}\label{expand}
E_\omega T_\nu&=& \nu+ 2 \sum_{i<0\leq j<\sigma-1} \Pi_{i,j}
+ 2 \sum_{0 \leq i \leq j<\sigma-1} \Pi_{i,j}
\nonumber\\
&&{} + 2 \sum_{\sigma\leq i\leq j < \nu} \Pi_{i,j}
+ 2 \sum_{i\leq\sigma-1 \leq j < \nu} \Pi_{i,j}
\nonumber\\[-8pt]
\\[-8pt]
&=& \nu+ 2 W_{-1}R_{0,\sigma-2} + 2 \sum_{j=0}^{\sigma-2} W_{0,j}
\nonumber\\
&&{} + 2 \sum_{i=\sigma}^{\nu-1} R_{i,\nu-1}
+ 2 W_{\sigma-1}(1+R_{\sigma,\nu-1}).
\nonumber
\end{eqnarray}
We will show that the dominant term in \eqref{expand} is the last
term: $2 W_{\sigma-1}(1+R_{\sigma,\nu-1})$. A few easy consequences of
Lemmas \ref{nutail} and \ref{Wtail} are that the tails of the first
three terms in the expansion \eqref{expand} are negligible. The
following statements are true for any $\delta >0$ and any $A>1$:
%
\begin{eqnarray}
Q(\nu> \delta x) &=& P(\nu>\delta x) = o(x^{-s}),
\label{first}\\
\qquad
Q(2W_{-1}R_{0,\sigma-2}>\delta x, \sigma< \nu)
&\leq & Q\bigl(W_{-1}> \sqrt{\delta x}\bigr)
\nonumber\\
&&{} + P\bigl(2R_{0,\sigma-2} > \sqrt{\delta x},\sigma<\nu\bigr)
\nonumber\\[-8pt]\label{second}
\\[-8pt]
&\leq & Q\bigl(W_{-1}> \sqrt{\delta x}\bigr)
+ P\bigl(2 \nu A > \sqrt{\delta x}\bigr)
\nonumber\\
&=& o(x^{-s}),
\nonumber\\
Q \Biggl(2 \sum_{j=0}^{\sigma-2} W_{0,j} > \delta x,
\sigma<\nu \Biggr) &\leq&
P \Biggl(2 \sum_{j=1}^{\sigma-1} j A > \delta x,
\sigma< \nu \Biggr)
\nonumber\\[-8pt]\label{third}
\\[-8pt]
&\leq & P( \nu^2 A > \delta x) = o(x^{-s}).
\nonumber
\end{eqnarray}
In the first inequality in \eqref{third}, we used the fact that $\Pi
_{i,j} \leq\Pi_{0,j}$ for any $0<i<\nu$ since $\Pi_{0,i-1} \geq1$.

The fourth term in \eqref{expand} is not negligible, but we can make
it arbitrarily small by taking $A$ large enough.
\begin{lem} \label{fourth}
For all $\delta >0$, there exists an $A_0=A_0(\delta )<\infty$ such that
\[
P \Biggl( 2 \sum_{\sigma_A\leq i < \nu} R_{i,\nu-1} > \delta x
 \Biggr) < \delta x^{-s}\qquad\forall A\geq A_0(\delta ).
\]
\end{lem}
\begin{pf}
This proof is essentially a copy of the proof of Lemma 3 in \cite{kksStable}.
\begin{eqnarray*}
P \Biggl( 2 \sum_{\sigma_A\leq i < \nu} R_{i,\nu-1} > \delta x \Biggr)
&\leq & P \Biggl( \sum_{\sigma_A \leq i < \nu} R_i > \frac{\delta
}{2} x  \Biggr)
\\
&=& P \Biggl( \sum_{i=1}^\infty\mathbf{1}_{\sigma_A \leq i <
\nu } R_i > \frac{\delta }{2} x \frac{6}{\pi^2}
\sum_{i=1}^\infty i^{-2} \Biggr)
\\
&\leq& \sum_{i=1}^\infty P \biggl( \mathbf{1}_{\sigma_A \leq i < \nu} R_i
> x \frac{3\delta }{\pi^2} i^{-2}  \biggr).
\end{eqnarray*}
However, since the event $\{\sigma_A \leq i < \nu\}$ depends only on
$\rho_j$ for $j<i$, and $R_i$ depends only on $\rho_j$ for $j\geq i$,
we have that
\[
P \Biggl( 2 \sum_{\sigma_A\leq i < \nu} R_{i,\nu-1} > \delta x
 \Biggr) \leq \sum_{i=1}^\infty P ( \sigma_A \leq i < \nu )
 P \biggl(R_i > x \frac{3\delta }{\pi^2} i^{-2}  \biggr).
\]
Now, from \eqref{Rtail}, we have that there exists a $K_1 > 0$ such
that $P(R_0 > x) \leq K_1 x^{-s}$ for all $x>0$. We then conclude that
%
\begin{eqnarray}\label{small}
P \Biggl( \sum_{\sigma_A\leq i < \nu} R_{i,\nu-1} > \delta x  \Biggr)
&\leq & K_1  \biggl( \frac{3\delta }{\pi^2} \biggr)^{-s}x^{-s}
\sum_{i=1}^\infty P ( \sigma_A \leq i < \nu ) i^{2s}
\nonumber\\
&=& K_1  \biggl( \frac{3\delta }{\pi^2} \biggr)^{-s}x^{-s} E_P
\Biggl[\sum_{i=1}^{\infty} \mathbf{1}_{\sigma_A\leq i < \nu} i^{2s}
\Biggr]
\\
&\leq& K_1  \biggl( \frac{3\delta }{\pi^2} \biggr)^{-s}
x^{-s} E_P[ \nu^{2s+1} \mathbf{1}_{\sigma_A < \nu}].
\nonumber
\end{eqnarray}
Since $E_P \nu^{2s+1}<\infty$ and $\lim_{A\rightarrow\infty}
P(\sigma_A<\nu) = 0$, we have that the right side
of~\eqref{small} can be made less than $\delta x^{-s}$ by choosing $A$
large enough.
\end{pf}

We need one more lemma before analyzing the dominant term in \eqref{expand}.
\begin{lem}\label{finitemoment}
$E_Q [ W_{\sigma_A-1}^s\mathbf{1}_{\sigma_A<\nu}  ]<\infty$
for any $A>1$.
\end{lem}
\begin{pf}
First, note that on the event $\{\sigma_A < \nu\}$, we have that
$\Pi_{i,\sigma_A-1} \leq\Pi_{0,\sigma_A-1}$ for any
$i\in[0,\sigma_A)$. Thus,
\[
W_{\sigma_A-1} = W_{0,\sigma_A-1} + \Pi_{0,\sigma_A-1}W_{-1} \leq
(\sigma_A + W_{-1})\Pi_{0,\sigma_A-1}.
\]
Also, note that $\Pi_{0,\sigma_A -1} \leq A \rho_{\sigma_A-1}$ by the
definition of $\sigma_A$. Therefore,
\[
E_Q [ W_{\sigma_A-1}^s\mathbf{1}_{\sigma_A<\nu}  ] \leq
E_Q [ (\sigma_A + W_{-1})^s A^s\rho_{\sigma_A-1}^s
\mathbf{1}_{\sigma_A<\nu} ].
\]
Therefore, it is enough to prove that both $E_Q [ W_{-1}^s\rho
_{\sigma_A-1}^s \mathbf{1}_{\sigma_A<\nu}  ]$ and
$E_Q [\sigma_A^s \times\break \rho_{\sigma_A-1}^s\mathbf{1}_{\sigma_A<\nu}  ]$ are finite
(note that this is trivial if we assume that $\rho$ has bounded
support). Since $W_{-1}$ is independent of $\rho_{\sigma_A-1}^s
\mathbf{1}_{\sigma_A<\nu}$ we have that
\[
E_Q [ W_{-1}^s\rho_{\sigma_A-1}^s \mathbf{1}_{\sigma_A<\nu} ]
= E_Q[W_{-1}^s]E_P[\rho_{\sigma_A-1}^s \mathbf{1}_{\sigma_A<\nu}],
\]
where we may take the second expectation over $P$ instead of $Q$
because\break the random variable only depends on the environment to the
right of zero.\break By Lemma \ref{Wtail}, we have that
$E_Q[W_{-1}^s]<\infty$. Also, $E_P[\rho_{\sigma_A-1}^s
\mathbf{1}_{\sigma_A<\nu}] \leq\break E_P[\sigma_A^s\rho_{\sigma_A-1}^s
\mathbf{1}_{\sigma_A<\nu}]$,
and so the lemma will be proved once we prove the latter is finite. However,
\begin{eqnarray*}
E_P [ \sigma_A^s \rho_{\sigma_A-1}^s \mathbf{1}_{\sigma_A<\nu}  ]
= \sum_{k=1}^\infty E_P  [ k^s \rho_{k-1}^s \mathbf{1}_{\sigma_A
= k <\nu}  ] \leq\sum_{k=1}^\infty k^s E_P [\rho_{k-1}^s
\mathbf{1}_{ k \leq\nu}  ],
\end{eqnarray*}
and since the event $\{k\leq\nu\}$ depends only on $(\rho_0,\rho
_1,\ldots\rho_{k-2})$ we have that $E_P [ \rho_{k-1}^s\mathbf
{1}_{k\leq\nu} ] = E_P\rho^s P(\nu\geq k)$ since $P$ is a
product measure. Then since $E_P \rho^s = 1$, we have that
\[
E_P [ \sigma_A^s \rho_{\sigma_A-1}^s \mathbf{1}_{\sigma_A<\nu
}  ]\leq\sum_{k=1}^\infty k^s P(\nu\geq k).
\]
This last sum is finite by Lemma \ref{nutail}.
\end{pf}

Finally, we turn to the asymptotics of the tail of
$2W_{\sigma-1}(1+R_{\sigma,\nu-1})$, which is the dominant
term in \eqref{expand}.
\begin{lem}\label{fifth}
For any $A>1$, there exists a constant $K_A\in(0,\infty)$ such that
\[
\lim_{x\rightarrow\infty} x^s Q \bigl( W_{\sigma-1}(1+R_{\sigma
,\nu-1})>x , \sigma<\nu \bigr) = K_A.
\]
%
\end{lem}
\begin{pf}
The strategy of the proof is as follows. First, note that on the event
$\{\sigma<\nu\}$ we have $W_{\sigma-1}(1 + R_\sigma) = W_{\sigma
-1}(1 + R_{\sigma,\nu
-1}) + W_{\sigma-1}\Pi_{\sigma,\nu-1}R_\nu$. We will begin by analyzing
the asymptotics of the tails of $W_{\sigma-1}(1 + R_\sigma)$ and
$W_{\sigma-1}\Pi_{\sigma,\nu-1}R_\nu$. Next, we will show that
$W_{\sigma-1}(1 + R_{\sigma,\nu-1})$ and $W_{\sigma-1}\Pi_{\sigma,\nu-1}R_\nu$
are essentially
independent in the sense that they cannot both be large. This will
allow us to use the asymptotics of the tails of $W_{\sigma-1}(1 +R_\sigma)$
and $W_{\sigma-1}\Pi_{\sigma,\nu-1}R_\nu$ to compute the
asymptotics of the tails of $W_{\sigma-1}(1 + R_{\sigma,\nu-1})$.

To analyze the asymptotics of the tail of $W_{\sigma-1}(1 + R_\sigma)$, we
first recall from~\eqref{Rtail} that there exists a $K>0$ such that
$P(R_0 > x)\sim K x^{-s}$. Let $\mathcal{F}_{\sigma-1} = \sigma(\ldots,\omega
_{\sigma-2}, \omega_{\sigma-1})$ be the $\sigma$-algebra generated
by the environment to the left of~$\sigma$. Then on the event $\{\sigma
<\infty\}$,
$R_\sigma$ has the same distribution as $R_0$ and is independent of
$\mathcal{F}_{\sigma-1}$. Thus,
%
\begin{eqnarray}\label{Rstail}
&& \lim_{x\rightarrow\infty} x^s Q\bigl(W_{\sigma-1}(1+R_\sigma) > x,
\sigma<\nu\bigr)
\nonumber\\
&&\qquad = \lim_{x\rightarrow\infty} E_Q \biggl[ x^s Q \biggl( 1+R_\sigma
>\frac{x}{W_{\sigma -1}} , \sigma<\nu\Big | \mathcal{F}_{\sigma-1}  \biggr)  \biggr]
\\
&&\qquad = K E_Q  [ W_{\sigma-1}^s \mathbf{1}_{\sigma<\nu}  ].
\nonumber
\end{eqnarray}
A similar calculation yields
%
\begin{eqnarray}\label{prodtail}
&& \lim_{x\rightarrow\infty} x^s Q ( W_{\sigma-1}\Pi_{\sigma
,\nu-1} R_\nu> x, \sigma<\nu )
\nonumber\\
&&\qquad = \lim_{x\rightarrow\infty} E_Q  \biggl[ x^s Q \biggl( R_\nu>
\frac{x}{W_{\sigma-1}\Pi_{\sigma,\nu-1}}, \sigma<\nu\Big |
\mathcal{F}_{\nu-1} \biggr)  \biggr]
\\
&&\qquad = E_Q  [ W_{\sigma-1}^s \Pi_{\sigma,\nu-1}^s \mathbf{1}_{\sigma<\nu} ] K.
\nonumber
\end{eqnarray}
Next, we wish to show that
%
\begin{equation}\label{both}
\qquad
\lim_{x\rightarrow\infty} x^s Q \bigl(W_{\sigma-1}(1+R_{\sigma,\nu
-1}) >\varepsilon x, W_{\sigma-1}\Pi_{\sigma,\nu-1}R_\nu> \varepsilon
x,\, \sigma <\nu \bigr) = 0.
\hspace*{-5pt}
\end{equation}
Since $\Pi_{\sigma,\nu-1} < \frac{1}{A}$ on the event
$\{\sigma<\nu\}$, we have for any $\varepsilon>0$ that
%
\begin{eqnarray}\label{bothbig}
\qquad
&& x^s Q \bigl(W_{\sigma-1}(1+R_{\sigma,\nu-1}) > \varepsilon x,
W_{\sigma -1}\Pi_{\sigma,\nu-1}R_\nu> \varepsilon x, \sigma<\nu \bigr)
\nonumber\\
&&\qquad\leq x^s Q \bigl(W_{\sigma-1}(1+R_{\sigma,\nu-1}) >
\varepsilon x, W_{\sigma-1}R_\nu> A \varepsilon x, \sigma<\nu \bigr)
\nonumber\\
&&\qquad= x^s E_Q  \biggl[ Q \biggl(1+R_{\sigma,\nu-1} > \frac
{\varepsilon x}{W_{\sigma-1}} \Big| \mathcal{F}_{\sigma-1} \biggr)
\nonumber\\[-8pt]
\\[-8pt]
&&\hspace*{64pt}{}\times
Q \biggl( R_\nu> A\frac{\varepsilon x}{W_{\sigma-1}}
\Big |\mathcal{F}_{\sigma-1}  \biggr)
\mathbf{1}_{\sigma<\nu}  \biggr]
\nonumber\\
&& \qquad\leq E_Q  \biggl[ x^s Q\biggl (1+R_{\sigma} > \frac{\varepsilon
x}{W_{\sigma-1}} \Big| \mathcal{F}_{\sigma-1} \biggr)
\nonumber\\
&&\hspace*{52pt}{}\times
Q \biggl( R_\nu> A\frac{\varepsilon x}{W_{\sigma-1}} \Big|
\mathcal{F}_{\sigma-1}  \biggr)
\mathbf{1}_{\sigma<\nu}  \biggr],
\nonumber
\end{eqnarray}
where the equality on the third line is because $R_{\sigma,\nu-1}$ and
$R_\nu$ are independent when $\sigma<\nu$ (note that
$\{ \sigma<\nu\} \in \mathcal{F}_{\sigma-1}$),
and the last inequality is because
$R_{\sigma,\nu-1} \leq R_\sigma$. Now, conditioned on $\mathcal{F}_{\sigma-1}$,
$R_\sigma$
and $R_\nu$ have the same distribution as $R_0$. Then since by
\eqref{Rtail} for any $\gamma\leq s$, there exists a $K_\gamma>0$ such that
$P(1+R_0 > x)\leq K_\gamma x^{-\gamma}$, we have that the integrand
in \eqref{bothbig} is bounded above by $K_\gamma^2\varepsilon
^{-2\gamma} W_{\sigma-1}^{2\gamma}\mathbf{1}_{\sigma<\nu}x^{s-2\gamma
}$, $Q$-a.s. Choosing $\gamma=\frac{s}{2}$ gives that the integrand
in~\eqref{bothbig} is $Q$-a.s. bounded above by $K^2_{\frac{s}{2}}
\varepsilon^{-s} W_{\sigma-1}^s \mathbf{1}_{\sigma<\nu}$ which by Lemma
\ref{finitemoment} has finite mean. However, if we choose $\gamma=s$,
then we get that the integrand of \eqref{bothbig} tends to zero
$Q$-a.s. as $x\rightarrow\infty$. Thus, by the dominated convergence theorem,
we have that \eqref{both} holds.

Now, since $R_\sigma= R_{\sigma,\nu-1} + \Pi_{\sigma,\nu-1}R_\nu$,
we have that for any $\varepsilon>0$,
\begin{eqnarray*}
&& Q\bigl(W_{\sigma-1}(1+R_\sigma)>(1+\varepsilon)x, \sigma<\nu\bigr)
\\
&&\qquad \leq Q\bigl(W_{\sigma
-1}(1+R_{\sigma,\nu-1})>\varepsilon x,
W_{\sigma-1} \Pi_{\sigma,\nu-1}R_\nu> \varepsilon x, \sigma<\nu\bigr)
\\
&&\qquad \quad  {} + Q\bigl(W_{\sigma-1}(1+R_{\sigma,\nu-1}) > x,\sigma<\nu\bigr)
 + Q\bigl(W_{\sigma-1}\Pi_{\sigma,\nu-1}R_\nu> x,\sigma<\nu\bigr).
\end{eqnarray*}
Applying \eqref{Rstail},
\eqref{prodtail} and \eqref{both}, we get that for any $\varepsilon>0$,
%
\begin{eqnarray}\label{lb}
&& \liminf_{x\rightarrow\infty} x^s Q\bigl(W_{\sigma-1}(1+R_{\sigma,\nu
-1}) > x, \sigma<\nu\bigr)
\nonumber\\[-8pt]
\\[-8pt]
&&\qquad \geq K E_Q[W_{\sigma-1}^s\mathbf{1}_{\sigma<\nu}](1+\varepsilon)^{-s}
- K E_Q[W_{\sigma-1}^s\Pi_{\sigma,\nu-1}^s\mathbf{1}_{\sigma<\nu}].
\nonumber
\end{eqnarray}
Similarly, for a bound in the other direction, we have
\begin{eqnarray*}
&& Q\bigl(W_{\sigma-1}(1+R_\sigma)>x, \sigma<\nu\bigr)
\\
&&\qquad \geq  Q\bigl(W_{\sigma-1}(1+R_{\sigma,\nu-1}) > x, \mbox{ or }
W_{\sigma -1}\Pi_{\sigma,\nu-1}R_\nu> x, \sigma<\nu\bigr)
\\
&&\qquad = Q\bigl(W_{\sigma-1}(1+R_{\sigma,\nu-1})>x, \sigma<\nu\bigr)
+ Q(W_{\sigma-1}\Pi_{\sigma,\nu-1}R_\nu>x, \sigma< \nu)
\\
&&\qquad \quad {} - Q\bigl(W_{\sigma-1}(1+R_{\sigma,\nu-1})>x, W_{\sigma-1}
\Pi_{\sigma,\nu-1}R_\nu>x, \sigma< \nu\bigr).
\end{eqnarray*}
Thus, again applying \eqref{Rstail}, \eqref{prodtail} and
\eqref{both}, we get
%
\begin{eqnarray}\label{ub}
&& \limsup_{x\rightarrow\infty} x^s
Q\bigl(W_{\sigma-1}(1+R_{\sigma,\nu-1}) > x, \sigma<\nu\bigr)
\nonumber\\[-8pt]
\\[-8pt]
&&\qquad \leq K E_Q[W_{\sigma-1}^s\mathbf{1}_{\sigma<\nu}] - K
E_Q[W_{\sigma-1}^s\Pi_{\sigma,\nu-1}^s\mathbf{1}_{\sigma<\nu}].
\nonumber
\end{eqnarray}
Finally, applying \eqref{lb} and \eqref{ub} and letting $\varepsilon
\rightarrow0$, we get that
\begin{eqnarray*}
&& \lim_{x\rightarrow\infty} x^s Q\bigl(W_{\sigma-1}(1+R_{\sigma,\nu-1})
> x, \sigma<\nu\bigr)
\\
&&\qquad =
K E_Q[W_{\sigma-1}^s(1-\Pi_{\sigma,\nu-1}^s)\mathbf{1}_{\sigma
<\nu}] =: K_A,
\end{eqnarray*}
and $K_A \in(0,\infty)$ by Lemma \ref{finitemoment} and the fact
that $1-\Pi_{\sigma,\nu-1} \in(1-\frac{1}{A}, 1)$.
\end{pf}

Finally, we are ready to analyze the tail of $E_\omega T_\nu$ under the
measure $Q$.
\begin{pf*}{Proof of Theorem \ref{Tnutail}}
Let $\delta >0$, and choose $A\geq A_0(\delta )$ as in Lemma~\ref{fourth}.
Then using \eqref{expand}, we have
\begin{eqnarray*}
Q(E_\omega T_\nu> x)
& =& Q(E_\omega T_\nu> x, \sigma>\nu)+ Q(E_\omega T_\nu> x, \sigma<\nu)
\\
& \leq & Q(E_\omega T_\nu> x, \sigma>\nu)+ Q(\nu> \delta x)
\\
&&{} + Q(2W_{-1}R_{0,\sigma -2}> \delta x, \sigma<\nu)
\\
&&{} + Q \Biggl(2\sum_{j=0}^{\sigma-2} W_{0,j} > \delta x,
\sigma<\nu  \Biggr)
+ Q \Biggl(2\sum_{\sigma\leq i < \nu} R_{i,\nu-1} > \delta x \Biggr)
\\
&&{} + Q\bigl(2W_{\sigma-1}(1+R_{\sigma,\nu-1})>(1-4\delta )x,\sigma<\nu\bigr).
\end{eqnarray*}
Thus, combining
equations \eqref{long}, \eqref{first}, \eqref{second}
and \eqref{third} and Lemmas \ref{fourth} and \ref{fifth}, we get that
%
\begin{equation}\label{ub2}
\limsup_{x\rightarrow\infty} x^s Q(E_\omega T_\nu> x)
\leq\delta+ 2^{s}K_A(1-4\delta)^{-s}.
\end{equation}
The lower bound is easier, since $Q(E_\omega T_\nu> x) \geq
Q(2W_{\sigma
-1}(1+R_{\sigma,\nu-1})>x, \sigma<\nu)$. Thus,
%
\begin{equation}\label{lb2}
\liminf_{x\rightarrow\infty} x^s Q(E_\omega T_\nu> x)
\geq 2^{s}K_A.
\end{equation}
From \eqref{ub2} and \eqref{lb2}, we get that $\overline{K}:=\limsup
_{A\rightarrow\infty} 2^s K_A < \infty$. Therefore, letting
$\underline{K}:= \liminf_{A\rightarrow\infty} 2^s K_A$, we have
from \eqref{ub2} and \eqref{lb2} that
\[
\overline{K} \leq\liminf_{x\rightarrow\infty} x^s
Q(E_\omega T_\nu > x) \leq
\limsup_{x\rightarrow\infty} x^s Q(E_\omega T_\nu> x) \leq\delta
+\underline{K}(1-4\delta )^{-s}.
\]
Then letting $\delta \rightarrow0$ completes the proof of the theorem with
$K_\infty= \underline{K}=\overline{K}$.
\end{pf*}

\section*{Acknowledgment}

We thank the referee for a very thorough and careful reading of the
paper and for the useful suggestions.


\printaddresses

\begin{thebibliography}{14}

\bibitem{dzLDTA}
\begin{bbook}[msn]
\bauthor{\bsnm{Dembo},~\bfnm{Amir}\binits{A.}} \AND
  \bauthor{\bsnm{Zeitouni},~\bfnm{Ofer}\binits{O.}}
(\byear{1998}).
\btitle{Large Deviations Techniques and Applications},
\bedition{2nd} ed.
\bseries{Applications of Mathematics (New York)}
\bvolume{38}.
\bpublisher{Springer}, \baddress{New York}.
\bmrnumber{MR1619036}
\end{bbook}
\endbibitem

\bibitem{ESZ}
\begin{btechreport}[vtex]
\bauthor{\bsnm{Enriquez},~\bfnm{N.}\binits{N.}},
  \bauthor{\bsnm{Sabot},~\bfnm{C.}\binits{C.}} \AND
  \bauthor{\bsnm{Zindy},~\bfnm{O.}\binits{O.}}
(\byear{2007}).
\btitle{Limit laws for transient random walks in random environment on $\mathbb Z$}.
\btype{Preprint}.
\href{http://arxiv.org/abs/math/0703660v3}{arXiv:math/0703660v3 [math.PR]}.
\end{btechreport}
\endbibitem

\bibitem{fpKMF}
\begin{barticle}[vtex]
\bauthor{\bsnm{Fitzsimmons},~\bfnm{P.~J.}\binits{P.~J.}} \AND
  \bauthor{\bsnm{Pitman},~\bfnm{Jim}\binits{J.}}
(\byear{1999}).
\btitle{Kac's moment formula and the {F}eynman--{K}ac formula for additive
  functionals of a {M}arkov process}.
\bjournal{Stochastic Process. Appl.}
\bvolume{79}
\bpages{117--134}.
\bmrnumber{MR1670526}
\end{barticle}
\endbibitem

\bibitem{gsMVSS}
\begin{barticle}[msn]
\bauthor{\bsnm{Gantert},~\bfnm{Nina}\binits{N.}} \AND
  \bauthor{\bsnm{Shi},~\bfnm{Zhan}\binits{Z.}}
(\byear{2002}).
\btitle{Many visits to a single site by a transient random walk in random
  environment}.
\bjournal{Stochastic Process. Appl.}
\bvolume{99}
\bpages{159--176}.
\bmrnumber{MR1901151}
\end{barticle}
\endbibitem

\bibitem{gQCLT}
\begin{barticle}[msn]
\bauthor{\bsnm{Goldsheid},~\bfnm{Ilya~Ya.}\binits{I.~Y.}}
(\byear{2007}).
\btitle{Simple transient random walks in one-dimensional random environment:
  The central limit theorem}.
\bjournal{Probab. Theory Related Fields}
\bvolume{139}
\bpages{41--64}.
\bmrnumber{MR2322691}
\end{barticle}
\endbibitem

\bibitem{iEV}
\begin{barticle}[msn]
\bauthor{\bsnm{Iglehart},~\bfnm{Donald~L.}\binits{D.~L.}}
(\byear{1972}).
\btitle{Extreme values in the {$GI/G/1$} queue}.
\bjournal{Ann. Math. Statist.}
\bvolume{43}
\bpages{627--635}.
\bmrnumber{MR0305498}
\end{barticle}
\endbibitem

\bibitem{kRDE}
\begin{barticle}[msn]
\bauthor{\bsnm{Kesten},~\bfnm{Harry}\binits{H.}}
(\byear{1973}).
\btitle{Random difference equations and renewal theory for products of random
  matrices}.
\bjournal{Acta Math.}
\bvolume{131}
\bpages{207--248}.
\bmrnumber{MR0440724}
\end{barticle}
\endbibitem

\bibitem{kksStable}
\begin{barticle}[msn]
\bauthor{\bsnm{Kesten},~\bfnm{H.}\binits{H.}},
  \bauthor{\bsnm{Kozlov},~\bfnm{M.~V.}\binits{M.~V.}} \AND
  \bauthor{\bsnm{Spitzer},~\bfnm{F.}\binits{F.}}
(\byear{1975}).
\btitle{A limit law for random walk in a random environment}.
\bjournal{Compositio Math.}
\bvolume{30}
\bpages{145--168}.
\bmrnumber{MR0380998}
\end{barticle}
\endbibitem

\bibitem{kGPD}
\begin{barticle}[msn]
\bauthor{\bsnm{Kobus},~\bfnm{Maria}\binits{M.}}
(\byear{1995}).
\btitle{Generalized {P}oisson distributions as limits of sums for arrays of
  dependent random vectors}.
\bjournal{J. Multivariate Anal.}
\bvolume{52}
\bpages{199--244}.
\bmrnumber{MR1323331}
\end{barticle}
\endbibitem

\bibitem{kmCLT}
\begin{barticle}[msn]
\bauthor{\bsnm{Kozlov},~\bfnm{S.~M.}\binits{S.~M.}} \AND
  \bauthor{\bsnm{Molchanov},~\bfnm{S.~A.}\binits{S.~A.}}
(\byear{1984}).
\btitle{Conditions for the applicability of the central limit theorem to random
  walks on a lattice}.
\bjournal{Dokl. Akad. Nauk SSSR}
\bvolume{278}
\bpages{531--534}.
\bmrnumber{MR764989}
\end{barticle}
\endbibitem

\bibitem{pThesis}
\begin{bmisc}[vtex]
\bauthor{\bsnm{Peterson},~\bfnm{J.}\binits{J.}}
(\byear{2008}).
Ph.D. thesis.
Awarded in 2008 by the University of Minnesota. Available at
\href{http://arxiv.org/abs/0810.257v1}{arXiv:0810.257v1 [math.PR]}.
\end{bmisc}
\endbibitem

\bibitem{rsBFD}
\begin{barticle}[msn]
\bauthor{\bsnm{Rassoul-Agha},~\bfnm{Firas}\binits{F.}} \AND
  \bauthor{\bsnm{Sepp{\"a}l{\"a}inen},~\bfnm{Timo}\binits{T.}}
(\byear{2006}).
\btitle{Ballistic random walk in a random environment with a forbidden
  direction}.
\bjournal{ALEA Lat. Am. J. Probab. Math. Stat.}
\bvolume{1}
\bpages{111--147 (electronic)}.
\bmrnumber{MR2235176}
\end{barticle}
\endbibitem

\bibitem{sRWRE}
\begin{barticle}[vtex]
\bauthor{\bsnm{Solomon},~\bfnm{Fred}\binits{F.}}
(\byear{1975}).
\btitle{Random walks in a random environment}.
\bjournal{Ann. Probab.}
\bvolume{3}
\bpages{1--31}.
\bmrnumber{MR0362503}
\end{barticle}
\endbibitem\vadjust{\goodbreak}

\bibitem{zRWRE}
\begin{bincollection}[msn]
\bauthor{\bsnm{Zeitouni},~\bfnm{Ofer}\binits{O.}}
(\byear{2004}).
\btitle{Random walks in random environment}.
In \bbooktitle{Lectures on Probability Theory and Statistics}.
\bseries{Lecture Notes in Math.}
\bvolume{1837}
\bpages{189--312}.
\bpublisher{Springer}, \baddress{Berlin}.
\bmrnumber{MR2071631}
\end{bincollection}
\endbibitem

\end{thebibliography}
\end{document}